\documentclass[11pt]{amsart}

\usepackage[utf8]{inputenc}

\usepackage{etoolbox}

\usepackage{mathtools}
\usepackage{color}
\usepackage{tikz}
\usepackage[makeroom]{cancel}
\usepackage{mathtools}
\usepackage[all,cmtip]{xy}

\usepackage{todonotes}

\usepackage{amsmath}
\usepackage{amsthm}

\usepackage{graphicx}
\graphicspath{ {figures/} }
\usepackage{amssymb}
\usepackage{epstopdf}
\usepackage{amsthm}

\usepackage{lingmacros}
\usepackage{tree-dvips}

\usepackage[numbers]{natbib}         
\usepackage[colorlinks]{hyperref}    

\DeclareMathOperator{\Span}{span}

\colorlet{shadecolor}{lightgray!25}

\newtheoremstyle{definitionsty}{3pt}{3pt}{\slshape}{}{\bfseries}{.}{.5em}{\addcontentsline{lod}{section}{#1~\protect\numberline{#2}{#3}}}
\theoremstyle{definitionsty}

\begin{document}

\newtheorem{thm}{Theorem}[section]

\theoremstyle{remark}
\newtheorem{theorem}{Theorem}[section]

\newtheorem{defi}[theorem]{Definition}
\newtheorem{cor}[theorem]{Corollary}
\newtheorem{conj}[theorem]{Conjecture}
\newtheorem{lemma}[theorem]{Lemma}
\newtheorem{remark}[theorem]{Remark}

\newcommand{\metasl}{\widetilde{SL_2}}

\newcommand{ \schnii }{ S \left (  V_{ \nu }^{ n } \right )  } 
\newcommand{ \jacquet }[1]{ #1_{ N, \psi_\beta } }
\newcommand{ \localquot }{ H_{ \nu }^{ \vv }  \backslash H_{ \nu } }
\newcommand{ \Ombeta }{ \Omega_{ \beta }  } 
\newcommand{ \coinv }[2]{ #1_{ H, \xi } }
\newcommand{ \anyCoinv }[1]{ #1_{ H , \xi } }

\newcommand{ \mmm }[2]{ m = #1, \mmp = #2 }
\newcommand{ \loc }{ _{ \nu }  }

\newtheorem{innercustomthm}{Theorem}
\newenvironment{customthm}[1]
  {\renewcommand\theinnercustomthm{#1}\innercustomthm}
  {\endinnercustomthm}

\newcommand{ \despegar }{ 
{\color{red} despegar - take off }
{\color{blue} despegar - take off }
{\color{magenta} despegar - take off }
{\color{red} despegar - take off }
{\color{blue} despegar - take off }
{\color{magenta} despegar - take off }
{\color{red} despegar - take off }
{\color{blue} despegar - take off }
{\color{magenta} despegar - take off }
{\color{red} despegar - take off }
{\color{blue} despegar - take off }
{\color{magenta} despegar - take off }
{\color{red} despegar - take off }
{\color{blue} despegar - take off }
{\color{magenta} despegar - take off }
 } 

\newcommand{ \om }{ \Omega }

\newcommand{ \ad }{ _{ \fA }  }
\newcommand{ \adpr }{  \left ( \fA \right )  } 
\newcommand{ \adata }{ \alpha =  \left (   \left ( q, V \right )  ,\xi
,
 \left ( q\p, V\p \right ),\xi\p
  \right ) } 
  \newcommand{ \adatanu }{ \alpha_{ \nu }  =  \left (   \left ( q_{ \nu } , V _{ \nu } \right )  ,\xi_{ \nu } 
,
 \left ( q\p_{ \nu } , V_{ \nu } \p \right ),\xi\p_{ \nu } 
  \right ) } 
\newcommand{ \conu }{ \cO_{ \nu }  } 

\newcommand{ \abcd }{ \begin{pmatrix} a & b \\ c & d \end{pmatrix} } 

\newcommand{ \snlam }{ \left ( SN \left ( h \right ) ,\lam \right ) _{ \Fnu } } 

\newcommand{ \localPi }{  \left ( \omega_{ \psivloc } \right ) _{ H_\nu, \xi_\nu }  } 
\newcommand{ \localPip }{  \left ( \omega_{ \psivlocp } \right ) _{ H\p_\nu, \xip_\nu }  } 

\newcommand{ \localFourier }{ \left ( 
\localPi
\right ) _{  N_{ \nu } , \psi_{ \beta }  } } 
\newcommand{ \localFourierp }{ \left ( 
\localPip
\right ) _{  N_{ \nu } , \psi_{ \beta }  } } 

\newcommand{ \Fn }{ F_{ \beta} }
\newcommand{ \qqquotin }{ H _{ \F }  \backslash H _{ \fA }   }

\newcommand{ \bktg }{ \aq G } 

\newcommand{ \seeapp }{ \proof
See Appendix \ref{SomeProofs} for the proof.

$ \square $ }

\newcommand{ \repeatNotations }{ Fix a maximal compact subgroup $ K \subseteq G\ad $  so that we have an Iwasawa decomposition 
$$
G\ad = P\ad K
 = M\ad N\ad K
$$
 and for 
$$
h = n\cdot m \left ( a \right ) k
,\qquad
n \in N\ad,a \in GL_{ n ,\fA} ,k \in K
$$
Let 
$$
\abs{  a \left ( g \right )  }  = 
\abs{  \det \left ( a \right )  } 
$$ 
and for
$$
\rhon = \dfrac{ n+1 } { 2 }
$$
let
$$
\Psi_{ s } \,\colon G\ad \to \fC,
\tn{ be given by }
\Psi_{ s } \left (  g \right )  = 
\abs{  a \left ( g \right )  } ^{ s + \rhon } 
$$ }

\newcommand{ \oop }{ \Omega \left (  \cO_{ \nu }  \right )  } 
\newcommand{ \vopp }{ V\p \left ( \cO_{ \nu }  \right )  } 
\newcommand{ \hop }{ H \left ( \cO_{ \nu }  \right )  } 
\newcommand{ \hopp }{ H\p \left ( \cO_{ \nu }  \right )  } 
\newcommand{ \pkop }[1]{  \qnu^{ #1 } \cO_{ \nu }  }
\newcommand{ \muc }[1]{ \, d \left (  #1 \right )  }
\newcommand{ \mucc }[1]{ \qnu^{ -#1 }   }
\newcommand{ \ps }{ \mucc s } 
\newcommand{ \pst }{ \mucc {  \left (  2s - 1 \right )  } } 
\newcommand{ \pk }[2]{ \qnu^{ -#1#2 }  }
\newcommand{ \ainv }{ a^{ -1 }  } 
\newcommand{ \gfac }{ \gamma^{ -1 }  } 
\newcommand{ \qnu }{ q_{ \nu } } 

\newcommand{ \phio }{ \phi_{ 0 } }
\newcommand{ \phiop }{ \phip_{ 0 } }

\newcommand{ \ovH }{ \overline{ H } } 
\newcommand{ \jackn }{ I_{ V } ^{ \phi } ( \xi ) }

\newcommand{ \rbeta }{  r \left ( \beta \right )  } 

\newcommand{ \betarepn }{ \beta\in Image \left ( \gram\cdot \right )  }
\newcommand{ \nicequotnv }{ H^{ \vv} _{ \fA }  \backslash H _{ \fA }   } 

\newcommand{ \weilAction }[3]{ \left (  \weilVn( #1, #2 ) \phi  \right ) ( #3 ) }
\newcommand{ \leviAction }{ \gamma^{ -1 } \abs{ \det( a ) } ( \det(a) , \left ( -1  \right )^{ \frac{ n(n-1) } { 2 } } ) \phi( ax )  } 
\newcommand{ \unipotentAction }{ \psi( \frac{ 1 } { 2 } tr( b\qv( x ) )   ) \phi( x ) } 
\newcommand{ \longWeylAction }{ \gamma^{ -1 } \int_{ V \otimes L } \phi( y ) \psi( b_{ V } ( x , y )  ) \,dy } 
\newcommand{ \Gt }{ \widetilde{ Sp }( W ) } 
\newcommand{ \SymQuot  }{ \Sym( \F ) \backslash \Sym( \fA )    } 
\newcommand{ \QuotF  }{ \F \backslash \fA  } 
\newcommand{ \alphat }{ \widetilde{ \alpha } }

\newcommand{ \glnf }{ GL_{ n }  \left (  \F \right )  } 

\newcommand{ \SymF }{ \Sym \left ( \F \right ) } 

\newcommand{ \N }{ N }

\newcommand{ \vox }{ V\otimes X } 
\newcommand{ \voxp }{ V\p \otimes X } 
\newcommand{ \voxpr }{  \left (  \vox \right )  } 
\newcommand{ \voxprp }{  \left (  \voxp \right )  } 
\newcommand{ \sqrClass }{ \F^{ \times }  /  \F^{ \times, 2 }  } 

\newcommand{ \sqrclass }{ \sqrClass } 
\newcommand{ \ovV }{ W } 
\newcommand{ \goodsec }{ f^{ \phi }_{ s, \choicew } } 
\newcommand{ \goodsecCH }[1]{ f^{ \phi }_{ s, #1 } } 
\newcommand{ \gooddist }{ Z_{ s, \choicew }  } 
\newcommand{ \gooddistCH }[1]{ Z_{ s, #1 }  } 
\newcommand{ \ovq }{ \overline{ q } } 
\newcommand{ \llam }{ \lam } 
\newcommand{ \choicew }{ w_{ 0 } } 
\newcommand{ \etaB }{ \eta } 
\newcommand{ \distH }{ Z \left (   \left ( h , h\p \right ) \cdot \phi \right ) =  \xi \left ( h \right )  \xi\p \left ( h\p \right )  Z \left ( \phi \right ) ,
\forall h\in H_{ \nu } , h\p\in H_{ \nu } \p,\phi\in  S \left ( W_\nu \right ) 
 }
\newcommand{ \distG }{ 
 Z \left ( k \cdot \phi \right ) = 
  Z \left ( \phi \right ) ,\forall \phi\in  S \left ( W_\nu \right ) ,
k \in K_{ \nu } 
 }
 \newcommand{ \pairH }{ Z \left (   \left ( h , h\p \right ) \cdot \phi \right ) =  \xi \left ( h \right )  \xi\p \left ( h\p \right )  Z \left ( \phi \right ) ,
\forall h\in H_{ \nu } , h\p\in H_{ \nu } \p,\phi\in  S \left ( W_\nu \right ) 
 }
\newcommand{ \pairG }{ 
 Z \left ( k \cdot \phi \right ) = 
  Z \left ( \phi \right ) ,\forall \phi\in  S \left ( W_\nu \right ) ,
k \in K_{ \nu } 
 }

\newcommand{ \mprime }{ m^\prime } 

\newcommand{ \cent }[1]{ \omega_{\sigma_{#1}} } 
\newcommand{ \centp }[1]{ \omega_{\sigma^{ \prime }_{#1}} } 
\newcommand{ \Fns }{ \tn{Fns} } 

\newcommand{ \ph }[1]{ \phi_{ #1 }  }
\newcommand{ \php }[1]{ \phi^{ \prime } _{ #1 }  } 

\newcommand{ \homj }{ Hom \left ( \F^{ n } ,V \right )  } 
\newcommand{ \homjp }{ Hom \left ( \F^{ n } ,V\p \right )  }


\newcommand{ \thliftii }[2]{ 
\fthLift{ #2 }
 }
\newcommand{ \ovvpi }{ \overline{ \varphi_\pi } } 
\newcommand{ \thliftpii }[2]{ I^{ #2  }_{ V^{ \prime } _{ #1 }, \xi\p_{ #1 } } }

\newcommand{ \thlift }[1]{ \thliftii{ #1 }{ \ph{ #1  } } }
\newcommand{ \thliftp }[1]{ \thliftpii{ #1 }{ \php{ #1 } } }


 \newcommand{ \homogii }[2]{ 
#1  \left (  \F \right )  \backslash #2 \left (  \fA \right )
 } 

\newcommand{ \homog }[1]{ 
\homogii{ #1 }{ #1 }
 } 
 
 \newcommand{ \doubleTri }{ \left < \left < \cdot, \cdot \right >  \right >  } 

\newcommand{ \slint }{ \homog{ SL_{ 2 } } } 
\newcommand{ \slinttilde }{ \homog{ \widetilde{ SL_{ 2 } }  } } 

\newcommand{ \spint }[1]{ \homog{ Sp \left ( W \right )  } } 
\newcommand{ \spinttilde }[1]{ \homog{ \widetilde{ \spint }  } }

\newcommand{ \nquad }{ \mkern-72mu } 
\newcommand{ \rhon }{ \rho_{ n }  }

\newcommand{ \gpa  }{ \qAA P G }

\newcommand{ \finalGlobalExpansionAux }{ \int_{ K } \int_{ \ga M }  W_{ \phi }  \left (  mk \right )  W\p_{ \phip }  \left (  mk \right ) \delb m \fPsi{ mk }{ s } \,dm \,dk } 

\newcommand{ \finalGlobalExpansion }[1]{ \int_{ K } \int_{ \ga M }  W_{ \phi }  \left (  mk \right )  W\p_{ \phip }  \left (  mk \right ) \delb m \fPsi{ mk }{ s = #1 } \,dm \,dk } 

\newcommand{ \theUnramifiedSetting }{  \begin{enumerate}
\item
 $ \vv \in \voxpr \left (  \cO_{ \nu } \right ),\vv \p\in \voxprp \left (  \cO_{ \nu } \right ) $
\item
Let  $ \phio\in  S \left ( \voxpr \left ( \cO_\nu \right ) \right ) ,\phiop\in  S \left ( \voxprp \left ( \cO_\nu \right ) \right )  $ be two Schwartz functions defined on a lattice, i.e. given by characteristic functions
 $$
\phio
= 
1_{ \voxpr \left (  \cO_{ \nu } \right )  } 
\tn{ and }
\phiop
= 
1_{ \voxprp \left (  \cO_{ \nu } \right )  } 
$$
\item
 The discriminants $ d,d\p $ of our quadratic spaces $ \left ( q, V \right ) ,\left ( q\p, V\p \right ) $, respectively, satisfy
 $ d,d\p\in \cO_{ \nu } ^{ \times }  $.
\item
$ W_{ \phi } \,\colon G_\nu \to \fC $ and $ W\p_{ \phip } \,\colon G_\nu \to \fC $ are $ K $-finite where $ K \subseteq G $ is a maximal compact subgroup.
\item
$ W_{ \phi } \left (  e \right ) = 1 $ and $ W\p_{ \phip } \left (  e \right ) = 1 $
where $ e $ is the identity element of $ G $.
\end{enumerate} }

\newcommand{ \glblZetaConsti }[1]{ 
b_{ n } ^{ S }  \left ( #1 \right ) 
  } 

\newcommand{ \glblZetaConst }{ 
b_{ n } ^{ S }  \left (  s \right ) 
  } 

\newcommand{ \glblZetaConstSub }{ 
b_{ n } ^{ S }  \left ( \rhon \right ) 
  }   

\newcommand{ \sxa }{ S( X_{ \fA }  ) }
\newcommand{ \sva }{ S( V_{ \fA } ) }
\newcommand{ \svoa }{ S( V_{ \fA } \p ) }
\newcommand{ \svap }{ S( V_{ \fA } \p ) }
\newcommand{ \schi }{ S( V\p_{ \nu } )  }
\newcommand{ \schii }{ S( V_{ \nu } )  }
\newcommand{ \schn }{ S( V_{ n,\nu } )  }

\newcommand{\SymVar}[1]{\rm Sym_{ #1 } }
\newcommand{\Sym}{\rm Sym_n}
\newcommand{\SymThree}{\rm Sym_3}
\newcommand{\SymQ}{\Sym(\F)}
\newcommand{\SymA}{\Sym(\fA)}
\newcommand{\SymAQ}{\SymQ \backslash \SymA}

\newcommand{ \qqquoti }{ H_{ 0 } _{ \F }  \backslash H_{ 0 } _{ \fA }   }
\newcommand{ \qqquotii }{ \F \backslash \fA }
\newcommand{ \psiVone }{ \psi,V_1 }
\newcommand{ \thetainti }{ I_{ V_3 } ^{ \phi_{ 3 }  } ( \xi_{V_3} )  }
\newcommand{ \thetaintii }{ I_{ V_1 } ^{ \phi_{ 1 }  } ( \xi_{V_1} )  }
\newcommand{ \bp }{ b }

\newcommand{ \vvbeta }{ \vv^{\beta} }
\newcommand{ \vvbetanu }{ \vv^{\beta_{ \nu } } }

\newcommand{ \vvbetao }{ \vv^{\betao} }
\newcommand{ \vvbetaonu }{ \vv^{\betao_{ \nu } } }

\newcommand{ \sumi }{ \sum_{ \gamma\in V_1(\F) } }
\newcommand{ \sumii }{ \sum_{ \substack { \gamma\in V_1(\F) \\ q_{ V_1 } (\gamma)=\beta }} }

\newcommand{ \nsumi }[1]{ \sum_{ \gamma\in V_{ #1 } (\F) } }
\newcommand{ \nsumii }[1]{ \sum_{ \substack { \gamma\in V_1(\F) \\ q_{ V_{ #1 }  } (\gamma)=\beta }} }
\newcommand{ \nmsum }{ \sum_{ \gamma\in \vox_{ \F } } }

\newcommand{ \betarep }{ 0 \neq \beta\in \text{image}(\qii) }
\newcommand{ \jack }{ I_{ V_1 } ^{ \phi } ( \xi_{V_1} ) }

\newcommand{ \outerfactorii }{ \gamma \absv{\alpha}^{-\frac{ 1 } { 2 }} ( \alpha , -d(V_1) ) }
\newcommand{ \Fone }{ F_{1,\beta} }
\newcommand{ \Fonebetazero }{ F_{1,\beta=0} }
\newcommand{ \Fthree }{ F_{3,\beta} }
\newcommand{ \fourierbeta }{ \mathcal F r_{\beta} }

\newcommand{ \pmset }{ \{\pm1\} }
\newcommand{ \vinFmodpm }{ v \in \pmset \backslash V_{ 1 } (\F) - \{0\}  }

\newcommand{ \nicequot }{ H_{ 0 }^{ \vvbeta } _{ \fA }  \backslash H_{ 0 } _{ \fA }   } 
\newcommand{ \nicequotnu }{ H_{ 0 }^{ \vvbeta } _{ \F }  \backslash H_{ 0 } _{ \nu }   } 
\newcommand{ \badquot }{ H_{ 0 }^{ \vvbeta } _{ \F }  \backslash H^{ \vvbeta }_{ 0 } _{ \fA }   } 

\newcommand{ \hr }{ {\color{red} HERE } } 

\newcommand{ \contra }[1]{ \hat{ #1 } } 

\newcommand{ \nnsumi }[1]{ \sum_{ \gamma\in \voxpr\left (  \F  \right ) } }
\newcommand{ \nnsumiF }[1]{ \sum_{ \gamma\in \voxpr_{ \F } } }
\newcommand{ \nnsumiFbeta }[1]{ \sum_{ \substack{
\gamma\in \voxpr_{ \F } 
\\
\gram \gamma  = \beta
} } }

\newcommand{ \qcdata }{ 
\eta = \left ( \left ( q, V \right ), \xi \right ) 
 } 
 
 \newcommand{ \goodqcData }{ GoodData } 
 
\newcommand{ \qcdatap }{ 
\eta\p = \left ( \left ( q\p, V\p \right ), \xi\p \right ) 
 } 

\newcommand{ \qqquotni }{ H _{ \F }  \backslash H _{ \fA }   }
\newcommand{ \qqquotnii }{ \Sym \F \backslash \Sym \fA }
\newcommand{ \qqquotniii }{ H ^{ \vvbeta } _{ \F } \backslash H _{ \F }  }
\newcommand{ \qqquotniiii }{ H ^{ \vvbeta } _{ \F } \backslash H _{ \fA }  }
\newcommand{ \qqquotniiiii }{ H ^{ \vvbeta } _{ \F } \backslash H _{ \fA }  }

\newcommand{ \nicequotn }{ H^{ \vvbeta } _{ \fA }  \backslash H _{ \fA }   } 
\newcommand{ \nicequotm }{ H^{ \vvbeta } _{ \fA }  \backslash H _{ \fA }   } 
\newcommand{ \nicequotOm }{ H^{ \vvbetao } _{ \fA }  \backslash H _{ \fA }   } 

\newcommand{ \nicequotOnnu }{ H^{ \vvbetao } _{ \nu }  \backslash H _{ \nu }   } 

\newcommand{ \nicequotnnu }{ H^{ \vvbeta } _{ \nu }  \backslash H _{ \nu }   } 
\newcommand{ \badquotn }{ H^{ \vvbeta } _{ \F }  \backslash H^{ \vvbeta } _{ \fA }   } 
\newcommand{ \badquotm }{ H^{ \vvbeta } _{ \F }  \backslash H^{ \vvbeta }_{ \fA }   } 
\newcommand{ \badquotOm }{ H^{ \vvbetao } _{ \F }  \backslash H^{ \vvbetao }_{ \fA }   } 

\newcommand{ \badquotOn }{ H^{ \vvbetao } _{ \F }  \backslash H^{ \vvbetao } _{ \fA }   } 

\newcommand{ \qqquotF }{ \F \backslash \fA }

\newcommand{ \thetaintni }{ I_{ V,\xi } ^{ \phi } }

\newcommand{\wni}{\omega_{ \vpsi } }

\newcommand{ \prodnu }{ \prod_{ \nu } }
\newcommand{ \prodnuIn }[1]{ \prod_{ \nu \in #1 } }
\newcommand{ \prodnuNotIn }[1]{ \prod_{ \nu \not \in #1} }

\newcommand{ \cppf }{ \cpffin }
\newcommand{ \pf }{ \cpf \cap S } 
\newcommand{ \si }{ \cpfInf \cap S }
\newcommand{ \sii }{ \cppf \cap S }

\newcommand{ \cpf }{ \Sigma_{ \F }  } 
\newcommand{ \cpfInf }{ \Sigma_{ \F, \infty } }
\newcommand{ \cpffin }{ \Sigma_{ \F, f } }

\newcommand\restr[2]{\ensuremath{\left.#1\right|_{#2}}}
\newcommand{ \qi }{ q_{ V_3  }   }
\newcommand{ \qii }{ q_{V_1} }

\newcommand{ \qinu }{ q_{ V_{ 3, \nu }  }  }
\newcommand{ \qiinu }{ q_{V_{ 1, \nu }  } }

\newcommand{ \xiv }{ \xi_{ V_{ 3 }  }  }
\newcommand{ \xivo }{ \xi_{V_{ 1 }} }
\newcommand{ \vv }{ \mathbf{ \overline{ v } } }
\newcommand{ \ww }{ \mathbf{ \overline{ w } } }
\newcommand{ \uu }{ \mathbf{ \overline{ v }_{ 0 } } }

\newcommand{\wi}{\omega_{\psi ,V }}
\newcommand{\wii}{\omega_{\psi ,V\p }}

\newcommand{\weilVn}{\omega_{ \vpsi }}
\newcommand{\wiii}{\omega_{\psi ,V}(n(b)n(b')g}

\newcommand{\HALF}{\frac{1}{2}}
\newcommand{\tr}{\rm tr}
\newcommand{\nuI}{\nu_i}

\newcommand{ \Fx }{ \F^{ \times }  } 

\newcommand{\Fnu}{\mathbb{F}_{\nu}}
\newcommand{\Fnux}{\mathbb{F}_{\nu}^{ \times } }
\newcommand{\FnuxSqr}{\mathbb{F}_{\nu}^{ \times, 2 } }

\newcommand{\cprod}[2]{#1 \times #2}
\newcommand{\pcprod}[2]{(\cprod{#1}{#2})}


\newcommand{ \good }{ {\color{blue} } } 

\newcommand{ \qvformleft }{ ( \qi,V_{ 3 }  )  } 
\newcommand{ \chil }{  \alpha _{ \lambda,\nu }  } 

\newcommand{ \xik }{ \xi_{ k,\nu } }
\newcommand{ \xil }{ \xi_{ \lambda,\nu } }
\newcommand{ \xileps }{ \xi_{ \lambda, \eps, \nu } }
\newcommand{ \xilamepsnu }{ \xi_{ \lambda, \eps,\nu } }

\newcommand{ \refv }{ \refl{v } }
\newcommand{ \refu }{ \refl{u } }
\newcommand{ \refup }{ \refl{\up } }
\newcommand{ \refui }{ \refl{u_{ 1 } } } 
\newcommand{ \refuii }{ \refl{u_{ 2 } } }
\newcommand{ \refvv }{ \refl{ \vvbeta } } 
\newcommand{ \refw }{ \refl{ w } } 

\newcommand{ \up }{ u^{ \prime }  }
\newcommand{ \bi }{ b_{V_3} }

\newcommand{ \subgp }[1]{ #1 \leq G } 
\newcommand{ \provedis }{ \tn{Prove/disprove: } } 

\newcommand{ \aV }{ V } 
\newcommand{ \squarex }{ \square^{ \times }  } 

\newcommand{ \qqqiii }{ q( u_1 ) q( u_2 )  } 
\newcommand{ \bbbiii }{ b{ ^2 } ( u_1 , u_2 )    }

\newcommand{ \rfield }{ \mathcal O_{ \nu } / \mathcal P_{ \nu }  } 		
\newcommand{ \squareclass }{ \Fnux / \left (  \Fnux  \right )^{ 2 }  } 
\newcommand{ \resq }{ \overline{ q_{ 2 }  } } 

\newcommand{ \mmupi }{  -\mu\cdot\pi } 	
\newcommand{ \mupi }{ \mu\cdot\pi } 	
\newcommand{ \blue }[1]{ \color{blue} #1 \color{black} } 
\newcommand{ \red }[1]{ \color{red} #1 \color{black} } 
\newcommand{ \green }[1]{ \color{green} #1 \color{black} } 
\newcommand{ \yellow }[1]{ \color{yellow} #1 \color{black} } 
\newcommand{ \magenta }[1]{ \color{magenta} #1 \color{black} } 
\newcommand{ \QF }[3]{ <#1,#2,#3> } 
\newcommand{ \hbnu }[2]{ ( #1  , #2 )_{ \nu } } 
\newcommand{ \hbfA }[2]{ ( #1  , #2 )_{ \fA } } 
\newcommand{ \qu }{ q_{ U }  } 
\newcommand{ \rii }{  \fR^{ 2 } } 
\newcommand{ \pto }{ \left (  u_{ 0 } , v_{ 0 }   \right ) } 

\newcommand{ \viiinu }{ V_{ 3, \nu } } 
\newcommand{ \vinu }{ V_{ 1, \nu } } 

\newcommand{ \xileft }{ \xi } 
\newcommand{ \xiright }{ \xi\p } 

\newcommand{ \xileftnu }{ \xi_\nu } 
\newcommand{ \xirightnu }{ \xi\p_\nu } 

\newcommand{ \localquad }{ \left (  q_{ \viiinu  } , v  , \xileftnu, \xirightnu  \right ) } 
\newcommand{ \localquadvv }{ \left (  q_{ \viiinu  } , \vv  , \xileftnu, \xirightnu  \right ) } 

\newcommand{\nb}{\begin{pmatrix} I_{ n } & b \\  & I_{ n } \end{pmatrix}}

\newcommand{ \spfour }{ \widetilde{ Sp }_{ 4 }  } 
\newcommand{ \spsix }{ \widetilde{ Sp }_{ 6 }  } 
\newcommand{ \spen }{ G_{ \fA } }
\newcommand{ \ofive }{ O( q_{ 5 }  ) } 
\newcommand{ \oem }{  H_{ \fA } } 
\newcommand{ \diag }{ \tn{diag} } 
\newcommand{ \rank }{ \tn{rank} } 
\newcommand{ \conf }[1]{}
\newcommand{ \tooDetail }{}

\newcommand{ \schad }{ S( \voxpr \ad ) }
\newcommand{ \schadp }{ S( \voxprp \ad ) }
\newcommand{ \schloc }{ S( \voxpr_{ \nu } ) }
\newcommand{ \schlocp }{ S( \voxprp_{ \nu } ) }

\newcommand{ \auxRange }{ \rhon }

 
 \newcommand{\abs}[1]{\left| {#1} \right|}

\newcommand{\cA}{ \mathcal {A}}
\newcommand{\cB}{ \mathcal {B}}
\newcommand{\cC}{ \mathcal {C}}
\newcommand{\cD}{ \mathcal {D}}
\newcommand{\cE}{ \mathcal {E}}
\newcommand{\cF}{ \mathcal {F}}
\newcommand{\cG}{ \mathcal {G}}
\newcommand{\cH}{ \mathcal {H}}
\newcommand{\cI}{ \mathcal {I}}
\newcommand{\cJ}{ \mathcal {J}}
\newcommand{\cK}{ \mathcal {K}}
\newcommand{\cL}{ \mathcal {L}}
\newcommand{\cM}{ \mathcal {M}}
\newcommand{\cN}{ \mathcal {N}}
\newcommand{\cO}{ \mathcal {O}}
\newcommand{\cP}{ \mathcal {P}}
\newcommand{\cQ}{ \mathcal {Q}}
\newcommand{\cR}{ \mathcal {R}}
\newcommand{\cS}{ \mathcal {S}}
\newcommand{\cT}{ \mathcal {T}}
\newcommand{\cU}{ \mathcal {U}}
\newcommand{\cV}{ \mathcal {V}}
\newcommand{\cW}{ \mathcal {W}}
\newcommand{\cX}{ \mathcal {X}}
\newcommand{\cY}{ \mathcal {Y}}
\newcommand{\cZ}{ \mathcal {Z}}

\newcommand{\fA}{\mathbb{A}}
\newcommand{\fB}{\mathbb{B}}
\newcommand{\fC}{\mathbb{C}}
\newcommand{\fD}{\mathbb{D}}
\newcommand{\fE}{\mathbb{E}}
\newcommand{\fF}{\mathbb{F}}
\newcommand{\fG}{\mathbb{G}}
\newcommand{\fH}{\mathbb{H}}
\newcommand{\fI}{\mathbb{I}}
\newcommand{\fJ}{\mathbb{J}}
\newcommand{\fK}{\mathbb{K}}
\newcommand{\fL}{\mathbb{L}}
\newcommand{\fM}{\mathbb{M}}
\newcommand{\fN}{\mathbb{N}}
\newcommand{\fO}{\mathbb{O}}
\newcommand{\fP}{\mathbb{P}}
\newcommand{\fQ}{\mathbb{Q}}
\newcommand{\fR}{\mathbb{R}}
\newcommand{\fS}{\mathbb{S}}
\newcommand{\fT}{\mathbb{T}}
\newcommand{\fU}{\mathbb{U}}
\newcommand{\fV}{\mathbb{V}}
\newcommand{\fW}{\mathbb{W}}
\newcommand{\fX}{\mathbb{X}}
\newcommand{\fY}{\mathbb{Y}}
\newcommand{\fZ}{\mathbb{Z}}


\newcommand{\pfA}{ \left ( \fA \right ) }
\newcommand{\pfB}{ \left ( \fB \right ) }
\newcommand{\pfC}{ \left ( \fC \right ) }
\newcommand{\pfD}{ \left ( \fD \right ) }
\newcommand{\pfE}{ \left ( \fE \right ) }
\newcommand{\pfF}{ \left ( \fF \right ) }
\newcommand{\pfG}{ \left ( \fG \right ) }
\newcommand{\pfH}{ \left ( \fH \right ) }
\newcommand{\pfI}{ \left ( \fI \right ) }
\newcommand{\pfJ}{ \left ( \fJ \right ) }
\newcommand{\pfK}{ \left ( \fK \right ) }
\newcommand{\pfL}{ \left ( \fL \right ) }
\newcommand{\pfM}{ \left ( \fM \right ) }
\newcommand{\pfN}{ \left ( \fN \right ) }
\newcommand{\pfO}{ \left ( \fO \right ) }
\newcommand{\pfP}{ \left ( \fP \right ) }
\newcommand{\pfQ}{ \left ( \fQ \right ) }
\newcommand{\pfR}{ \left ( \fR \right ) }
\newcommand{\pfS}{ \left ( \fS \right ) }
\newcommand{\pfT}{ \left ( \fT \right ) }
\newcommand{\pfU}{ \left ( \fU \right ) }
\newcommand{\pfV}{ \left ( \fV \right ) }
\newcommand{\pfW}{ \left ( \fW \right ) }
\newcommand{\pfX}{ \left ( \fX \right ) }
\newcommand{\pfY}{ \left ( \fY \right ) }
\newcommand{\pfZ}{ \left ( \fZ \right ) }

\newcommand{ \bin }[1]{ \frac{ #1 \left (  #1 - 1 \right )  } { 2 } }

\newcommand{\tn}{\textnormal}

 \newcommand*\from{\colon}

\newcommand{\set}[2]{ \{ #1 | #2 \}}

\newcommand{\ra}{\rightarrow}
\newcommand{\raabove}[1]{\overset{#1}\ra}
\newcommand{\rabelow}[1]{\underset{#1}\ra}

\newcommand{\resf}[2]{\left.#1\right|_{#2}}

\newcommand{\defeq}{\vcentcolon=}
\newcommand{\eqdef}{=\vcentcolon}

\newcommand{\mativ}[4]{\left[
 \begin{array}{ c c }
     #1 & #2 \\
     #3 & #4
  \end{array} \right]
}

\newcommand{ \pip }{  \pi^{ \prime }  } 

\newcommand{\st}[1]{\{ #1 \} }
\newcommand{\limply}[2]{#2 \Leftarrow #1}
\newcommand{\lequiv}[2]{#2 \Leftrightarrow #1}

\providecommand\given{} 
\newcommand\SetSymbol[1][]{\nonscript\,#1\vert \allowbreak \nonscript\,\mathopen{}}
\DeclarePairedDelimiterX\Set[2]{\lbrace}{\rbrace}%
 { #1 \,\delimsize| \,\mathopen{} #2 }

 \newcommand{ \gboxi }[2]{ \parbox{#1em}{#2} }
\newcommand{ \gbox }[1]{ \gboxi{15}{#1} }
\newcommand{ \xipnu }[1]{ \xi^{ \prime } _{ \nu }  } 
\newcommand{ \xinu }[1]{ \xi_{ \nu }  }

\newcommand{ \co }[2]{ \Theta_{ V_{ \nu }, W_{ \nu }, #1 }  \left (  \xi_{ \nu }  \right ) }
\newcommand{ \cop }[2]{ \Theta_{ V\p_{ \nu }, W_{ \nu }, #1 }  \left (  \xip_{ \nu } \right ) }

\newcommand{ \gco }[2]{ \Theta_{ #1 }  \left( #2 \right) }
\newcommand{ \gcop }[2]{ \Theta_{ #1 }  \left( #2 \right) }
\newcommand{ \codef }[1]{ Co_{ \psi_{ \nu }  }  \left (  #1 \right ) \defeq \cocoa{ V_{ #1 }  }{ \psi_{ \nu }  }{ O \left (  V_{ #1 }  \right )  }{ \xi_{ #1 }  } } 
\newcommand{ \codefp }[1]{ Co_{ \psi_{ \nu }  }^{ \prime }   \left (  #1 \right ) \defeq \cocoa{ \vp_{ #1 }  }{ \psi_{ \nu }  }{ O \left (  \vp_{ #1 }  \right )  }{ \xip_{ #1 }  } }
 \newcommand{ \kh }[2]{  \left (  #1, #2 \right ) _{ \F_{ \nu }  }  }

\newcommand*\squared[1]{\tikz[baseline=(char.base)]{
            \node[shape=rectangle,draw,inner sep=2pt] (char) {#1};}}


\newcommand{ \ai }{ a_1 } 
\newcommand{ \aii }{ a_2 } 
\newcommand{ \aiii }{ a_3 }

\newcommand{ \zero }{ \mathbf{ \underline 0 } }

\newcommand{ \thetap }{ \theta^{ \prime }  } 
\newcommand{ \tempH }{ H \left ( \F \right ) \backslash H \left ( \fA \right )  } 
\newcommand{ \tempHp }{ \hp \left ( \F \right ) \backslash \hp \left ( \fA \right )  } 

\newcommand{ \prodhh }{ H\times \hp } 
\newcommand{ \hhp }{  \left ( \prodhh \right )  } 
\newcommand{ \ovh }{  \overline{ H } } 
\newcommand{ \bigq }{ q_{ W }  }

\newcommand{ \ip }[1]{  \left (  #1 \right )  }
\newcommand{ \extendedWeil }[3]{ \omega_{ #3, #1 } }
\newcommand{ \weil }[2]{ \extendedWeil{ #1 }{ #2 }{ \psi } }
\newcommand{ \sextendedWeil }[2]{ \omega_{ #2, #1 } }
\newcommand{ \sweil }[1]{ \extendedWeil{ #1 }{ \psi } }
\newcommand{ \ma }{ \begin{pmatrix} a &  \\  & a^\star \end{pmatrix} } 
\newcommand{ \tgi }{ \widetilde{ { G }_{ 1 }  } } 
\newcommand{ \tgii }{ \widetilde{ { G }_{ 2 }  } } 
\newcommand{ \inth }[1]{ \int_{ \homog{ H_{ #1 }  } }  } 
\newcommand{ \intg }{ \int_{ \homog{ G_{ n }  } }  } 
\newcommand{ \inthh }{ \int_{ \homog{ H } }  } 
\newcommand{ \inthhp }{ \int_{ \homog{ \hp } }  } 
\newcommand{ \ccl }[1]{ \co{ \psi }{ #1 }  }
\newcommand{ \ccr }[1]{ \cop{ \psi }{ #1 } }
\newcommand{ \ccla }[1]{ \ccl #1  \left (  \fA \right ) }
\newcommand{ \ccra }[1]{ \ccr #1  \left (  \fA \right ) }

\newcommand{ \eps }{ \epsilon } 
\newcommand{ \epsp }{ \epsilon^{ \prime }  } 
\newcommand{ \lam }{ \lambda } 
\newcommand{ \lamp }{ \lambda^{ \prime }  } 
\newcommand{ \inj }{ \hookrightarrow }

\newcommand{ \clred }[1]{ \color{red} \cancel{ #1 } \color{black} }
\newcommand{ \clbluered }[1]{ \color{blue} \cancel{ #1 } \color{red} }

\newcommand*\circled[1]{\tikz[baseline=(char.base)]{
            \node[shape=circle,draw,inner sep=2pt] (char) {#1};}}

\newcommand{\LDOTSvec}[1]{(#1_1,\ldots,#1_n)}

\newcommand{\resfunc}{\rm F}

\newcommand{ \horizline }[2]{ \noindent\rule{#1cm}{#2pt} }

\newcommand{ \phip }{ \phi^{ \prime }  } 
\newcommand{ \omom }[2]{ \omega_{ { #1 }, #2 }  } 
\newcommand{ \cocoa }[4]{  \left (  \omom{ #1 }{ #2 } \right )_{ #3, #4 }  }

\newcommand{ \orb }[1]{ \cO_{ #1 }  }

\newcommand{ \bb }{ \widetilde{ B } } 

\newcommand{ \slsl }{ \widetilde{ SL }_{ 2 }  } 

\newcommand{ \arithFunc }{ \tn{Fns}( \widetilde{ sl_2 }_{ \F }  \backslash \widetilde{ sl_2 }_{ \fA } )  }


\newcommand{ \VV }{ V_{ n } ^{ \prime } } 
\newcommand{ \qq }{ q^{ \prime }  } 
\newcommand{ \chch }{ \xi_{ n + 2 } ^{ \prime }  } 
\newcommand{ \F }{ \fF } 
\newcommand{ \K }{ \fK } 

\newcommand{ \mmp }{ m^{ \prime }  } 
\newcommand{ \sigmap }{ \sigma^{ \prime }  } 
\newcommand{ \p }{ ^{\prime} } 


\newcommand{ \thepair }{ H \times G } 
\newcommand{ \gammai }[1]{ \gamma \left (  #1, \psi^{ \frac{ 1 } { 2 } } 
 \right ) ^{ -1 }  }

\newcommand{ \Vp }{ V^{ \prime }_{ n }   } 
\newcommand{ \vi }{ v_{ 1 }  } 
\newcommand{ \vn }{ v_{ n }  } 



\newcommand{ \spnnn }{ \widetilde{ Sp }_{ 2n }  } 
\newcommand{ \vnnn }{ V_{ n }  } 
\newcommand{ \qnnn }{ q_{ n }  } 
\newcommand{ \vmmm }{ V } 
\newcommand{ \qmmm }{ q } 
\newcommand{ \xinnn }{ \xi_{ n }  } 
\newcommand{ \ximmm }{ \xi } 
\newcommand{ \ximmmnu }{ \xi_{ m, \nu }  } 
\newcommand{ \hnnnA }{ H_{ n } _{ \fA }   } 
\newcommand{ \hnnnF }{ H_{ n } _{ \F }   } 
\newcommand{ \hnnnFnu }{ H_{ 5 } ( \Fnu )  } 
\newcommand{ \hmmm }{ H } 
\newcommand{ \hmmmA }{ H_{ \fA }  } 
\newcommand{ \hmmmAp }{ H\p_{ \fA }  } 
\newcommand{ \hmmmF }{ H_{ \nu } } 
\newcommand{ \hmmmFp }{ H\p_{ \nu } } 
\newcommand{ \hmmmFvvp }{ H^{ \prime, \vv } _{ \nu } } 

\newcommand{ \summmm }{ \sum_{ \substack { \gamma\in \vmmm_{ \F }  \\ \gramqm{ \gamma } =\beta }} }

\newcommand{ \ranklinen }[5]{ \rank \beta = #1 \implies \vv = ( #2 , #3, #4, \ldots, #5 ) \in \left (  \vmmm \otimes X \right )_{ \F } \cong \vmmm^{ n }_{ \F } }

\newcommand{ \Ttilde }{ \widetilde{ T } } 
\newcommand{ \Ztilde }{ \widetilde{ Z } } 
\newcommand{ \SLtwonu }{ \SLtwo,_{ \nu }  } 
\newcommand{ \irrSLtwo }{ Irr_{ gen } { \left (   \SLtwonu  \right ) } } 
\newcommand{ \presentpair }{ \ofive \times \spsix } 
\newcommand{ \dualpair }[2]{ { O( q_{ #1 }  ) } \times \widetilde{ Sp }_{ #2 }  }

\newcommand{ \HHHfive }{ H_{ 5 }  }

\newcommand{ \booyaquotfive }{ \HHHfive^{ \vv }_{ \fA }  \backslash \HHHfive_{ \fA }   } 

\newcommand{ \intformulaPARfive }[3]{  \int_{ \booyaquotfive } #1(h^{ -1 } #2) \xi_{ V_3 } ( h ) \, dh  } 
\newcommand{ \formulaPARfive }[3]{ \absv { #3 } ( #3 , -d(q_{ 5 } ) d(q_{ 3 } ) )_{ \nu }  \intformulaPARfive {#1}{#2}{#3} }

\newcommand{ \localquadfive }{ \left (  q_{ _{ 5 }   } , \vv  , \xi_{ 5 } , \eta_{ 3 }   \right ) } 
\newcommand{ \localquadfivevv }{ \left (  q_{ \viiinu  } , \vv  , \xileftnu, \xirightnu  \right ) }

\newcommand{ \schfivei }{ S( V_{ 3,\nu } ^{ 3 } ) } 
\newcommand{ \schfiveii }{ S( V_{ 5,\nu } ^{ 3 } ) } 

\newcommand{ \hfivenF }{  H_{ n+2 } \left (  \fF  \right ) }
\newcommand{ \hfivenA }{  H_{ n+2 } \left (  \fA  \right ) } 

\newcommand{ \statei }{ \ensuremath{\int_{ H^{\vv}_{ \fA }  \backslash H_{ \fA } }\left (  \wi( n(b)g,h ) \phi  \right )(\vv)\xiv( h ) \,dh } }
\newcommand{ \uptosqr }{ (\F_{ \nu }^{ \times })^{ 2 }\backslash\F_{ \nu } }
\newcommand{ \iotav }{ \iota_{ X,V } }
\newcommand{ \iotavo }{ \iota_{ X,V_{ 0 } } }
\newcommand{ \tilb }{ \tilde{B} }

\newcommand{ \rofi }{ \mathcal{O} }
\newcommand{ \locrofi }{ \rofi_{ \nu } }
\newcommand{ \VOO }{ V(\locrofi) }
\newcommand{ \locQuot }{ H^{\vv}_{\nu} \backslash H_{\nu} }
\newcommand{ \volume }{ Vol_{ \locQuot }( \VOO )  }
\newcommand{ \hgood }{ \tilde h }
\newcommand{ \nicelimit }[4]{ #1\lim_{ #2 \to #3 } #4 }
\newcommand{ \oneOverFraction }[1]{ \frac{ 1 } { #1 } }

\newcommand{ \mi }{ \begin{pmatrix} 0 & 1 \\ 0 & 0 \end{pmatrix} }
\newcommand{ \mii }{ \begin{pmatrix} 1 & t \\ 0 & 1 \end{pmatrix} }
\newcommand{ \miii }{ \begin{pmatrix} 0 & t \\ 0 & 0 \end{pmatrix} }
\newcommand{ \expi }{ e^{ 2\pi i t x^2 } }

\newcommand{ \Whit }[1]{ W_{ \psi }( #1,\xi_{ #1,\nu } )  }
\newcommand{ \lxiv }{ \xi_{ V,\nu } }
\newcommand{ \lxivo }{ \xi_{ V_{ 0 },\nu }}

\newcommand{ \induced }[1]{ Ind_{ B }^{ \widetilde{SL_2} }( #1 )  }

\newcommand{ \resprod }{ \otimes^{ \prime }_{ \nu } }

\newcommand{ \pinu }{  \ensuremath{\pi_{ \nu } }}
\newcommand{ \pionu }{ \pi_{ 0,\nu } }

\newcommand{ \homi }[1]{ Hom_{ G_{ \nu }\times H_{ \nu } }(\wi,#1\otimes1_{ H_{\nu} }) }
\newcommand{ \homii }[1]{ Hom_{ G_{ \nu }\times H_{ 0,\nu } }(\wii,#1\otimes1_{ H_{0,\nu} }) }
\newcommand{ \vquoto }{ V_{0,j-1} \backslash V_{ 0,j } }
\newcommand{ \vquot }{ V_{j-1} \backslash V_{ j } }

\newcommand{ \WI }[1]{ \wi( #1 ) }
\newcommand{ \WII }[1]{ \wii( #1) }

\newcommand{ \hh }{ H_{ \F }  }
\newcommand{ \hho }{ H^{ 0 }_{ \F }  }
\newcommand{ \qf }{ q_{V} }
\newcommand{ \omg }{  \Omega_{\kappa}^{\times} }
\newcommand{ \thetaIntegral }{ I_V^{ \deltaA( \phi )  }( \xi_V ) }
\newcommand{ \frthetaIntegral }{ I_{V ,\beta}^{ \deltaA( \phi )  }( \xi_V ) }
\newcommand{ \deltaA }{ \Delta }
\newcommand{ \ivdelphi }{ I_V^{ \deltaA(\phi) } }
\newcommand{ \wpsiv }{ \omega_{\psi,V} }
\newcommand{ \mainquot }{ H (\F ) \backslash H _{ \fA }  }


\newcommand{ \opair }{ ( q , V ) } 
\newcommand{ \SLtwo }{ \widetilde{ SL }{ _2 }}
\newcommand{ \delo }{ \delta_0 } 
\newcommand{ \ompsi }{ \omega_{ \psi, V }  } 
\newcommand{ \invdist }{ \left (  S( V ) ^{ * }   \right ) ^{ O\opair } } 
\newcommand{ \Btilde }{ \widetilde{ B } } 
\newcommand{ \induce }{ Ind_{ \Btilde } ^{ \SLtwo } \mu } 
\newcommand{ \half }{ \frac{1}{2} }

\newcommand{ \hquot }{  H_{ \nu }  / H_{ \nu } ^{ v }  } 
\newcommand{ \ppp }[2]{ ( #1 , #2 )_{ \nu }  } 

\newcommand{ \ideles }{  \fA^{ \times }   } 
\newcommand{ \schbig }{ S( V_3( \F_\nu ) )  } 
\newcommand{ \schsmall }{ S( V_1( \F_\nu ) )  } 

\newcommand{\absv}[1]{\left| {#1} \right|_{\nu}}

\newcommand{ \intformula }[3]{  \int_{ \booyaquot } \left (  #1  \right ) (h^{ -1 } #2) \xi_{ V_3 } ( h ) \, dh  } 
\newcommand{ \formula }[3]{ \absv { #3 } ( #3 , -d(q_{ U } ) )_{ \nu }  \intformula {#1}{#2}{#3} } 
\newcommand{ \suffixformula }[4]{ \gamma^{ -1 } \absv{#1}^{ \frac{ #4 } { 2 } } ( #1 , #2d(#3 ))_{ \nu } } 
\newcommand{ \intformulaPAR }[3]{  \int_{ \booyaquot } #1(h^{ -1 } #2) \xi_{ V_3 } ( h ) \, dh  } 
\newcommand{ \formulaPAR }[3]{ \absv { #3 } ( #3 , -d(q_{ U } ) )_{ \nu }  \intformulaPAR {#1}{#2}{#3} } 

\newcommand{ \pad }{ \left (  \fA \right ) } 


 
\pagenumbering{roman}
\thispagestyle{empty}

\title[Morphisms of Theta Lifts]{
\large
Morphisms of global theta lifts
of Non-trivial Automorphic Characters \\ of Orthogonal Groups
 }
\begin{titlepage}
\begin{center}

\vspace*{2\baselineskip}

\end{center}

\setlength{\parindent}{0pt}
\setlength{\parskip}{0pt}
\vspace*{2\baselineskip}

\maketitle

\vspace{0.25cm}
\begin{center}
Ron Erez
\end{center}
\vspace{0.25cm}
\begin{abstract}
 This work is largely inspired by the 2003 Ph.D. thesis \cite{snitz} of Kobi Snitz.
In his thesis, Snitz constructed two irreducible, automorphic, cuspidal representations $ \pi $ and $ \pi' $ of the metaplectic group $ G\left (  \mathbb A \right ) = \widetilde{ SL }_{ 2 }  \left (  \mathbb A \right )  $ where each representation is obtained from a different global theta lifts of certain non-trivial automorphic characters $ \xi $ and  $ \xi' $ of the orthogonal groups $ H_{ \mathbb A } = O \left (  q, V   \right )   \left (  \mathbb A \right )    $ and  $ H_{ \mathbb A } '= O \left (  q', V' \right ) \left (  \mathbb A \right )   $, respectively, where $ \mathbb A = \mathbb A_{ \mathbb F } $ is the adele ring of a number field $ \mathbb F $. Snitz shows that for certain matching data of quadratic spaces and automorphic quadratic characters, that these two representations of 
$ G \left (  \mathbb A \right )  $
 are isomorphic, i.e. $ \pi\cong\pi'$. Moreover, he constructs an Eulerian global isomorphism 
and obtains explicit formulas, given by a certain orbital integral, for the corresponding local isomorphisms $ \pi_{ \nu } \cong \pi\p_{ \nu } $.
\par
The goal of this work is to reformulate and generalize Snitz's work to higher rank groups.
Namely we wish to determine for which admissible data $\left (   \left ( q, V \right )  ,\xi
,
 \left ( q', V' \right ),\xi'
  \right )$ satisfying certain local necessary conditions
 could an isomorphism possibly exist between two global theta lifts $ \pi $ and $ \pi'$ with respect to two reductive dual pairs 
$ H_{ \mathbb A } \times G_{ \mathbb A }  $
and
$ H'_{ \mathbb A } \times G_{ \mathbb A }  $
and two non-trivial automorphic quadratic characters $ \xi $ and  $ \xi'$ of the orthogonal groups $ H_{ \mathbb A } = O \left (  q, V   \right )   \left (  \mathbb A \right )    $ and  $ H_{ \mathbb A } '= O \left (  q', V'\right )  \left (  \mathbb A \right )    $, respectively and the group $ G $ which is the symplectic or the metaplectic group.
\par
In case these local necessary conditions are satisfied we also determine a global condition that will ensure that the two representations $ \pi $ and $ \pi'$  are isomorphic.
We then construct an explicit isomorphism between $ \pi $ and $ \pi'$ given by a certain non-degenerate pairing. The formulas obtained are similar to the Rallis inner product formula.

\end{abstract}

\vspace{1cm}
\noindent
\small{\emph{Keywords:} Theta lift; Automorphic Characters; CAP Representations}

\newpage
\tableofcontents

\newcommand{\Jn}{\nJ_n}
\newcommand{\wn}{\nw_n}

\end{titlepage}

\newcommand{ \iterateText }{ SECOND ITERATION } 
\newcommand{ \iterate }{
\boxed{  \colorbox{black!90}{
{\color{green} 
\iterateText }
{\color{cyan} \iterateText }
{\color{red} \iterateText }
}
}
\\
\despegar
} 

\newcommand{ \homc }[3]{ Hom_{ #1 }  \left ( #2 , #3 \right )  }
\newcommand{ \gen }{ \overline{ \omega } } 
\newcommand{ \adele }{ad\`{e}le } 
\newcommand{ \adelic }{ad\`{e}lic } 
\newcommand{ \adeles }{ad\`{e}les } 
\newcommand{ \dx }{ \tn d x } 
\newcommand{ \dd }[1]{  \tn d #1 } 
\newcommand{ \gammafactor  }{ \gamma \left (  a, \psi^{ -1 }  \right ) } 
\newcommand{ \gammafactorx  }{ \gamma \left (  x, \psi^{ -1 }  \right ) } 
\newcommand{ \hp }{ H^{ \prime }  } 

\newcommand{ \mquad }{ \left ( q_{ m } ,\xi_{ m } ,\qp_{ \mmp } ,\xip_{ \mmp }   \right ) } 
\newcommand{ \qp }{ q^{ \prime }  } 
\newcommand{ \xip }{ \xi^{ \prime }  }
\newcommand{ \vp }{ V^{ \prime }  } 
\newcommand{ \vpn }{ V^{ \prime,  n }  } 
\newcommand{ \dto }{ Z_{ \nu } ^{ \times }  } 
\newcommand{ \dt }{ Z_{ \nu } } 
\newcommand{ \oms }{ \Omega^{ \times }  }
\newcommand{ \theGroup }{ \widetilde{ SL }_{ 2 }  } 
\newcommand\blank[1]{\rule[-.2ex]{#1}{.4pt}}

\newcommand{ \zs }{ Z_{ s, \nu }  } 
\newcommand{ \zo }{ Z_{ o, \nu }  } 
\newcommand{ \tzs }{ \widetilde{ \zs } } 
\newcommand{ \tzo }{ \widetilde{ \zo } } 
\newcommand{ \pv }{ v^{ \prime }  } 
\newcommand{ \hi }{ H_{ 1 }  } 
\newcommand{ \hii }{ H_{ 2 }  } 
\newcommand{ \pref }[3]{ $ #3 $ \,
{ \small #2  }
\hfill
{\color{red}\pageref{#1}} \\
}

\newcommand{ \vpsi }{ V, \psi } 
\newcommand{ \vpsip }{ V\p, \psi } 
\newcommand{ \vpsiloc }{ V_{ \nu } , \psi_{ \nu }  } 
\newcommand{ \vpsilocp }{ V\p_{ \nu } , \psi_{ \nu }  } 
\newcommand{ \psiv }{ \vpsi } 
\newcommand{ \psivlocp }{ \vpsilocp } 
\newcommand{ \psivloc }{ \vpsiloc } 
\newcommand{ \psivp }{ \vpsiploc } 
\newcommand{ \INVpsiv }{ V, \psi^{ -1 } }
\newcommand{ \INVpsivp }{ V, \psi^{ -1 } }

\newcommand{ \smallRedComment }[1]{ \tiny{\color{red} (#1)}}
\newcommand{ \changeGen }{ { \tiny{\color{red} (Change this in case you do prove things in full generality.) }  } } 
\newcommand{ \review }{ \smallRedComment{ 
This looks wrong. I think you should review this.
  } } 

\newcommand{ \localPairingSig }{ B_{ \nu } \,\colon \schloc  \times \schlocp   \to \fC  } 
\newcommand{ \localPairingSigs }{ B_{ \nu, s } \,\colon \schloc  \times \schlocp   \to \fC  } 
\newcommand{ \localPairingSigsVal }[1]{ B_{ \nu, s = #1 } }
\newcommand{ \absA }[1]{ \abs{ #1 } _{ \fA }  } 
\newcommand{ \Hintegral }{ \int_{ \qAA { H^{ \vv } } H }  } 
\newcommand{ \Hintegralp }{ \int_{ \qAA { H^{\prime, \vv\p } }  { H\p } } } 
\newcommand{ \Hintegralnu }{ \int_{ \qFnFn { H^{ \vv } } H }  } 
\newcommand{ \Hintegralnup }{ \int_{ \qFnFn { H^{\prime, \vv\p } }  { H\p } } } 
\newcommand{ \HintegralF }{ \int_{ \qFnFn { H^{ \vv } } H }  } 
\newcommand{ \HintegralFp }{ \int_{ \qFnFn { H^{\prime, \vv\p } }  { H\p } } } 

\newcommand{ \HintegralFW }{ 
\int_{
{ H^{ \vv } } _{ \F }  \backslash H _{ \F } / K
 }
 } 

\newcommand{ \vali }{ 3/2 }
\newcommand{ \consti }{  } 
\newcommand{ \constiINV }{  } 
\newcommand{ \valii }{ 1/2 }
\newcommand{ \constii }{  } 
\newcommand{ \constiIiiINV }{  } 
\newcommand{ \wfunc }{ \omega_{ \vpsi } } 
\newcommand{ \wfuncp }{ \omega_{ \vpsip } } 
\newcommand{ \wfuncpi }{ \omega_{ V\p,\psi^{ -1 } } } 
\newcommand{ \vvi }{ v_{ 1 }  }
\newcommand{ \vvip }{ v_{ 1 } \p }
\newcommand{ \delb }[1]{  \delta_{ P }^{ -1 } \left ( #1 \right ) }
\newcommand{ \delBorel }[1]{  \delta_{ B }^{ -1 } \left ( #1 \right ) }
\newcommand{ \delbc }[1]{  \delta_{ P } \left ( #1 \right ) }
\newcommand{ \delbp }[2]{  \delta_{ P }^{ #2 }  \left ( #1 \right ) }
\newcommand{ \com }[1]{ \,\,\,\,\,\,\,\,\left [ \tn{ #1 } \right ] }
\newcommand{ \fthLift }[1]{ 
I_{ V, \xi }^{ #1 }
 } 
 \newcommand{ \fthLiftp }[1]{ 
I_{ V\p, \xi\p }^{ \prime, #1 }
 } 
\newcommand{ \fthLiftBeta }[1]{ 
I_{ V, \xi, \beta }^{ #1 }
 } 
 \newcommand{ \fthLiftBetao }[1]{ 
I_{ V, \xi, \betao }^{ #1 }
 }

 \newcommand{ \thLiftBeta }[2]{ 
\fthLiftBeta{ #1 } \left ( #2 \right )
  }

\newcommand{ \thLiftBetao }[2]{ 
\fthLiftBetao{ #1 } \left ( #2 \right )
}
  
\newcommand{ \thLift }[2]{ 
 \fthLift{ #1 } \left ( #2 \right )
 } 
\newcommand{ \thLiftp }[2]{ 
 \fthLiftp{ #1 } \left ( #2 \right )
 } 
\newcommand{ \ga }[1]{  #1 _{ \fA }   }
\newcommand{ \fa }[1]{  #1 _{ \nu }   }

\newcommand{ \betao }{ \beta_{ 0 } } 

\newcommand{ \Vn }{ V } 

\newcommand{ \signm }[1]{ \left ( -1 \right )^{ \frac{ #1 \left ( #1 - 1 \right ) } { 2 } } }

\newcommand{ \tB }{ \widetilde{ B } } 
\newcommand{ \ttB }{ \widetilde{ \widetilde{ B } } } 
\newcommand{ \spquot }{ G_{ \F }  \backslash G\ad } 
\newcommand{ \ovquot }[1]{ H#1_{ \F }  \backslash H#1\ad } 
\newcommand{ \shspquot }{ \left [ G \right ] } 
\newcommand{ \wf }[5]{ \int_{ \ovquot{#3} } 
\sum_{ 
\overset{ \gamma \in \Vn_{ \F } }{ #5 }
 }
#4
 \left (  
\omega_{ \vpsi } 
 \left ( #1, #2 \right ) 
\phi
  \right ) 
 \left (  \gamma \right ) 
 \xi \left ( h \right ) 
\,dh
 }
\newcommand{ \wfNosum }[9]{ \int_{ \ovquot{#3} } 
#6
#4
 \left (  
\omega_{ \vpsi } 
 \left ( #1, #2 \right ) 
\phi
  \right )
\left ( #7 \right )
 \xi \left ( 
#8
  \right )
  #9
\,dh
 }
\newcommand{ \lw }[2]{ \omega_{ \vpsi } 
 \left ( #1, #2 \right ) 
 }
\newcommand{ \Vi }{ v_{ 1 }  } 
\newcommand{ \Vii }{ v_{ 1 } \p } 
\newcommand{ \al }[1]{   \alpha \left ( #1 \right ) }
\newcommand{ \as }[2]{ 
\abs{  \alpha \left ( #1 \right ) } ^{ #2 } 
  } 
\newcommand{ \asp }{  \alpha\p \left ( s\p \right )  } 
\newcommand{ \mm }{  m \left ( a \right )  } 
\newcommand{ \mmi }{  m \left ( a^{ -1 } \right ) }
\newcommand{ \nn }{ n \left ( b \right ) } 
\newcommand{ \mn }{  n \left ( b \right )  \mm  }

\newcommand{ \bphi }{ \Phi } 
\newcommand{ \psq }{  \psi \left ( b  q \left ( av_1 \right )  \right )  } 
\newcommand{ \pl }[1]{ \phi \left (  #1 a v_{ 1 } \right )  }
\newcommand{ \plp }[1]{ \phi\p \left (  #1 a v_{ 1 }\p \right )  }
\newcommand{ \ppl }[1]{ \psq \phi \left (  #1 a v_{ 1 } \right )  }
\newcommand{ \pplp }[1]{ \phi\p \left (  #1 a v_{ 1 }\p \right )  }
\newcommand{ \psibionly }{ \psi_{ \beta } ^{ -1 }  \left (  \nn  \right ) } 
\newcommand{ \psibi }{ \chisa \psibionly \,db \,d^{ \times }  a } 
\newcommand{ \chisa }{ \chi_{ s }  \left (  a \right )  } 
\newcommand{ \rquot }[2]{ #1 _{ \F }  \backslash #2 _{ \F }  }
\newcommand{ \aq }[1]{ #1 _{ \F }  \backslash #1 _{ \fA }  }

\newcommand{ \qFA }[2]{ #1 _{ \F }  \backslash #2 _{ \fA }  }
\newcommand{ \qFF }[2]{ #1 _{ \F }  \backslash #2 _{ \F }  }
\newcommand{ \qFnFn }[2]{ #1 _{ \nu }  \backslash #2 _{ \nu }  }
\newcommand{ \qAA }[2]{ #1 _{ \fA }  \backslash #2 _{ \fA }  }

\newcommand{ \fdot }{ \bullet }

\newcommand{ \cc }{ c } 
\newcommand{ \cci }{ \cc^{ -1 } } 
\newcommand{ \ddi }{ d } 
\newcommand{ \fPsi }[2]{ \Psi_{ #2 }  \left ( #1 \right ) }
\newcommand{ \aeis }[2]{ \cE \left ( #1, #2 \right ) }
\newcommand{ \aeisb }[2]{ \cE_{ \beta } \left ( #1, #2 \right ) }
\newcommand{ \aeisSec }[3]{ \cE \left ( #1, #2, #3 \right ) }
\newcommand{ \aeisSecb }[3]{ \cE_{ \beta } \left ( #1, #2, #3 \right ) }

\newcommand{ \nbs }{  n \left ( b \right )  } 
\newcommand{ \mas }{  m \left ( a \right )  } 
\newcommand{ \umu }{ \underline \mu } 
\newcommand{ \locSat }[2]{ \abs{ a_{ { #1#1 }  }   }^{ #2 }  }
\newcommand{ \locSatp }[2]{ \abs{ a_{ { #1,#1 }  }   }^{ #2 }  }
\newcommand{ \nuc }[1]{ \mu_{ #1, \nu }  }
\newcommand{ \slam }{ \lam }

\newcommand{ \recallForMoron }{ Recall that if $ g \in G $ where $ G $ is the symplectic or metaplectic group then by the Iwasawa decomposition 
$$
G = PK = MNK\cong GL_{ n } NK
$$
where $ P  $ is the Siegel parabolic of $ G $ and $ M $ is the Levi part of $ P $ and $ N $ is the unipotent radical of $ P $. The group $ M $ is of the form 
$$
M = \Set*{  m \left ( a \right )  = \begin{pmatrix} a & 0 \\ 0 & a^\star \end{pmatrix} } { a \in GL_{ n } }
$$
so for $ m =  m \left ( a \right )  \in M $ we define $  a \left ( g \right )  = a $, in other words $a \,\colon G \to GL_n$. Finally we denote
$$
\abs{  a \left ( g \right )  } 
 = 
 \abs{  \det \left ( a \right )  } 
$$
and 
$$
\begin{cases}
   \Psi_{ s }  \,\colon GL_n(\fA) \to \F \\
   \fPsi{ g }{ s }=\absA{  a \left ( g \right ) }^{ s + \rhon } 
   \\
\end{cases}
$$
where $ \rhon = \dfrac{ n+1 } { 2 } $. }

\newcommand{ \naeis }[2]{ \cE^{ \star } \left ( #1, #2 \right ) }
\newcommand{ \apair }{ B_{ s } }

\newcommand{ \czetaN }{ \Lambda } 
\newcommand{ \czeta }[1]{ \czetaN \left ( #1 \right )  }

\newcommand{ \glnfnu }{ GL_{ n } \left ( \F_{ \nu } \right )  } 
\newcommand{ \dpairs }[2]{
\centerline{
\xymatrix{
& \widetilde{ Sp \left ( W \right )  } \left (  \fA \right )  \ar@{-}[dl]\ar@{-}[dr]
& 
\\
O \left (  #1 \right )  \left (  \fA \right ) 
& 
& O \left (  #2 \right )  \left (  \fA \right ) 
}
} }

\newcommand{ \dpairss }[2]{
 \centerline{
\xymatrix{
& \widetilde{ SL }_{ 2 } \left (  \fA \right )  \ar@{-}[dl]\ar@{-}[dr]
& 
\\
O \left (  #1 \right )  \left (  \fA \right ) 
& 
& O \left (  #2 \right )  \left (  \fA \right ) 
}}
 }

\newcommand{ \ovP }{ \overline{ P } }
\newcommand{ \lp }[1]{ \left(  #1 ,#1\p \right ) }
\newcommand{ \hhpLowerSpace }{ \hhp^{ w_{ 0 }  } }


\makeatletter
\def\@tocline#1#2#3#4#5#6#7{\relax
  \ifnum #1>\c@tocdepth 
  \else
    \par \addpenalty\@secpenalty\addvspace{#2}%
    \begingroup \hyphenpenalty\@M
    \@ifempty{#4}{%
      \@tempdima\csname r@tocindent\number#1\endcsname\relax
    }{%
      \@tempdima#4\relax
    }%
    \parindent\z@ \leftskip#3\relax \advance\leftskip\@tempdima\relax
    \rightskip\@pnumwidth plus4em \parfillskip-\@pnumwidth
    #5\leavevmode\hskip-\@tempdima
      \ifcase #1
       \or\or \hskip 1em \or \hskip 2em \else \hskip 3em \fi%
      #6\nobreak\relax
    \hfill\hbox to\@pnumwidth{\@tocpagenum{#7}}\par
    \nobreak
    \endgroup
  \fi}
\makeatother

\newcommand{ \newtitle }{ {Matchings of theta lifts associated\\
to Non-trivial Automorphic Characters\\ of Odd Orthogonal Groups} } 

\newcommand{ \DEPinput }[1]{}

\newcommand{ \ininput }[1]{\input{#1}}
\newcommand{ \gLift }{ \Theta \left ( q,V,\xi \right ) } 
\newcommand{ \gLiftp }{ \Theta \left ( q\p,V\p,\xi\p \right ) } 
\newcommand{ \egLift }{ \Theta \left ( \eta \right ) } 
\newcommand{ \egLiftp }{ \Theta \left ( \eta \p \right ) }

\newcommand{ \tH }{ H } 
\newcommand{ \tG }{ G } 

\newcommand{ \wwpsi }[1]{  \omega_{ \psi_{ #1 } } }
\newcommand{ \repwpsi }[1]{ \rep{ Mp \left ( W \right ) }{ \wwpsi{ #1 } }{ S } }
\newcommand{ \weilcond }{  A \left ( g \right ) ^{ -1 }  \rhop \left ( h \right )  A \left ( g \right ) = 
  \rhopg \left ( h \right )  } 
\newcommand{ \tspww }{ \widetilde{ \spww } } 
\newcommand{ \tspwwp }{ \tspww_{ \psi }  } 
\newcommand{ \spww }{ Sp \left ( \fW \right ) } 
\newcommand{ \heisww }{ H \left ( \fW \right ) } 
\newcommand{ \rep }[3]{   \left ( #1, #2, #3 \right )  }
\newcommand{ \rhop }{ \rho_{ \psi }  }
\newcommand{ \rhopg }{ \rho_{ \psi }^{ g }   }
\newcommand{ \ts }[1]{ x_{ #1 } \otimes y_{ #1 } }

\newcommand{ \chiv }[1]{ \chi_{ V_\nu }  \left ( #1 \right )  } 
\newcommand{ \chivp }[1]{ \chi_{ V_\nu\p }  \left ( #1 \right )  }

\newcommand{ \nicequotmv }{ H^{ \vv} _{ \fA }  \backslash H _{ \fA }   } 

\newcommand{ \gha }{ G_{ \fA } \times H_{ \fA }  } 
\newcommand{ \ghnu }{ G_{ \nu } \times H_{ \nu }  } 
\newcommand{ \omegavpsi }{ \omega_{ \vpsi } } 
\newcommand{ \omegapsip }{ \omega_{ \vpsi } } 

\newcommand{ \T }{ r }
\newcommand{ \chForInduced }{ \alpha } 
\newcommand{ \soo }[3]{ 
\centerline{
\xymatrix{
& 
#1
 \ar@{-}[dl]\ar@{-}[dr]
& 
\\
#2
& 
& 
#3
}
} 
 }

\newcommand{ \piat }{Piatetski-Shapiro } 
 
 \newcommand{ \skewPairing }{ 
S \left (  \voxpr _{ \fA }  \right ) \times V_{ \pi^\vee } 
 }
 
\newcommand{ \mainFquot }{ \hhp_{ \F } / \hhpLowerSpace_{ \F } }
\newcommand{ \nz }{ Z_{ new } }

\newcommand{ \forex }[5]{ For example if  $ \beta = \begin{pmatrix} #1 & #2 \\ #3 & #4 \end{pmatrix} $ then 
$$
G_{ \beta } = 
\left \{
#5
\right \}
$$ }

\newcommand{ \distYP }[3]{ \begin{align*}
 Y\p \left (  \phi,\phi\p \right ) & =  \int_{ N }\int_{ M }  \int_{ H / H^{ v_1 }  }  \int_{ H\p / H^{\prime, v_1\p }  }   \left (  \mn \cdot \phi \right )   \left (  h^{ -1 } v_{ 1 }  \right )   \left (  \mn \cdot \phi\p \right )   \left (  h^{ \prime, -1 } v_{ 1 }  \right )  
 \\
&
  \xi \left ( h \right )   \xi\p \left ( h\p \right )  \,dh \,dh\p \,dndm \\
& =
\int_{ M }  \int_{ H / H^{ v_1 }  }  \int_{ H\p / H^{\prime, v_1\p }  }   \left (  \mm \cdot \phi \right )   \left (  h^{ -1 } v_{ 1 }  \right )   \left (  \mm \cdot \phi\p \right )   \left (  h^{ \prime, -1 } v_{ 1 }  \right )   \xi \left ( h \right )   \xi\p \left ( h\p \right )  \,dh \,dh\p \,dn \\
\\
& =
\int_{ \F^{ \times }  }   
\abs{ a } ^{ 
\frac{ #1 } { 2 }
+ 
\frac{ #2 } { 2 }
 }  \chi \left ( a \right ) 
 \int_{ H / H^{ v_1 }  }  \int_{ H\p / H^{\prime, v_1\p }  } \phi  \left (  h^{ -1 } av_{ 1 }  \right )   \phi\p\left (  h^{ \prime, -1 } av_{ 1 }  \right )   \xi \left ( h \right )   \xi\p \left ( h\p \right )  \,dh \,dh\p \,d^{ \times } a \\
\\
& =
\int_{ \F^{ \times }  }
\abs{ a } ^{ #3 }  \chi \left ( a \right ) \int_{ H / H^{ v_1 }  }  \int_{ H\p / H^{\prime, v_1\p }  } \phi  \left (  h^{ -1 } av_{ 1 }  \right )   \phi\p\left (  h^{ \prime, -1 } a v_{ 1 }  \right )   \xi \left ( h \right )   \xi\p \left ( h\p \right )  \,dh \,dh\p \,d^{ \times } a \\
\\
& =
\gooddist\left (  \phi,\phi\p \right )
\end{align*}
}

\newcommand{ \vol }[1]{ vol \left ( #1 \right )  }

\newcommand{ \mw }{  m \left ( \varpi^{ J }  \right )  } 
\newcommand{ \wj }{ \varpi^{ J }  } 
\newcommand{ \wji }[1]{ \varpi^{ #1 } }

\newcommand{\omi}{\omega_1}
\newcommand{\omii}{\omega_2}

\newcommand{\V}{\rm V}
\newcommand{\refl}[1]{\tau_{#1}}
\newcommand{\refli}{\refl{\omi}}
\newcommand{\reflii}{\refl{\omii}}

\newcommand{ \chik }{ \chi_{ \kappa,\nu }  }

\newcommand{ \uuu }[1]{\underline{#1}}
\newcommand{ \bbb }[1]{{\bfseries{#1}}}

\newcommand{ \liftc }[1]{ \gco{ \psi, V , W_{ #1 }  }{ \xi } }
\newcommand{ \lift }[1]{ \gco{ \psi, V , W_{ #1 }  }{ \sigma } }

\newcommand{ \bdy }[1]{
H\ad
\times 
G_{ #1, \fA } 
   }
\newcommand{ \plift }[2]{ 
 $ \left (  
H \left (  \fA \right ) ,
 Sp \left ( W_{ #1 }  \right ) 
  \right ) $
& - &
 $  \liftc{ #1 }  $
 }
 
 \newcommand{ \wint }{
{ q\left [ \ww \right ] = 
  q\left [ \vv \right ]  }
}

\newcommand{ \afac }[1]{ \abs{  \det \left ( a \right ) } ^{ #1 / 2 }  }

\newcommand{ \glno }{ GL_{ n }  \left (  \cO_{ \nu } \right ) }

\newcommand{ \cFunc }[2]{ \begin{cases}
  1 & , #1 \in \left (  #2 \vox \right )  \left (  \cO_{ \nu }  \right ) \\
   0  & , #1 \not\in \left (  #2 \vox \right )  \left (  \cO_{ \nu }  \right ) \\
\end{cases}
 } 
 
\newcommand{ \wwp }{ W W\p \left ( a \right ) }
\newcommand{ \av }{ \abs{ a } } 

\newcommand{ \charFunc }{    1_{  \left (  \vox \right )  \left (  \cO_{ \nu }  \right )  }  } 

\newcommand{ \di }{ m } 
\newcommand{ \dii }{ m\p } 

\newcommand{ \mnm }[2]{ $ m = #1 $ and $ m\p = #2 $ }

\newcommand{ \boo }[1]{   $$
\tn{ The theta lift }
 \liftc{ #1 }
\tn{ with reductive dual pair }
 \bdy{ #1 }$$ }
 \newcommand{ \vsp }{ $$
\!\!\!\!\!\!\!\!\!\!\!\!\!
\vdots
\qquad\qquad\qquad\qquad\qquad
\vdots
$$
 } 

\newcommand{ \FC }{  \left (  FC \right )  } 
\newcommand{ \CC }{ \left (  CC \right ) } 

\newcommand{ \bad }{{\color{magenta} - BAD SECTION}} 

\newcommand{ \jcond }{ 0\leq j_{ n } \leq \cdots  \leq j_{ 1 }  } 

\newcommand{ \gram }[1]{ \gramdef{ q }{ #1 } }
\newcommand{ \gramqm }[1]{ \gramn{ m }{ #1 } }
\newcommand{ \gramn }[2]{ \gramdef{ q_{ #1 } }{ #2 } }
\newcommand{ \gramp }[1]{ \gramdef{ q\p }{ #1 } }
\newcommand{ \gramdef }[2]{ #1 \left [ #2 \right ] }

\newcommand{ \here }{ {\color{red} HERE } } 

\newcommand{ \myH }{ H / H^{ \vv }  } 
\newcommand{ \ct }{ \dfrac{  \chi \left ( 2 \right )  } { 2 } } 
\newcommand{ \overAlpha }{ \overline{ \alpha  } } 
\newcommand{ \q }{ p } 
\newcommand{ \f }{ \chi_{ 0 }  } 
\newcommand{ \ff }[2]{ \f \left ( #1, #2 \right )  }
\newcommand{ \lRing }[1]{ 
\q^{ #1 } \cO_{ \nu } 
 }
\newcommand{  \lRingDiff }[2]{ \lRing{ #1 }
 \setminus
 \lRing{ #2 }  } 

\newcommand{ \chiNice }[2]{ \left [ 
\chi \left ( 
#1  + 1
\right ) 
 + 
 \chi \left (  #2
-
1 \right ) 
  \right ] }

\newcommand{ \chiNicer }[1]{ \chiNice{#1}{#1} }

\newcommand{ \chiBetter }[2]{ \left [ 
\chi \left ( 
#1  + #2
\right ) 
 + 
 \chi \left (  #1
-
#2 \right ) 
  \right ]  }
  
\newcommand{ \dble }[3]{ #1 \backslash #2 / #3 }
\newcommand{ \HvHK }{ H^{ \vv }  _{ \F }  \backslash H _{ \F } / K } 
\newcommand{ \pa }{ P_{  w }  } 
\newcommand{ \kk }{ \kappa } 
\newcommand{ \Hv }{ H^{ \vv }  } 
\newcommand{ \HintegralFK }{ 
\int_ \HvHK
 } 
\newcommand{ \bq }[2]{ \left < #1, #2 \right > }
\newcommand{ \ws }{ w^{ \star }  }

\addcontentsline{toc}{section}{Abstract}

\pagenumbering{arabic}

\newpage

\section{ Introduction }

\subsection{ Motivation }
This work is inspired by the 2003 Ph.D. thesis \cite{snitz} of Kobi Snitz.
In his work Snitz constructed two irreducible, automorphic, cuspidal representations $ \pi $ and $ \pi\p $ of the metaplectic group $ G_{ \fA } = \widetilde{ SL_2 } \left ( \fA \right )  $ over the \adele ring $ \fA $ via the global theta correspondence by lifting certain automorphic nontrivial quadratic characters of two different orthogonal groups $ H\ad $ and $ H\p\ad $ to the group $ G\ad $. Snitz proved that for certain defining data, namely for certain pairs of quadratic spaces and non-trivial automorphic quadratic characters that the representations 
$ \pi $ and $ \pi\p $
given by this data are isomorphic.
In addition, Snitz constructed explicit global and local isomorphisms between these two representations.
\par
The goal of this work is to reformulate Snitz’s work so that similar or equivalent results may be generalized to higher rank groups, i.e. we would like to 
consider certain trivial \underline{local} necessary conditions on our defining data such that the corresponding global theta lifts $ \pi $ and $ \pi\p $ could potentially be isomorphic.
 This will be given by a non-trivial \underline{global} condition given in terms of a partial global $ L $-function. Finally, in case these representations are isomorphic we will construct explicit local and global isomorphisms between them.
\subsubsection{ The global theta correspondence in a nutshell }
The {\bfseries{global theta correspondence}} is a correspondence between automorphic representations of a group $ H_{ \fA }  $ and automorphic representations of another group $ G_{ \fA }  $, where a typical case to keep in mind is $ H $ being an orthogonal group and $ G $  is either the symplectic group $ Sp \left (  W \right )  $ or a double cover of the symplectic group $ \widetilde{ Sp } \left (  W \right )  $ called the \bbb{metaplectic group}.
The image of this correspondence is often referred to as a “global theta lift” from $ H_{ \fA }  $ to $ G_{ \fA } $ and at times we will simply say a “lift”.
\par
Note that the pair of groups $ H_{ \fA } \times G_{ \fA }  $ are not arbitrarily chosen and are called a \bbb{reductive dual pair} (see Definition \ref{sec:ReductiveDualPairs}, page~\pageref{sec:ReductiveDualPairs}, Section \ref{sec:redDualPair} for more details). 
In this work, the only reductive dual pairs that will be considered are
$
 O \left ( q,V \right )  \ad  \times Sp \left (  W \right )  \ad  
$
when $  \dim \left ( V \right )  = m $ is even
or
$
 O \left ( q,V \right )  \ad  \times \widetilde{ Sp } \left (  W \right )  \ad  
$ 
when $  \dim \left ( V \right )  $ is odd,
where $  \left (  q,V \right ) $ is a quadratic space
and $  \left (  \left < \cdot,\cdot \right >, W \right )  $ is a symplectic vector space with $ \dim\,W = 2n $.
Moreover $ W $ has a polar decomposition $ W = X + Y $ where $ X $ and $ Y $ are maximal isotropic subspaces of $ W $. A standard symplectic basis $  \left \{ e_{ 1 } , \ldots, e_{ n } ,f_{ 1 } , \ldots, f_{ n }   \right \} $ can be chosen so that for every $ 1\leq i,j\leq  n $
$ 
\left < e_{ i } ,e_{ j }   \right > = 
\left < f_{ i } ,f_{ j }   \right > = 
0,
\left < e_{ i } ,e_{ j }   \right > = \delta_{ ij } 
 $.

\par
The defining data in this work will be denoted by $ \adata $ where the defining data of the representation $ \pi $ of $ G_{ \fA } $, coming from the theta correspondence, is the quadratic space $  \left (  q,V \right )  $ together with a non-trivial automorphic quadratic character $ \xi \,\colon H_\fA \to \fC^\times $ and then we will denote the \bbb{global theta lift} corresponding to the data
$ \qcdata $
 by 
$$
\pi = \gLift = \egLift
$$
\par
Similarly, for the quadratic space $\left ( q\p, V\p \right ) $ together with a non-trivial automorphic quadratic character $ \xi\p \,\colon H_\fA\p \to \fC^\times $ we obtain a representation $ \pi\p $ of $ G_{ \fA } $ and we denote the global theta lift corresponding to the data $ \qcdatap $ by
$$
\pi\p = \gLiftp = \egLiftp
$$
\par
Certain elements generating $ \pi $ and $ \pi\p $ will be called \bbb{theta integrals} and will be denoted by $ I^{ \phi }_{ V, \xi }  $ and $ I^{  \phip  }_{ V\p, \xi\p } $, respectively, where $ \phi $ and $ \phi\p $ are Schwartz functions. The theta integrals are functions on $ \aq G $.
\par
Moreover, we will prove that the representations $ \pi $ and $ \pi\p $ satisfying certain conditions in this work are irreducible, cuspidal, automorphic representations of the group $ G\ad $. By Flath's theorem \cite{flath1979decomposition} these representations can be written as a restricted tensor product of local representations of the group $ G\loc $, i.e. 
$$ \pi \cong 
\otimes \loc ^{ \p } \pi\loc 
\tn{ and }\pi\p \cong 
\otimes \loc ^{ \p } \pi\p\loc  
$$

\subsubsection{ An example of lifting a trivial automorphic character }
For example one could lift a trivial automorphic character $ \xi \,\colon H_\fA \to \fC^\times $ (that is $ \xi\equiv 1 $) of the orthogonal group $ H\ad $ to an automorphic representation of the symplectic (or metaplectic) group $ G\ad $. This lift is characterized by what is known as the Siegel-Weil formula and has deep arithmetical significance, generalizes work of Siegel on representation numbers of quadratic forms and has many important applications in the theory of automorphic forms.
See Appendix \ref{sec:ExampleLiftingCharacters}, page \pageref{sec:ExampleLiftingCharacters} for more details.

\subsubsection{ An example of lifting a non-trivial automorphic quadratic character }
Another natural example would be to lift a non-trivial automorphic quadratic character $ \xi \,\colon H_\fA \to \fC^\times $ of an orthogonal group $ H_{ \fA }  $ to an automorphic representation of the symplectic (or metaplectic) group $ G\ad $. For instance $ \xi $ could be taken to be trivial at almost every place outside of a finite set $ S $ (such that $ \abs{ S } $ is even) and equal to the determinant character for every place $ \nu\in S $, i.e. 
$
\xi =
 \left (  \otimes_{ \nu\in S }
\,
det\loc  \right ) 
\otimes 
 \left (  \otimes_{ \nu \not \in S }
\,
1\loc  \right ) 
$.
\par
Another example would be to write the global character $ \xi \,\colon H_\fA \to \fC^\times $ as the composition of a Hilbert symbol $  \left (  \cdot, \lam \right )\ad $

and the spinor norm $ SN $
of the orthogonal group.
In this case the local characters $ { \xi\loc  \,\colon H\loc \to \fC^\times  }$ are not necessarily trivial at almost every place.
\par
Lifts of non-trivial automorphic quadratic characters of the orthogonal group are the primary objects of study in this work and Snitz addressed the question of when two such lifts are isomorphic for certain low rank orthogonal groups over the \adeles of dimensions three and one.
The type of characters he considered are precisely the two types of characters mentioned in the previous two paragraphs.
See Appendix \ref{sec:ExampleLiftingCharacters} on page \pageref{sec:ExampleLiftingCharacters} for more details on Snitz's work and how the problem of lifting non-trivial automorphic characters of the orthogonal group in both Snitz's and our case is related to CAP representations and the generalized Ramanujan conjecture.

\subsubsection{ Goal of this work }

Philosophically, this work can be divided into three main steps, briefly described as follows. Of course in practice many parts of these steps are intermingled.
\begin{enumerate}
\item
\bbb{Local necessary conditions.}
Find necessary conditions for $ \pi $ and $ \pi\p $ to be isomorphic. In particular if there were a global isomorphism $ \pi\cong\pi\p $ this would imply the existence of local isomorphisms $ \pi\loc \cong\pi\p\loc  $ at every place $ \nu $.
\par
Hence we consider necessary local conditions on our local data $ \adatanu $ so that one could possibly have an isomorphism between $ \pi\loc  $ and $ \pi\p\loc  $ at every place $ \nu $.
\par
The two necessary conditions we consider for $ \pi\loc $ and $ \pi\p\loc  $ to be
isomorphic are called
\begin{enumerate}
\item
The \bbb{central character condition  $ \CC $}.
\item
The \bbb{Fourier coefficient condition $  \FC  $}.
\end{enumerate}
These conditions $ \CC $ and $ \FC $ will determine the local data $ \adatanu $. Such local data $ \alpha \loc  $ satisfying both $ \CC $ and $ \FC $ will be called a \bbb{locally admissible quadruple}.
\par
Moreover, the local data $ \alpha \loc  $ will determine conditions on our global data $ \adata $. Such global data will be referred to as a \bbb{globally admissible quadruple}.
Note that the local representations $ \pi\loc  $ and $ \pi\p\loc  $  are called \bbb{local theta lifts}.
\item
\bbb{Global conditions.}
Given a globally admissible quadruple $ \adata $, assume that a certain global condition is satisfied where this global condition is given in terms of a global partial $ L $-function. Under these conditions, prove that the corresponding global automorphic representations $ \pi $ and $ \pi\p $ are isomorphic.
\par
Note that the global condition will be given in terms of a global partial $ L $-function which will come from the \bbb{unramified computation}. This is a computation done in a generic, i.e. unramified, setting.
\item
\bbb{Construction of explicit morphisms.}
Construct explicit global morphisms $ \pi\to\pi\p $ and for every place $ \nu $ construct explicit local morphisms $ \pi\loc \to\pi\p\loc $.
More precisely, we will find an explicit non-trivial local pairing
$
\pi\loc \times \contra\pi\loc\p\to\fC
$
 and an explicit non-trivial global pairing
 $
\pi\times \contra\pi\p\to\fC
$
 between the local and global representations, respectively, satisfying certain equivariance properties.
 \par
This is equivalent to finding local and global isomorphisms between the representations under consideration.
It turns out that constructing such pairings is a much more accessible and natural problem than constructing explicit isomorphisms.
\end{enumerate}

\subsubsection{ Simplifying Assumptions }
Throughout this work we will be working under the simplifying assumption that the automorphic quadratic characters $ \xi \,\colon H_\fA \to \fC^\times $ and $ \xi\p \,\colon H_\fA\p \to \fC^\times $ are non-trivial and that the quadratic spaces $  \left (  q,V \right ) $ and $ \left (  q\p,V\p \right ) $ are globally anisotropic.
\par
Moreover, as aforementioned throughout this work we will be working under the simplifying assumption that the automorphic quadratic characters are non-trivial and that the quadratic spaces of interest are globally anisotropic. This will not be mentioned again.

\subsection{ Main Statement }

It is important to note that the main results of this work can only be stated for the case $ n = 1 $ since we were unable to do the unramified computation for $ n $ larger than one. Besides the unramified computation all results are presented and proved in full generality.
\par
Note that although we have not proven that the local factors are Eulerian for \( n>1 \) we do believe that such is the case.

\begin{theorem}[Main statement]
\label{sec:mainStatement}

Given a globally admissible quadruple $ \adata $.
If a certain global partial $ L $-function $ \Delta_{ S }  $
has a pole at
$ s = \rhon := 
\frac{ n+1 } { 2 }
 $
then the global theta lifts 
$$
\pi = \gLift\tn{ and }
\pi\p = \gLiftp
$$
are isomorphic.

\end{theorem}


\subsubsection{ Main Statement in terms of pairings }
\label{sec:alternativeMainStatement}
\par
For a globally admissible quadruple 
$ \adata $ we will construct two global pairings and then compare them. We will see that they are equal up to a constant where the constant is given by a global partial $ L $-function.
\par
Let $ X $ be a maximal isotropic subspace of a symplectic space
 $  \left (  W,\left < \cdot, \cdot \right > \right )  $.
The first global pairing we consider is a natural global pairing 
$$
B \,\colon \schad \times \schadp  \to \fC
$$
given by the Petersson inner product of two global theta integrals.
Such global theta integrals generate the representations $ \pi $ and $ \pi\p $ of $ G\ad $.
The other global pairing will be given by a product of certain explicit local pairings
$$
B\loc  \,\colon \schloc \times \schlocp  \to \fC 
$$

A key point in the proof is the existence or non-existence of a certain partial $ L- $function $ \Delta_S \left ( s \right ) $ which comes from the unramified computation of certain local pairings $ B_{ s,\nu } $. Unfortunately we do not know how to compute or even formulate a conjecture as to what this should be equal to in the higher rank setting. For $ n = 1 $ and for an admissible quadruple $ \Delta_S \left ( s \right ) $ is given by 
$$
\Delta_S \left ( s \right ) = 
\prod_{ \nu \not \in S } 
\left ( 1-q_{ \nu } ^{ -s } \right )^{ -1 } 
$$

\begin{theorem}[Alternative main statement in terms of pairings]
\label{sec:Alternative}
Suppose $ \adata $ is a globally admissible quadruple.
Then
\begin{enumerate}
\item
There exists a natural global pairing 
$$
B \,\colon \schad  \times \schadp  \to \fC
$$
which is Eulerian (i.e. factors as a product of local pairings) and satisfies for every factorizable
$ 
\phi = \otimes \loc ^{ \p } \phi\loc \in \schad  
$ and $
\phip = \otimes \loc ^{ \p } \phi\p\loc \in \schadp
 $
where S is a sufficently large set of primes that may depend on $ \phi $ and $ \phi\p $.

 \begin{align*}
    B \left ( \phi, \phip \right )  & =
\kappa
\cdot
    Res_{ s = \rhon }
\Delta_S \left ( s \right ) 
\cdot
\prod_{ \nu \in S } 
 B\loc 
\left ( \phi_{ \nu } , \phip _{ \nu } \right ) 
 \end{align*}
where $ \rhon = \dfrac{ n+1 } { 2 } $, $  \Delta_S $ is a certain global partial $ L- $function, $ S $ is a finite set of places
which depends on $ \phi $ and $ \phi\p $
, $ B\loc  $ is a certain explicit local pairings, and $ \kappa $ is a non-zero constant coming from the residue of a certain simple normalized auxiliary Eisenstein series $ \cE^{ \star } $.

\item

Moreover, if the global partial $ L $-function, $ \Delta_S \left ( s \right ) $, has a 
simple
 pole at $ s = \rhon $
then $ \pi $ and $ \pi\p $ are isomorphic.
\end{enumerate}
\end{theorem}
\begin{remark}
The global factor $ \Delta_S \left ( s \right ) $ is used to determine whether or not a non-zero global isomorphism exists between $ \pi $ and $ \pi\p $. This is a global condition that is not accounted for by the local necessary conditions. The formula for $ \Delta_S \left ( s \right ) $ follows from the unramified computation, i.e. it follows from a "generic" computation.
\end{remark}

\vspace{0.5cm}
\begin{center}
    {\bfseries{Acknowledgements}}
\end{center}

I would like to thank Prof. Joseph Bernstein for numerous conversations regarding the reformulation of the problem at hand.
Any errors or inaccuracies in this work are solely my responsibility.

\newpage
\section{ Overview of proof }

In this section we give a succinct yet slightly more detailed account of the main steps in the proof of the main statements Theorems \ref{sec:mainStatement} and \ref{sec:Alternative}.

\subsection{ Local necessary conditions and admissible quadruples }
Let us be a little more precise about the data involved in this work. We will study for which globally anisotropic quadratic spaces $  \left (  q,V \right )  $ and $  \left (  q\p,V\p \right )  $ defining orthogonal groups $ H = O \left (  q,V \right )  $\label{notation:orthoLeft} and $  H\p = O\left (  q\p,V\p \right ) $\label{notation:orthoRight}, respectively, and for which non-trivial automorphic quadratic characters $ \xi \,\colon H_\fA \to \fC^\times $\label{notation:charLeft} and $ \xi\p \,\colon H_\fA\p \to \fC^\times $\label{notation:charRight} are the corresponding global theta lifts 
$ \pi= \gLift  $\label{notation:liftLeft} and $ \pi\p= \gLiftp $\label{notation:rightLeft} isomorphic. Moreover we will choose our data 
$ \adata $
  so that each of the representations $ \pi $ and $ \pip $ of $ G_{ \fA }  $ are irreducible, cuspidal, automorphic.
For the data considered in this work these representations satisfy first occurrence and in particular these theta lifts do not vanish.
\par
In order to narrow down our search for such admissible data
$ \alpha $, 
a natural starting point would be to consider the corresponding local question. So if we suppose that we have two irreducible, automorphic, cuspidal representations. See Theorem \ref{nonvanish} for the details.
$ \pi \cong \otimes \p_{ \nu } \pi_{ \nu } $ and $ \pip \cong \otimes \p_{ \nu } \pip_{ \nu }$ of $ G\ad $ then we consider necessary local conditions for the local representations 
$ \pi_{ \nu }  $ and $ \pip_{ \nu }  $ of $ G_{ \nu }  $ to be isomorphic.

Data 
$ \adatanu $
satisfying these local conditions will be call
a {\bbb{locally admissible quadruple}} and its global counterpart $ \adata $ will be called a \bbb{globally admissible quadruple}.
\label{notation:globallyAdmiQuad}
\begin{remark}
Note that the two local representations  $ \pi_{ \nu } $ and $ \pi_{ \nu } \p $ are the representations coming from the local theta correspondence. They are equal to a certain space of co-invariants of the Weil representation. That is $ \pi $ and $ \pi\p $ are given by
\label{notation:locThetaLiftLeft}\label{notation:locThetaLiftRight}
\begin{enumerate}
   \item  $ \pi_{ \nu } = 
 \left (  \omega_{ \psi_{ \nu } ,V_{ \nu }  } \otimes \xi_{ \nu } ^{ \vee }  \right ) _{ H_{ \nu }  } = 
    \left ( \omega_{ \psi_{ \nu } ,V_{ \nu }  }  \right )_{ H_\nu,\xi_\nu }  $ \\
   \item  $ \pi\p_{ \nu } = 
 \left (  \omega_{ \psi_{ \nu } ,V\p_{ \nu }  } \otimes \xi_{ \nu } ^{ \prime, \vee }  \right ) _{ H\p } = 
    \left ( \omega_{ \psi_{ \nu } ,V\p_{ \nu }  }  \right )_{ H\p_\nu,\xi\p_\nu }  $ \\
\end{enumerate}
See Section \ref{sec:localThetaAndHoweDuality}, page \pageref{sec:localThetaAndHoweDuality} for more details about the definition and properties of the local theta correspondence.
\end{remark}
\subsubsection{ Local necessary conditions }
We consider two necessary conditions for an isomorphism to exist between 
$ \pi_{ \nu } $ and $ \pi_{ \nu } \p $ :
\begin{enumerate}
   \item
The first local necessary condition is that
    $ \pi_{ \nu } $ and $\pi_{ \nu } \p $ must have the same central character. We will refer to this necessary condition as the \bbb{central character condition} and we will denote this condition by $ \CC   $. \\
   \item
The second local necessary condition is that the space of $  \left (  N_{ \nu } ,\psi_{ \nu }  \right )  -$ co-invariants  of the local representations $ \pi_{ \nu }  $ and $ \pi\p_{ \nu }  $  with respect to the unipotent radical of the Siegel parabolic in $ G $ are equal to zero simulataneously.
One can think of these $  \left (  N_{ \nu } ,\psi_{ \nu }  \right )  -$ co-invariants as local Fourier coefficients of the representations.
\par
 Namely we have the condition
$$  \left ( \pi_{ \nu } \right )  _{ N_{ \nu } ,\psi_\beta } = 0 \tn{ if and only if } \left ( \pi_{ \nu }\p \right ) _{ N_{ \nu } ,\psi_\beta } = 0 $$ for every 
$ \beta\in \Sym \left (  \F \right ) $, where $ \psi_{ \beta }  $ is a character of the unipotent radical $ N_{ \nu }  $ of the Siegel parabolic $ P_{ \nu }  \subseteq G_{ \nu }  $.
We will refer to this necessary condition as the \bbb{Fourier coefficient condition} and we will denote this condition by $  \FC $. \\
\end{enumerate}
\subsubsection{ An example of a locally admissible quadruple in the low rank setting }
For a locally admissible quadruple $ \adatanu $ we can explicitly define the characters $ \xi_{ \nu }  $ and $ \xi_{ \nu } \p $ in terms of local invariants of the two local quadratic spaces, namely in terms of the Hasse invariants,
$ h( q_{ \nu }  ) $ and $ h( q_{ \nu } \p  ) $,
 and discriminants,
$ d( q_{ \nu }  ) $ and $ d( q_{ \nu } \p  ) $,
of the local quadratic spaces
$ 
\left ( q_{ \nu } , V_{ \nu }  \right )
 $
and
$ 
\left ( q_{ \nu }\p , V_{ \nu } \p \right )
 $.
\par

For instance in the low rank case corresponding to Snitz's thesis we consider quadratic spaces of dimensions three and one. In this case we have the following description of the defining data $ \adatanu $.
\par
The character $ \xi_{ \nu }  $ of
the orthogonal group $ H_{ \nu }  $ corresponding to the three dimensional space $ \left ( q, V \right )  $
is given by parameters $  \lam \in \Fnu / \Fnu^{ \times, 2 }, \eps \in \mu_{ 2 } $
 and 
the character $ \xi\p_{ \nu }  $
of the orthogonal group $ H_{ \nu }\p  $ corresponding to the one dimensional space $ \left ( q\p, V\p \right )  $ is given by a parameter $ \epsp \in \mu_{ 2 } $.
\par

Then we have the following example in the low rank setting.

\begin{lemma}[Example of a locally admissible quadruple]
\label{sec:lowRankExampleOfAdmissibleQuadruples}

Suppose we are given a locally admissible quadruple 
$$
\adatanu
$$
Then the parameters
$ 
\lam \in \Fnu / \Fnu^{ \times, 2 }, \eps, \epsp \in \mu_{ 2 }
 $
are given by
  $$ 
\begin{cases}
 \lam=- d( q_{ \nu }   ) \cdot d( q_{ \nu } \p  )  \\
       \eps = ( -1 ,  d( q_{ \nu }   )  ) \cdot h( q_{ \nu } \p  )    \\
       \epsp = ( -1 , -d(q_{ \nu } \p ) ) \cdot  h( q_{ \nu } \p  )     \\
\end{cases}
 $$

\end{lemma}

\begin{remark}
The description in the above example is consistent with the data in Snitz's thesis. However it is a slight generalization of his work since Snitz used the quaternions to a large extent in his thesis and this restricted his results (albeit a minor restriction) to a three dimensional quadratic space with discriminant $  d \left ( q_\nu \right ) $ equal to one. 

See Appendix \ref{sec:SnitzsWork} on page \pageref{sec:SnitzsWork}  for more details about Snitz's thesis.
\end{remark}
\begin{remark}
For a similar description of the defining data of our characters in the higher rank odd case see

Theorem \ref{theDataDescription}.

\end{remark}
\subsection{ Determining isomorphisms via global and local pairings }
In this section we define the global and local pairings. We explain the transition from isomorphisms to pairings. 
We describe the properties of the representations $ \pi $ and $ \pi\p $ which are proved via Rallis' theory of towers where one needs global Fourier analysis with respect to the locally compact abelian group $ N \subseteq G $, i.e. with respect to the unipotent radical of the Siegel parabolic in $ G $ in order to prove what is referred to as "first occurrence" in Rallis' theory of towers and this in turn will prove that our global representations are cuspidal. Moreover the Fourier expansion of the theta integrals will also be used in the unfolding of a certain auxiliary global pairing. Next we move on to the problem of unfolding the global pairing. It turns out that this is too difficult to compute directly. Instead we construct a family of pairings given by an integral over $ \aq G $ of the two theta integrals against a relatively simple Eisenstein series. Thus we end up with a Rankin-Selberg type integral which is much easier to unfold and at a later point we take a certain residue in order to obtain formulas for the original global pairing of interest.
After unfolding the global pairing we obtain a product of local integrals.
We compute these integrals at almost every place where the data satisfies certain generic conditions. This generic computation is referred to as the unramified computation. 
Finally in Appendix \ref{sec:CAPRamanujan} we explain the problem we are addressing with respect to other interesting problems, such as the Rallis inner product formula, the Siegel-Weil formula and CAP representations.

\subsubsection{ The Global Pairings }
In this section we define the natural global pairing $ B \,\colon \schad  \times \schadp  \to \fC $ of interest which turns out to be impossible to unfold. In addition we define a family of global pairings $ B_{ s }  $ depending on a complex parameter $ s \in \fC $. It turns out that unfolding $ B_{ s }  $ is relatively easy. By taking a certain residue of the pairing $ B_{ s }  $ we obtain a formula for the pairing of interest $ B $.
\par
We define a natural global pairing 
\label{notation:globalNaturalPairing}
\begin{equation}
B \,\colon \schad  \times \schadp  \to \fC
\end{equation}

given by the Petersson inner product of the two integrals coming from the two global theta lifts. This is given by an inner product
$$
B \left ( \phi, \phip \right ) = 
\left < \thlift {   }
,
\thliftp{   }
\right >
_{ \bktg   } 
$$
that is 
\begin{equation}
B \left ( \phi, \phip \right ) = 
\int_{ 
\bktg 
 } 
\thlift {   } \left (  g \right ) 
\cdot
\overline{ \thliftp{   } \left (  g \right )  }
\,dg
\end{equation}

The functions appearing in the integrand will be precisely defined in Appendix \ref{sec:TransferringAutomorphicForms} on page \pageref{sec:TransferringAutomorphicForms}. The only thing one needs to know for now is that $ \thlift {   } $ and $ \thliftp {   } $ are vectors in $ \pi $ and $ \pi\p $, respectively,  coming from the global theta correspondence. 
We will also call the vectors $ \thlift {   } $ and $ \thliftp {   } $ \bbb{theta integrals}.
\begin{remark}
Note that the integral defining the pairing $ B $ convergences since the global theta lifts $ \pi $ and $ \pi\p $ are cuspidal and the theta integrals are elements of these spaces. The proof of cuspidality is given in Section \ref{sec:proofOfCuspidality}.
\end{remark}
\par
In theory it would be desirable to factor the pairing $ B $ into a product of local pairings. However this turns out to be impossible. Therefore in addition to the above pairing we consider a family of pairings $ B_{ s } $ whose unfolding is tractable and then $ B $ can be viewed as a limit or more precisely a residue of the family of pairings 
$$
B_{ s }  \,\colon \schad  \times \schadp  \to \fC
$$
which is similar to the above pairing but in addition we take the two theta integrals against a certain simple Eisenstein series $ \cE^{ \star }  = \cE^{ \star } _{ s }  $, i.e. 
\begin{equation}
\apair
\left ( \phi, \phip \right )
=
\int_{ \aq G } 
\thLift \phi g
\overline{ \thLiftp \phip g }
\naeis g s
\,dg
\end{equation}

The actual formula is of marginal significance. Our main concern is that unfolding $ B_{ s }  $ is relatively straightforward and the residue of $ \cE^{ \star }  $ at a certain point $ s_{ 0 }  $ is a constant function which after a normalization is identically equal to one.
\subsubsection{ The local pairing }
\label{sec:localPairing}
In this section we present the actual explicit formulas for the local pairing of interest.
One must keep in mind that the actual formula for this local pairing is of marginal significance. The important thing to understand is that this pairing is equivariant with respect to the diagonal action of $ G $ and the actions of $ H $ and $ H\p $ and it is a pairing between two functions which are elements of $ \pi_{ \nu } $ and $ \pi\p_{ \nu } $.
Moreover these local pairings will depend on a choice of certain vectors $ \vv $ and $ \vv\p $. This is an important point to keep in mind.
\par
Before we define the local pairing we must set some notations.
Let
$ P  $ be the Siegel parabolic of $ G $ and $ M $ be the Levi subgroup in the Siegel parabolic $ P $. The Levi subgroup $ M \subseteq P \subseteq G $ of the Siegel parabolic $ P $ is given by
$$
M = \Set*{  m \left ( a \right )  = \begin{pmatrix} a & 0 \\ 0 & a^\star \end{pmatrix} } { a \in GL_{ n } }
$$

Fix a maximal compact subgroup $ K \subseteq G\ad $  so that we have an Iwasawa decomposition 
$$
G\ad = P\ad K
 = M\ad N\ad K
$$
 and for 
$$
h = n\cdot m \left ( a \right ) k
,\qquad
n \in N\ad,\,\,a \in GL_{ n ,\fA} ,\,\,k \in K
$$
Let 
$$
\abs{  a \left ( g \right )  }  = 
\abs{  \det \left ( a \right )  } 
$$ 
and for
$$
\rhon = \dfrac{ n+1 } { 2 }
$$
let the section
$$
\Psi_{ s } \,\colon G\ad \to \fC
$$
be given by
$$
\Psi_{ s } \left (  g \right )  = 
\abs{  a \left ( g \right )  } ^{ s + \rhon } 
$$
Likewise, locally we can define the character
$$
\Psi_{ s,\nu } \left (  g \right )  = 
\abs{  a \left ( g \right )  }_{ \nu }  ^{ s + \rhon } 
$$
\begin{remark}
Usually we will not write the place $ \nu $ since it should be clear from the context whether we are working locally or globally.
Thus when working locally we will continue to denote $ \Psi_{ s } $ instead of $ \Psi_{ s,\nu } $.
\end{remark}

\par
Before we present the definition of our local pairing we must present a few more definitions.
We make certain choices of vectors $ \vv $ and $ \vv\p $ which our local pairing will depend upon. Namely let $ \vv \in  \left (  \vox \right )   _{ \F }  \cong V_{ \F } ^{ n }
\cong
\homj
,
\,
\vv\p \in  \left (  \voxp \right )   _{ \F }  \cong V_{ \F } ^{\prime, n }
\cong
\homjp
 $
 and
 $$
\gram\vv  =  \gramp{ \vv\p  }
$$
i.e. the Gram matrices of $ \vv $ and $ \vv\p $ are equal. 
\par
 We will consider vectors $ W_{ \phi } $ and $ W\p_{ \phip } $ coming from the local theta correspondence.
These vectors are given by
\begin{equation}
W_{ \phi }  \left ( g \right ) = 
\HintegralF
\left (
\wfunc  \left ( g, h \right ) 
\phi \right )
\left (  
\vv
\right ) 
\xi \left ( h \right ) 
\,dh
\end{equation}

and
\begin{equation}
W\p_{ \phip }  \left ( g \right ) = 
\HintegralFp
\left (
\wfuncpi \left ( g, h \p  \right ) 
\phip \right )
\left (  
\vv\p
\right ) 
\xi\p \left ( h \p \right ) 
\,dh\p
\end{equation}

\begin{defi}[The local pairing]
\label{sec:defLocalPairingTop}
\label{notation:localNaturalPairingTop}
We define a local pairing
$$
B_{ \nu,s } \,\colon \schloc  \times \schlocp   \to \fC 
$$
given by
\begin{align*}
B_{ \nu,s } \left ( \phi, \phip \right )  & =  
 \int_{ K_{ \nu }  } \int_{ M_{ \nu }  }  W_{ \phi }  \left (  mk \right )  W\p_{ \phip }  \left (  mk \right ) \delb m \fPsi{ mk }{ s } \,dm \,dk
\end{align*}
\begin{remark}
It is important to note that the pairing is between Schwartz spaces and not on the space of coinvariants defining the local representations.
\end{remark}

\end{defi}
\begin{remark}
This local pairing naturally appears from unfolding the family of auxiliary global pairings $ B_{ s }  \,\colon \schad  \times \schadp  \to \fC  $.
\end{remark}

\subsubsection{ From isomorphisms to pairings }
We would like to consider two local theta lifts which are given by coinvariant spaces. These spaces are of interest but difficult to describe.
Since in general quotient spaces are difficult to construct.
In Snitz's thesis, Snitz attempted to construct an isomorphism between such spaces. 
\par
The main idea in this section is that it is easier to understand pairings between such coinvariant spaces. And this can be replaced with an even simpler problem of finding pairings satisfying certain equivariance properties between certain Schwartz spaces. Schwartz spaces are already spaces that are well understood and not given by a quotient.
\par
The bottom line is that it turns out that a more accessible problem is finding non-zero pairings instead of constructing isomorphisms. In this section we sketch the transition from isomorphisms to pairings.
See Section \ref{sec:morphismsToPairings} for more details.
In what follows we only give the general idea of the transition.
\par
Upon closer reading of Snitz's work we see that he constructed  a local morphism between two local theta lifts
$$
\Delta\,\colon \co{ \psi_{ \nu }   }{ 3 }  \to \cop{ \psi_{ \nu }  }{ 1 } 
$$
where 
$ \co{ \psi_{ \nu }   }{ 3 }  $ and $ \cop{ \psi_{ \nu }  }{ 1 } $ denote the local theta lifts coming from our local quadratic spaces and local characters.
In general constructing such morphisms is highly non-trivial. Instead we construct a pairing
\label{notation:CoinvariantPairing}
$$
\widetilde{ B_{ \nu }  } \,\colon  \co{ \psi _{ \nu }  }{ 3 } \times \cop{ \psi_{ \nu } ^{ -1 }   }{ 1}  \to \fC
$$
However the spaces $ \co{ \psi_{ \nu }   }{ 3 }  $ and $ \cop{ \psi_{ \nu }  }{ 1 } $ are still relatively difficult to understand.
Instead we constructing a pairing on the level of Schwartz spaces which we denote by  $ B_{ \nu }  $

$$
B_{ \nu }  \,\colon \schloc  \times 
\schlocp \to \fC
$$
which satisfies certain equivariance properties relative to the actions of the groups $ H_{ \nu } ,H\p_{ \nu } ,G_{ \nu }  $.
\subsection{ Properties of the representations $ \pi $ and $ \pi\p $ }
In this section, we state the properties of the representations $ \pi $ and $ \pi\p $ of $ G\ad $. We state the result that these representations are irreducible, automorphic, cuspidal representations of the group $ G\ad $.
\par
The defining data $ \adata $ of the representations $ \pi $ and $ \pi\p $ will be defined so that both of these representation will be irreducible, automorphic, cuspidal representations.
Cuspidality will be an important technical point when computing residues of a certain normalized auxiliary pairing 
$ B_{ s }  
 \,\colon \schad  \times \schadp  \to \fC
$
.
Cuspidality will follow from Theorem \ref{sec:IrrCuspidalAutomorphic}.
\par
In order to prove that
the global theta lifts $ \pi= \gLift $ and $ \pi\p= \gLiftp $ are irreducible, automorphic, cuspidal representations of $ G\ad $ we use Rallis' theory of towers and the notion of first occurrence. Roughly speaking we can consider a sequence, i.e. a tower, of reductive dual pairs $ H_{ \fA } \times G_{ k, \fA } ,\,\, k\geq 1 $ where $ G_{ k, \fA } $ are a sequence of symplectic (or metaplectic) groups of rank $ k $ and we fix the non-trivial automorphic quadratic character $ \xi \,\colon H_\fA \to \fC^\times $ of the orthogonal group $ H_{ \fA }  $. This gives a sequence of automorphic representations, coming from the global theta correspondence, of the group $ G_{ k,\fA }  $.
\label{notation:thetaLiftTower}
$$
\pi_{ k } 
= 
\Theta^{ \left ( k \right ) } 
 \left ( q,V,\xi \right )
$$
Many of these representations might be zero. By the work of Rallis it is known that at some point one of these representations will be non-zero. The first lift that is nonzero is called the \bbb{first occurrence} and such a representation is known to be cuspidal.
Moreover there is a definite upper bound on the first $ k $ such that $ \pi_{ k } \neq 0 $, for instance if $  \left (  q,V \right )  $ is globally anisotropic and $ \dim\,V = n $ then first occurrence will take place for some $ 1\leq k\leq 2n $, i.e. the first $ k \in \fN  $ such that $ \pi_{ k } \neq 0 $ will be at most $ k = 2n $.
The above is referred to as \bbb{Rallis' tower property}. We have briefly described this property for our setting where the quadratic spaces are assumed to be anisotropic.
For more details on the more general setting see Rallis' work \cite{rallis1984howe}.
\par
For an admissible quadruple $ \adata $ and under the simplifying assumptions that the characters $ \xi $ and $ \xi\p $ are non-trivial and the quadratic spaces $ \left ( q, V \right )  $ and $ \left ( q\p, V\p \right )  $
are
anisotropic
our representations $ \pi $ and $ \pi\p $ of $ G_{ \fA }  $ will be constructed so that they both satisfy first occurrence and lift to the same symplectic (or metaplectic) group.
We will show this by some standard Fourier analysis on the unipotent radical of the Siegel parabolic of $ G $. Using Fourier analysis we will prove that each of our automorphic representations
$ \pi $ and $ \pi\p $ 
satisfy first occurrence. For more details on Rallis' theory of towers see Appendix
 \ref{sec:RallisTheoryOfTowers} page~\pageref{sec:RallisTheoryOfTowers}. 
The proof of cuspidality can be found in Section \ref{sec:zeroLifts} page~\pageref{sec:zeroLifts} and the Fourier analysis used in this proof and also used in the unfolding of the global auxiliary pairings can be found in Lemma \ref{sec:globalFourierExpansion} page~\pageref{sec:globalFourierExpansion} and Lemma \ref{sec:FourierExpLemma} page~\pageref{sec:FourierExpLemma}.
\subsection{ Global Fourier analysis }
Let $ M \subseteq P $ be the Levi subgroup of the Siegel parabolic $ P \subseteq G $ and let $ N \subseteq P $ be the unipotent radical in $ P $, i.e. $ P = MN $.
Fourier analysis is a key tool in proving first occurrence and also in the unfolding of a certain family of auxiliary global pairings
$$
B_{ s }  \,\colon \schad  \times \schadp  \to \fC
,\qquad s\in \fC
$$
 We present the notations and main results of the global Fourier analysis here.
\par
We define the \bbb{theta integral} by
\label{notation:thetaIntegral}
\begin{equation}
\thLift \phi g =
\int_{ \aq H }
\sum_{ \gamma \in  \left (  \vox \right ) _{ \F } } 
 \left ( 
\omega_{ \vpsi } 
 \left ( g, h \right ) 
\phi
  \right ) 
 \left (  \gamma \right ) 
\xi \left ( h \right )  
\,dh
\end{equation}

where $ X $ is a maximal isotropic subspace of the symplectic space $ W $, $ g\in G\ad $, and $ \phi\in \schad $.
\label{notation:thetaIntegralFourier}

\begin{defi}[Fourier coefficient]
\label{sec:defFourierCoefficient}
The  $ \psi_{ \beta }  $-Fourier coefficient of  $ \fthLift \phi $ is a function
$$
W_{ \phi, \beta } \,\colon  S \left ( V \right )  \to Fn \left ( \bktg  \right ) 
$$
given by 
$$
W_{ \phi, \beta } \left (  g \right ) = 
\thLiftBeta \phi g = 
\int_{ \qFA N N } 
\thLift \phi { ng }
\psi^{ -1 }_{ \beta }  \left ( n \right ) 
\,dn
$$
where  $ \beta\in  \SymF  $ and $ \phi \in S\voxpr_{ \fA } $ and 
$ Fn \left ( \bktg  \right )  $
denotes the space of functions on $ \aq G $.

\end{defi}
Here are some notations that will be needed.
Suppose we are given a quadratic space $  \left (  q,V \right )  $ over some field $ \F $.
Consider a vector $ \vv =  \left (  v_{ 1 } , \ldots, v_{ n }  \right ) \in  V ^{ n } \cong
\homj
$
we denote the \bbb{Gram matrix} of $ q $ by
$$ \gram\vv
=  
 \left (  b_{ q } \left (   v_{ i } , v_{ j } 
 \right ) \right )_{ i,j } 
$$
where  $ b_{ q } $ is the bilinear form associated to the quadratic form $ q $, i.e. 
$$
b_{ q }  \left ( v, w \right ) = 
\dfrac{ 1 } { 2 }
\left [ 
q \left (  v + w \right ) 
-
 q \left ( v \right ) 
-
 q \left ( w \right ) 
\right ]
$$

We will say that  \bbb{$ q $ represents $ \beta\in  \SymF  $} if and only if there exists $ \vvbeta \in 
 V_{ \F } ^{ n }  $ such that
$$
\gram\vvbeta = \beta
$$
\begin{lemma}[General Fourier coefficient of theta integral]
\label{sec:globalFourierExpansionOverview}
\begin{enumerate}
\item
Suppose that the quadratic space $ \left ( q, V \right )  $ represents a non-degenerate $ \beta\in  \SymF  $.
Then the  $ \psi_{ \beta }  $-th Fourier coefficient of the theta integral $ \fthLift \phi  $ is given by
$$
\thLiftBeta \phi { \nn g }
= 
 \int_{ \nicequotm } \left (  \wni ( g, h ) \phi  \right ) (\vvbeta) \xi  (h)\,dh
 \int_{ \badquotm }  \xi  (h_{ 0 } )\,dh_{ 0 } 
$$
where $ \vvbeta \in  \left (  \vox \right )   \left ( \F \right ) $ represents  $ \beta $ and $ g\in G_{ \fA }, n \left ( b \right )  = { \small\begin{pmatrix} I & b \\ 0 & I \end{pmatrix}   }\in N\ad ,b  \in \SymA $.
\item
If $ q $ does not represent $ \beta $ or if $ \beta $ is degenerate then $ \fthLiftBeta{ \phi } \equiv 0 $.
\end{enumerate}
\end{lemma}
Let us fix choices of vectors $ \vv\in \voxpr_{ \F }  $ 
and 
$  \vv\p\in \voxprp_{ \F } $
where $ \restr { \xi } { H^{ \vv } } \equiv 1 $ and the Gram matrix of $ \vv $ is invertible.
Similarly, we are assuming that $ \restr { \xi } { H^{ \vv\p } } \equiv 1 $ and the Gram matrix of $ \vv\p $ is invertible.
Under these assumptions we define the Fourier function.
\begin{defi}[Fourier function]
\label{sec:defFourierFunction}
\begin{enumerate}
\item
The Fourier function corresponding to the form
 $ \fthLift \phi $ is a function of the form
 $$
W_{ \phi }  \left ( g \right ) = 
\Hintegral
\left (
\wfunc  \left ( g, h \right ) 
\phi \right )
\left (  
\vv
\right ) 
\xi \left ( h \right ) 
\,dh
$$
where $ \phi \in S\voxpr_{ \fA }.
$
\item
For every $ \phi \in S\voxpr_{ \nu } $ define the local Fourier function by 
$$
W_{ \phi, \nu }  \left ( g \right ) = 
\Hintegralnu
\left (
\wfunc  \left ( g, h \right ) 
\phi \right )
\left (  
\vv
\right ) 
\xi \left ( h \right ) 
\,dh
$$
where $ \phip \in S\voxprp_{ \fA }.
$
\end{enumerate}
\end{defi}

Next we write down the Fourier expansion for the global theta integrals considered in this work.
Recall that $ M \subseteq P $ denotes the Levi subgroup of the Siegel parabolic $ P \subseteq G $
\begin{lemma}(Fourier expansion of our theta integral)
\label{sec:lemmaFourierExpansionOfThetaIntegral}
Given data $ \qcdata  
$
such that the corresponding global theta lift 
$ \pi = \gLift $ satisfies first occurrence then
\begin{align*}
   \thLift \phi g  
	 &   =  \sum_{ \alpha \in M_\F } 
	  W_{ \phi } \left ( \alpha g \right )
\end{align*}
where 
$ \phi \in \schad $ and $ g \in G\ad $.
\end{lemma}
\begin{remark}
Similar results hold for $ \thliftp {   } $.
\end{remark}
\subsection{ Unfolding the auxiliary global pairing }
Next we want to unfold our global pairing $ B_{ \nu }  \,\colon \schad \times 
\schadp \to \fC $. Unfortunately it is impossible to do so directly. Instead we consider a family of global pairings 
$ B_{ s }  \,\colon \schad  \times \schadp  \to \fC $
where we take our two theta  integrals and multiply them by what we refer to as the \bbb{normalized auxiliary Eisenstein series} $ \cE^{ \star } $.
\par
This Eisenstein series $ \cE^{ \star } $ will be a normalization of a very simple Eisenstein series $ \cE $ such that the residue of the normalized auxiliary Eisenstein series $ \cE^{ \star } $ at $ s = \rhon = \frac{ n+1 } { 2 } $ is equal to $ 1 $.
\begin{remark}
Note that a crucial point for the validity of our argument involving residues is the fact that the representations $ \pi $ and $ \pi\p $ which we will construct are cuspidal.
Actually it is enough that only one of these are cuspidal.
\end{remark}
The main result of interest in this section is as follows.
\begin{theorem}[Unfolding of the global auxiliary pairing]
For  $  Re \left ( s \right ) > \rhon $ and 
$ \phi\in \schad
$ and $
\phip\in \schadp
 $
we have
$$
\glblZetaConst^{ -1 } 
\apair
\left ( \phi, \phip \right ) 
= 
\finalGlobalExpansionAux
$$
\end{theorem}

\begin{cor}[Factorization of pairing]
\label{sec:corTakingResidue}
For $ \phi\in \schad $ and $
\phip\in \schadp
 $ we have
$$
 B \left ( \phi, \phip \right ) = 
 Res_{ s=\rhon } \apair \left ( \phi, \phip \right ) 
$$

\end{cor}

\subsubsection{ Definition of the auxiliary Eisenstein series for the curious reader }
In this part we define the auxiliary Eisenstein series used to make the unfolding of the global pairing feasible. Note that this part can be skipped without hindering the understanding of the main results.
\par
Recall that the natural global pairing is given by 
$$
B \,\colon \schad  \times \schadp  \to \fC
$$
where
$$
B \left ( \phi, \phip \right ) = 
\int_{ 
\aq G
 } 
\thlift {   } \left (  g \right ) 
\cdot
\overline{ \thliftp{   } \left (  g \right )  }
\,dg
$$
We define an auxiliary global pairing 
\label{notation:auxPairing}
$$
\apair \,\colon \schad \times \schadp \to \fC 
$$
given by 
$$
\apair
\left ( \phi, \phip \right )
=
\int_{ \aq G } 
\thLift \phi g
\overline{ \thLiftp \phip g }
\naeis g s
\,dg
$$
i.e. we integrate the theta integrals against a certain Eisenstein series described as follows.
First we define a simple Eisenstein series which will be used to facilitate in the unfolding.
\label{notation:eisSeries}
\begin{defi}[The Auxiliary Eisenstein series]
$$
\aeis g s = 
\sum_{ \gamma \in \rquot P G } 
\fPsi{ \gamma g }{ s }
$$
where 
$
\fPsi{ g }{ s }
= 
\absA{  a \left ( g \right ) }^{ s + \rhon }
$ where $ \rho_{ n } = \dfrac{ n+1 } { 2 } $.

\end{defi}
\begin{remark}
We have already precisely defined $ \Psi_{ s } $ and $  a \left ( g \right )  $ in Section \ref{sec:localPairing} page~\pageref{sec:localPairing}.
\end{remark}
Next we define the normalized auxiliary Eisenstein series.
The crucial fact about the auxiliary Eisenstein series is that it has a simple pole at $ s= \rhon $ whose residue is a constant function, i.e. its value is independent of $ g \in G $.
\begin{defi}[Normalized Auxiliary Eisenstein series]
\label{notation:auxEisSeries}
Define
$$
\naeis
g s
= 
\glblZetaConst
\cdot
\aeis g s
$$
which is normalized such that $ \cE^{ \star }  $ has a residue of  $ 1 $ at the simple pole $ s_{ 0 }  = \rhon = \frac{ n + 1 } { 2 }
 $
where
$ \glblZetaConst $
is a certain product of $ L $-functions. 
$ b_{ n } ^{ S }  $
 is the 
 following function
$$
\glblZetaConst = 
\prod_{ v\not\in S } 
b_{ n,\nu }  \left (  s \right ) 
$$
and 
$$
b_{ n,\nu }  \left (  s \right ) = 
L_{ v }  \left ( s + \rhon \right ) 
\cdot
\prod
_{ k=1 } 
^{ \left [ \frac{ n } { 2 } \right ] } 
L_{ v }  \left ( 2s + n + 1 - 2k \right ) 
$$
\end{defi}

\begin{remark}
The factor $ \glblZetaConst $ is of marginal significance.
The only important fact is that such a normalization exists and with this normalization the residue of the normalized auxiliary Eisenstein series is a constant function identically equal to one.
Explicit formulae for $ \glblZetaConst $ together with properties of this factor can be found in \cite{kudla1990poles}.
\end{remark}

\subsection{ Statement of the unramified computation }
In the section we state the unramified computation.
\par
The unfolding of the auxiliary global pairing $ B_{ s } $ naturally leads to the definition of the local pairing
$$
\localPairingSigs
$$
given by
\begin{align*}
B_{ \nu } \left ( \phi, \phip \right )  & =  
 \int_{ K} \int_{ \ga M }  W_{ \phi }  \left (  gk \right )  W\p_{ \phip }  \left (  gk \right ) \delb m \fPsi{ mk }{ s } \,dm \,dk
\end{align*}

\begin{remark}
An important point to note is that the above pairings $ B_{ \nu }  $ and $ B_{ \nu, s }  $ depend on choices of vectors $ \vv $ and $ \vv\p $.
\end{remark}
In order to have $ B= \otimes_{ \nu }  B_{ \nu }  $ we still need to deal with the locally generic set up of the local pairing. This is referred to as the unramified computation.
In the unramified setting we assume that our characters, quadratic spaces and lattices in the corresponding vector spaces are unramified.
It is important to stress that we have only treated the unramified computation in the case $ n = 1 $.
See Section \ref{sec:unrSetting} page~\pageref{sec:unrSetting}
for more details on the unramified computation.

\begin{lemma}[Unramified computation]
\label{sec:lemmaUnrComputation}
In the unramified setting for $ n = 1 $
$$
 B_{ \nu, s } 
\left ( \phi, \phip \right )
 = 
\left ( 1-q_{ \nu } ^{ -s } \right )^{ -1 } 
$$

\end{lemma}
\begin{remark}
This is the only point in this work which is not stated in full generality.
In order to generalize all of our results to higher rank groups one must complete the unramified computation for general $ n \in \fN $.
Note that although we have not proven that the local factors are Eulerian for \( n>1 \) we do believe that such is the case.

\end{remark}
\subsection{ Idea of the proof of the main statement }
In this section we sketch the proof of the main statements Theorems \ref{sec:mainStatement} and \ref{sec:Alternative}. Essentially Theorem \ref{sec:mainStatement} is the main result and Theorem \ref{sec:Alternative} is the way we go about proving the main result in terms of showing that the global pairing $ B $ is a non-degenerate pairing.

\begin{remark}
Note that since we have not completed the unramified computation then technically speaking this result is valid only for $ n = 1 $. However every other aspect of the proof, including the unfolding of the global auxiliary pairing, taking of residues and proof of irreducibility and cuspidality of our global automorphic representations $ \pi $ and $ \pi\p $ of $ G\ad $ are valid in complete generality.
\end{remark}
\begin{customthm}{\ref{sec:mainStatement}}[{Main statement}]

Given a globally admissible quadruple $ \adata $.
If a certain global partial $ L $-function $ \Delta_{ S }  $
has a pole at
$ s = \rhon := 
\frac{ n+1 } { 2 }
 $
then the global theta lifts 
$$
\pi = \gLift\tn{ and }
\pi\p = \gLiftp
$$
are isomorphic.
\end{customthm}
\proof
[Idea of proof]
This will follow immediately from Theorem \ref{sec:Alternative}.

$ \square $

\begin{customthm}{\ref{sec:Alternative}}[{Alternative main statement in terms of pairings}]
Suppose $ \adata $ is a globally admissible quadruple.
Then
\begin{enumerate}
\item
There exists a natural global pairing 
$$
B \,\colon \schad  \times \schadp  \to \fC
$$
which is Eulerian (i.e. factors as a product of local pairings) and satisfies for every factorizable
$ 
\phi = \otimes \loc ^{ \p } \phi\loc \in \schad  
$ and $
\phip = \otimes \loc ^{ \p } \phi\p\loc \in \schadp
 $
 the equality
 \begin{align*}
    B \left ( \phi, \phip \right )  & =
\kappa
\cdot
    Res_{ s = \rhon }
\Delta_S \left ( s \right ) 
\cdot
\prod_{ \nu \in S } 
 B\loc 
\left ( \phi_{ \nu } , \phip _{ \nu } \right ) 
 \end{align*}
where $ \rhon = \dfrac{ n+1 } { 2 } $, $  \Delta_S $ is a certain global partial $ L- $function, $ S $ is a finite set of places, $ B\loc  $ are certain explicit local pairings, and $ \kappa $ is a non-zero constant coming from the residue of a certain simple normalized auxiliary Eisenstein series $ \cE^{ \star } $.
\item
If the global partial $ L $-function, $ \Delta_S \left ( s \right ) $, has a pole at $ s = \rhon $
then
 $ \pi $ and $ \pi\p $ are isomorphic.
\end{enumerate}
\end{customthm}
\proof[Sketch of proof]
We will prove this in two steps.
\begin{enumerate}
\item
The factor 
$$
\prod_{ \nu \in S } 
 B_{ \nu } 
\left ( \phi_{ \nu } , \phip _{ \nu } \right ) 
$$
is a finite product coming from the unfolding of our global pairing.
\item
The factor 
$$
    Res_{ s = \rhon }
\Delta_S \left ( s \right ) 
$$
comes from the unramified computation, i.e. 
$$
\Delta_S \left ( s \right ) 
 = 
\prod_{ \nu \not \in S } 
\left ( 1-q_{ \nu } ^{ -s } \right )^{ -1 } 
$$
\end{enumerate}
$ \square $

\subsection{ Relation to other work }
It is interesting to note that our main result bears a close resemblance to the Rallis inner product formula. This is not surprising since our main object of study in the global setting is an inner product of two theta integrals as is in the Rallis inner product formula.
Amongst other things the Rallis inner product formula is an important tool in determining first occurrence of a theta lift and relating this to poles of certain $ L $-functions.
\par
Note that another interpretation of the Rallis inner product formula is to describe when an automorphism exists between a theta lift and itself, hence we are generalizing this result to the existence of non-zero morphisms between theta lifts in general. However one must keep in mind the following caveat of this work. We are lifting highly trivial and concrete automorphic representations of the orthogonal groups $ H_{ \fA }  $ and $ H_{ \fA } \p $. Giving a general conjecture and proof of some kind of generalized Rallis inner product formula for general cusp forms $ f^{ H }  $ and $ f^{ H\p }  $ of $ H_{ \fA }  $ and $ H_{ \fA } \p $ might be much more difficult.
\par
Also note that a common step in proofs of the Rallis inner product formula use the doubling method which usually depends on a suitable Siegel-Weil formula. 
Our proof avoids the use of the Siegel-Weil formula. Instead we consider a family of global pairings 
$ B_{ s }   \,\colon \schad  \times \schadp   \to \fC $
depending on a complex parameter $ s\in \fC $ which is a Rankin-Selberg convolution of our two theta integrals against a certain simple normalized auxiliary Eisenstein series $ \cE^{ \star }  $. The main steps of our proofs bear a strong resemblance to Jacquet's \cite{jacquet2006automorphic}
 $  \left ( 1972 \right )  $ proof of the $ GL_{ 2 } \times GL_{ 2 } $ case using the Rankin-Selberg method. In this sense our proof is very different from the standard proofs of the Rallis inner product formula. 
As an aside, it would be interesting to see if our methods could be applied to furnish an alternative proof of the Rallis inner product formula without recourse to the Siegel-Weil formula or whether our results could be applied to give an alternative proof of the Siegel-Weil formula. 
\par
Another interesting observation worth mentioning is that if we apply the formulas of the Weil representation in the definition of our local pairings then we obtain a function closely related to Godement-Jacquet's $ L $-function which is a generalization to $  GL \left ( n \right )  $ of the $ L $-function considered in Tate's thesis. Note that an alternative way of thinking of $  GL \left ( n \right )  $ is as the Levi part of the Siegel parabolic.
\par
Initially we made an educated guess based on Godement-Jacquet's $ L $-function to define our local pairing, then we realized a certain averaging operator would be a natural substitute of this pairing and finally we discovered the idea of the auxiliary pairing which was relatively easy to unfold and this led us yet again to a similar yet slightly different definition of our local pairing.
\par
This previously defined local natural pairing differs from the local pairing we are using in the sense that we do not integrate over the compact group and the modular character of the Siegel parabolic does not appear in the formula.
We will not examine this alternative local pairing although it could be interesting to see how these alternative local pairs are related.


\newpage
\newcommand{ \jii }{ J_{ 2n }  } 
\newcommand{ \sgp }{ Sp_{ 2n } } 
\newpage
\section{ Notations }

\subsection{ Measure }

Let  $ \F $\label{fm:f} be a number field and let  $ \fA = \fA_{ \F }  $\label{fm:adele ring} denote its \adele ring. For any place  $ \nu $ of  $ \F $, let  $ \Fnu $\label{fm:fun} be the completion of  $ \F $ with respect to $ \nu $. If $ \K $\label{fm:loccompact} is a locally compact, nonarchimedian field then we denote its ring of integers by $ \cO_{ \K }  $\label{fm:ringing} and the maximal ideal of  $ \cO_{ \K }  $ by  $ \cP_{ \K } = \cO_{ \K } \gen $\label{fm:macideal}
where  $ \gen \in \cP_{ \K } -\cP_{ \K } ^{ 2 }  $\label{fm:generator} is a generator for the maximal ideal  $ \cP_{ \K }  $. Let  \label{fm:cardresiduefield}$ q = \abs{ \cO_{ \K } / \cP_{ \K }  } $ denote the cardinality of the residue field. Denote by  $ \abs{ \cdot } _{ \K }  $\label{fm:absval} the absolute value of  $ \K $, normalized so that  $ \abs{ \gen } = q^{ -1 }  $. We fix a Haar measure  $ \dx $ on  $ \K $ such that 
$$
\int_{ \cO_\K } \dx = 1
$$
and let 
$$
\tn d^{ \times } x = 
\dfrac{ \dx } { \abs{ x } _{ \K }  }
$$
be the Haar measure on  $ \K^{ \times }  $.
\par
Note that if  $ \K = \Fnu $ then we will use the simpler notation  $ \cO_{ \nu }  =  \cO_{ \Fnu }  $ and  $ \cP_{ \nu }  =  \cP_{ \Fnu }  $.

\subsection{ The groups }
\subsubsection{ The orthogonal group }
For a quadratic space $  \left (  q,V \right )  $ one can construct a corresponding bilinear form $ b_q \,\colon V\times V \to \F $. The orthogonal group with respect to this form is defined to be 
$$
 O \left ( q,V \right )  = 
\Set*{  h \in GL_{ 2n } } {  b_q \left ( hv,hw \right )  =  b_q \left ( v,w \right ) ,\forall v,w \in V }
$$
In general in this work $ H $ and $ H\p $ will be used to denote the orthogonal group corresponding to the quadratic space 
$ \left ( q, V \right )  $ and $ \left ( q\p, V\p \right )  $, respectively.
In our definition of the orthogonal group the action is a left action.
\subsubsection{ The symplectic group }
We define the symplectic group which acts on the right of a symplectic vector space $  \left (  W_{ 2n } ,\left < \cdot, \cdot \right > \right )  $ as follows 
\begin{align*}
\sgp 
& = 
\Set*{ g\in GL_{ 2n }  }
{ 
\left < vg,wg  \right > = 
\left < v,w \right >
 ,\forall v,w \in V
  }
 \\
&   = 
\Set*{ g\in GL_{ 2n }  } { 
^{ t } g
\jii
g
= 
\jii
 }
\end{align*}
where 
$ J_{ 2n }  = \begin{pmatrix} 0_n & I_n \\ -I_n & 0_n \end{pmatrix}  $.
\par
In general in this work $ G $ will be used to denote the symplectic or metaplectic group.
\subsubsection{ The metaplectic group }
The definition of the metaplectic group (see Theorem \ref{sec:defWeil} page~\pageref{sec:defWeili}) is more involved and one must define the Weil representation first. 
Thus we devote a separate section for the definition. See Appendix \ref{sec:HeisStoneWeilRep} for the precise definition of the metaplectic group.
\subsubsection{ Subgroups of $ G $ }
\label{sec:ImportantSubgroups}
We set notations for the Siegel parabolic $ P $ of $ G $ and the Levi subgroup $ M $ and unipotent radical $ N $ of $ P $, that is $ P = MN $.
We also define an important character $ \Psi_{ s } $ which will be used later to define the auxiliary Eisenstein series.
By the Iwasawa decomposition 
$$
G = PK = MNK\cong GL_{ n } NK
$$
where $ P  $ is the Siegel parabolic of $ G $ and $ M $ is the Levi subgroup of $ P $. The group $ M $ is of the form 
$$
M = \Set*{  m \left ( a \right )  = \begin{pmatrix} a & 0 \\ 0 & a^\star \end{pmatrix} } { a \in GL_{ n } }
$$
and the group $ N $ is of the form 
$$
N = 
\Set*{  n \left ( b \right )  = \begin{pmatrix} I & b \\ 0 & I \end{pmatrix} } { b \in \Sym }
$$
where 
$$
\Sym = 
\Set*{ a \in GL_{ n }  } { a^t = a }
$$
We also define a function $ a \,\colon G \to GL_{ n } $ as follows. For $ g \in G $ we have $ g =  m \left ( a \right ) nk $ where
$ a \in GL_{ n } ,n \in N,k \in K $ then we define $  a \left ( g \right )  = a $.
\par
In addition 
$
\abs{  a \left ( g \right )  } 
$
is short hand for the absolute value of the determinant of $  a \left ( g \right ) $, namely 
$$
\abs{  a \left ( g \right )  }  = 
\abs{  \det \left ( a \left ( g \right ) \right )  } 
$$
Finally we denote
$$
\begin{cases}
   \Psi_{ s }  \,\colon GL_n(\fA) \to \fC^{ \times }  \\
   \fPsi{ g }{ s }=\absA{  a \left ( g \right ) }^{ s + \rhon } 
   \\
\end{cases}
$$
where $ \rhon = \dfrac{ n+1 } { 2 } $. 	
\begin{remark}
The section $ \Psi_{ s } $ will be used to define a certain auxiliary Eisenstein series $ \cE $ which will assist us in unfolding the natural global pairing $ B \,\colon \schad  \times \schadp  \to \fC 
$.
See Section \ref{sec:unfoldGlobal} page~\pageref{sec:unfoldGlobal} for more details on the global unfolding.
\end{remark}

\subsubsection{ Quadratic Spaces }
In this work we will be concerned with two quadratic spaces. These quadratic spaces will be denoted by
 
\begin{enumerate}
   \item   $  \left (  q ,V  \right )  $  an $ m $-dimensional quadratic space. \\
   \item   $  \left (  \qp ,\vp  \right )  $ an $ \mmp $-dimensional quadratic space. \\
\end{enumerate}
We also recall the definition of the Gram matrix and the notion of a quadratic space representing a value.
\begin{defi}[Gram matrix]
\label{sec:GramMatrix}
Suppose we are given a quadratic space $  \left (  q,V \right )  $ 
over some field $ \F $ .
Consider a vector $ \vv =  \left (  v_{ 1 } , \ldots, v_{ n }  \right ) \in  V ^{ n } 
 \cong
\homj
$
we denote the Gram matrix of $ q $ by
$$ \gram\vv
=  
 \left (  b_{ q } \left (   v_{ i } , v_{ j } 
 \right ) \right )_{ i,j } 
$$
where  $ b_{ q } $ is the bilinear form associated to the quadratic form $ q $, i.e. 
$$
b_{ q }  \left ( v, w \right ) = 
\dfrac{ 1 } { 2 }
\left [ 
q \left (  v + w \right ) 
-
 q \left ( v \right ) 
-
 q \left ( w \right ) 
\right ]
$$
\end{defi}
\begin{defi}
We will say that  $ q $ represents $ \beta\in  \SymF  $ if and only if there exists $ \vvbeta \in 
 V_{ \F } ^{ n }  $ such that
$$
\gram\vvbeta = \beta
$$
\end{defi}

\subsubsection{ The specific orthogonal groups of interest in this work }
 We will be concerned with the orthogonal groups corresponding to the quadratic spaces $ \left ( q, V \right )  $ and $ \left ( q\p, V\p \right ) $ which we will denote by 
\begin{enumerate}
   \item   $ H = O \left (  q, V  \right )  $ i.e. the orthogonal group on the "left-hand side". \\
   \item   $ \hp = O \left (  \qp, V\p  \right )  $ i.e. the orthogonal group on the "right-hand side". \\
\end{enumerate}
\subsubsection{ Quadratic Characters }
We will be concerned with certain automorphic characters of the orthogonal groups $ H\ad $ and $ H \ad \p $ throughout this work. Here are our notations for these characters.
\begin{enumerate}
   \item  Let $ \xi \,\colon H \ad \to \fC $ be a non-trivial quadratic automorphic character of the group  $ H \ad  $. \\
   \item  Let $ \xi\p \,\colon H\p \ad \to \fC  $ be a non-trivial quadratic automorphic character of the group  $ \hp \ad  $. \\
\end{enumerate}

\subsection{ 
Quadratic character of $ \F^{ \times }  $ }

Fix a non-trivial additive character $ \psi \,\colon \F \to \fC $. This character appears in the definition of the Weil representation.
Let
$ \gamma $ denotes the \bbb{Weil index} as described in Definition \ref{def:weilindex} page~\pageref{def:weilindex}.
\par
For every $ m $-dimensional quadratic space  $  \left (  q, V \right )  $ we will associate a quadratic character of $ \F^{ \times }  $ given in terms of the Hilbert symbol
$$
 \chi_{ V }  \left ( x \right ) = 
\left (  
 det \left ( x \right ) ,
 \left ( -1 \right )^{ \frac{ m \left ( m - 1 \right )  } { 2 } 
}
d \left ( V \right )
\right )_{ \F } 
$$
where $  d \left ( V \right )  $ is the discriminant of the quadratic space $ \left ( q, V \right )  $ and $ m = \dim \,V $.
Finally, we will associate another character of $ \F^{ \times } $ given by 
$$
\widetilde{ \chi }_{ V }  \left ( x \right ) = 
\begin{cases}
   \gammafactorx\chi_{ V }  \left ( x \right )  & \text{if} $ m $\tn{is odd } \\
   \chi_{ V }  \left ( x \right )  & \text{if} $ m $\tn{is even } \\
\end{cases}
$$
\subsection{ Coinvariants }
\label{sec:coinvariants}
For a representation $  \left (  \pi,V_{ \pi } ,G \right )  $ of a group $ G $ the space of $ G $-coinvariants is given by the quotient representation
$$
V_{ G }  = V /
\left < 
\pi \left ( g \right ) \cdot v
-
v
\mid
g \in G,
v \in V
 \right >
$$
We will use this to define the local theta correspondence in Appendix \ref{notation:LocalThetaLift}.

\subsection{ Additive characters }
In this section we present the notations of the additive characters appearing in this work.
Namely an additive character of $ \F $ and an additive character of the unipotent radical $ N $ in the Siegel parabolic subgroup $ P $ of $ G $.
\par
Fix a non-trivial additive character $ \psi \,\colon \F \to \fC $.
\subsubsection{ Additive characters of $ \F $ }
\label{addChar}
Given an additive character $ \psi \,\colon \F \to \fC $ then
for $ a \in \F^{ \times } $ we define another additive character of $ \F $ given by
$$
\psi_{ a }  \left (  x \right )  =  \psi \left ( ax \right ) 
$$
\begin{remark}
These describe all non-trivial additive characters of the additive group of $ \F $.
\end{remark}
\subsubsection{ Characters of the unipotent radical of $ G $}
\begin{lemma}[Characters of $ N $]
Characters of the unipotent radical of the Siegel parabolic $ N \subseteq G $ are characters
$$
\psi_\beta \,\colon N \to \fC
$$
given by 
$$
 \psi_\beta
\begin{pmatrix} I & x \\ 0 & I \end{pmatrix} 
 = 
\psi \left (  
tr \left (  \beta x \right ) 
 \right ) 
$$
where $ x \in \Sym \left ( \Fnu \right ),  \beta \in \Sym \left ( \F \right ) $ 

\end{lemma}
\proof
This is a standard result.
$ \square $

\subsection{ Notations for the Weil representation }
For more details about the Weil representation see Appendix \ref{sec:ScrodModel}.
\subsubsection{ Notation for the Weil representation restricted to a reductive dual pair }
For a fixed nontrivial additive character  $ \psi \,\colon \F \to \fC $  we denote the Weil  representation
by
$$
\omega_{ \vpsi } \defeq
\restr { \omega_{ \psi }  } { \widetilde{ Sp } \left (  W \right )  \left (  \fA \right ) \otimes O \left (  V \right ) \left (  \fA \right )   }
$$
\subsection{ Notation for the theta kernel }
We can form a natural equivariant function
$$
\theta_{ \vpsi } \colon \omega_{ \vpsi } \to
\cA \left (  G\ad \times H\ad \right ) 
$$
which we will refer to as the theta kernel. 
Moreover this function can be given explicitly by
$$
\theta_{ \vpsi } ^{ \phi }  \left (  g,h \right ) = 
\sum_{ \gamma \in  \left (  V \otimes X \right )_{ \F }  } 
 \left (  
\omega_{ \vpsi } \left (  g,h \right )
\phi
  \right ) 
 \left (  \gamma \right ) 
$$

\newcommand\numberthis{\addtocounter{equation}{1}\tag{\theequation}}
\newcommand{ \weilNoN }[1]{ \omega_{ \psi, #1 } }
\newpage

\section{ Admissible quadruples }
\label{sec:admissibleQuadrupleExamples}
\good
In this section we define the quadratic spaces and quadratic characters which will be of interest in this work. Such data will be called a \bbb{locally admissible quadruple}. The data will need to satisfy certain local necessary conditions for an isomorphism to possibly exist between the representations $ 
 \left (  \pi_{ \nu }, S  \left( \voxpr_{ \F_{ \nu }  } \right), G_{ \nu }   \right ) 
  $ and $  \left (  \pi_{ \nu }\p, S  \left( \voxprp_{ \F_{ \nu }  } \right), G_{ \nu }   \right )  $.
\par
We are interested in defining data $ \adatanu $ such that the two representations
 $ \pi_{ \nu }  $ and $ \pi_{ \nu } \p $ of the group $ G_{ \nu } $ will be isomorphic. A necessary condition is that they will have the same central character and certain coinvariants will be nonzero at the "same time".
 This is the motivation behind our definition of the $  \left ( CC \right )  $ central character condition and the  $  \left ( FC \right )  $ Fourier coefficient condition which together will give us the definition of a locally admissible quadruple $ \adatanu $.
\begin{remark}[Simplifying Assumptions]
Throughout this work we will be working under the simplifying assumption that the automorphic quadratic characters $ \xi \,\colon H_\fA \to \fC^\times $ and $ \xi\p \,\colon H_\fA\p \to \fC^\times $ are non-trivial and that the quadratic spaces $  \left (  q,V \right ) $ and $ \left (  q\p,V\p \right ) $ are globally anisotropic.
\end{remark}
\begin{defi}[QC-data]
\label{sec:QCData}\label{notation:QCData}
For a quadratic space $  \left ( q, V \right )  $ and a character  $ \xi $ of the orthogonal group $ O \left ( q, V \right )  $ we call the pair  $$ \qcdata $$ \bbb{QC-data}, i.e. the quadratic space and character pair data.
\end{defi}
\subsection{ The central character condition }
\begin{defi}[CC - Central character condition]
\label{fm:CC}
A quadruple
$$
\adatanu 
$$
is said to satisfy the \bbb{central character condition} if $ \pi_{  \nu }  $ and  $ \pi\p_{ \nu }  $ have the same central character.
\end{defi}
\begin{lemma}[Characterization of the central character condition]
The quadruple $ \adatanu $  satisfies the central character condition 
if and only if
$$ 
\chiv{ -I_{ n }  }
\xi_{ \nu }  ( -I_{ n }  )
 =
\chivp{ -I_{ n }  }
\xi_{ \nu } \p ( -I_{ n }  )
$$
\end{lemma}
\proof
We present a proof in the case that the dimensions of our quadratic spaces are odd. The proof is the same for the even case, however the gamma factor (which cancels anyways in the odd case) does not appear in the even case.
Recall that 
$$
\pi_{ \nu } =  \left (  \omega_{ \psivloc } \right )_{ H_{ \nu } ,\xi_{ \nu }  } 
$$
and
$$
\pi\p_{ \nu } =  \left (  \omega_{ \psivlocp } \right )_{ H_{ \nu } \p,\xi_{ \nu } \p } 
$$
\label{notation:centralCharLeft}\label{notation:centralCharRight}
Let $ z_{ \pi_\nu } $ denote the central character of $ \pi_{ \nu } $. Similarly, let $ z_{ \pip_\nu } $ denote the central character of $ \pip_{ \nu } $.
For  $ \phi_{ \nu } \in  \pi_{ \nu } $ we have 
$$
\omega_{ \psivloc } \left (  h \right ) \phi_{ \nu }  =  \xi_{ \nu }  \left ( h \right ) \phi_{ \nu } 
$$
therefore 
$$
 \left (  \omega_{ \psivloc } \left (  h \right ) \phi_{ \nu }  \right )  \left (  v \right ) =  \xi_{ \nu }  \left ( h \right )  \phi_{ \nu }  \left ( v \right ) 
$$
thus
$$
\phi_{ \nu }  \left ( h^{ -1 } v \right ) =  \xi_{ \nu }  \left ( h \right )  \phi_{ \nu }  \left ( v \right ) 
$$
In particular for  $ h = -I_{ n }  $ we have
\begin{align*}
   \phi_{ \nu }  \left ( -v \right )  & =  \xi_{ \nu }  \left ( -I_{ n }  \right ) \phi_{ \nu }  \left ( v \right )
\end{align*}
Consider the matrix  $  m \left ( -I_{ n }  \right ) = \begin{pmatrix} -I_{ n }  & 0 \\ 0 & -I_{ n }  \end{pmatrix} = -I_{ 2n } $ in the center of  $ G $  therefore
we have
$$
 \left ( \omega_{ \psivloc } \begin{pmatrix} -I_{ n }  & 0 \\ 0 & -I_{ n }  \end{pmatrix}\phi_{ \nu }  \right ) \left ( v \right ) 
 = 
 \left (  \pi_{ \nu } \begin{pmatrix} -I_{ n }  & 0 \\ 0 & -I_{ n }  \end{pmatrix}\phi_{ \nu }  \right ) \left ( v \right ) 
 = 
 z_{ \pi_\nu } \begin{pmatrix} -I_{ n }  & 0 \\ 0 & -I_{ n }  \end{pmatrix} \phi_{ \nu }  \left ( v \right ) 
$$
where  $ z_{ \pi_\nu }  $ is the central character of  $ \pi_{ \nu }  $.
On the other hand by the action of the Weil representation (see Lemma \ref{lemma:weilFormulas}) we have
\begin{align*}
     \left ( \omega_{ \psivloc } \begin{pmatrix} -I_{ n }  & 0 \\ 0 & -I_{ n }  \end{pmatrix}\phi_{ \nu }  \right ) \left ( v \right ) 
        & =   \gammafactor  \abs{ -1 } ^{ m/2 }  \chiv{ -I_{ n }  }  \phi_{ \nu }  \left ( -v \right )    \\
	 &  
= 
 \gammafactor \chiv{ -I_{ n }  }  \xi_{ \nu }  \left ( -I_{ n }  \right ) \phi_{ \nu }  \left ( v \right )    
\end{align*}
therefore
$$
z_{ \pi_\nu } \left ( m \left ( -I_{ n }  \right )  \right ) = 
\gammafactor \chiv{ -I_{ n }  }  \xi_{ \nu }  \left ( -I_{ n }  \right )
$$
Similarly we have
$$
z_{ \pip_\nu } \left ( m \left ( -I_{ n }  \right )  \right ) = 
\gammafactor \chivp{ -I_{ n }  }  \xi_{ \nu } \p \left ( -I_{ n }  \right )
$$
therefore the central character condition is equivalent to
$$
\chiv{ -I_{ n }  }
\xi_{ \nu }  ( -I_{ n }  )
 =
\chivp{ -I_{ n }  }
\xi_{ \nu } \p ( -I_{ n }  )
$$
$ \square $
\subsection{ The Fourier coefficient condition }
In this section we describe the Fourier coefficient necessary condition for local isomorphisms between $ \pi_{ \nu } $ and $ \pi\p_{ \nu } $ to exist.
\begin{remark}
Recall the notations
\begin{enumerate}
\item
$ \pi_{ \nu }  = \localPi $
\\
\item
$ \pi_{ \nu }\p  = \localPip $
\end{enumerate}
\end{remark}

\begin{defi}[Local Fourier Coefficients]
\label{notation:localFourierCoefficients}
Let $ N \subseteq G $ be the unipotent radical of the Siegel parabolic in $ G $ and let $ \psi_{ \beta }  $ be the non-trivial additive character of $ N $ determined by the symmetric matrix $ \beta $. Then we define the \bbb{local Fourier coefficient} 
of 
$ \pi_{ \nu }  = \localPi $
to be the co-invariant space
$
\localFourier
$.
\end{defi}

\begin{defi}[FC - Fourier coefficient condition]
A quadruple $ \adatanu $\label{fm:FC} is  said to satisfy the \bbb{Fourier coefficient condition} if and only if
$$
\localFourier = 0
\iff
\localFourierp = 0
$$
for every  $ \beta \in \Sym \left ( \F \right )
    $.
\end{defi}
\subsection{ Admissible quadruples }
In this section we describe the data of interest in this work. Namely data satisfying certain local necessary conditions in order for local isomorphisms between the representations $ \pi_{ \nu }  $ and $ \pi\p_{ \nu }  $ to exist.
\begin{defi}[Locally admissible quadruple]
A quadruple $ \adatanu $\label{fm:locallyadmquad} is called \bbb{locally admissible} if and only if it satisfies the central character $ \CC $ and the Fourier coefficient $ \FC $ conditions.

\end{defi}
\begin{lemma}[Example]
For the dual pairs  $ O \left ( q,V \right ) \times \widetilde{ SL_2 } $ and $ O \left ( q\p,V\p \right ) \times \widetilde{ SL_2 } $
a locally admissible quadruple 
 $  \adatanu $
is characterized by the data
  $$ 
\boxed{ \begin{cases}
 \lam=- d( q_{ \nu }   ) \cdot d( q_{ \nu } \p  )  \\
       \eps = ( -1 ,  d( q_{ \nu }   )  ) \cdot h( q_{ \nu } \p  )    \\
       \epsp = ( -1 , -d(q_{ \nu } \p ) ) \cdot  h( q_{ \nu } \p  )     \\
\end{cases} }
 $$
where  $ d $ denotes the discriminant and $ h $ denotes the Hasse invariant and  $  \lam \in \Fnu / \Fnu^{ \times, 2 }, \eps, \epsp \in \mu_{ 2 } $.
\end{lemma}
\proof
This follows directly from the $  \left ( CC \right )  $ and $  \left (  FC \right )  $ together with properties of the Hilbert symbol and the definition of the Hasse invariant.
$ \square $
\newcommand{ \dcase }[2]{  $ 
 dim \left ( V \right ) = #1,
 dim \left ( V\p \right ) = #2
  $ }

\begin{remark}
Note that if $ \adata $ is a globally admissible quadruple such that the corresponding theta lifts are non-zero then necessarily we have
 $ 
 dim \left ( V_\nu \right ) , dim \left ( V\p_\nu \right ) \in  \left \{ n, n + 2 \right \} $.
See Lemma \ref{sec:dimRestrictions} page~\pageref{sec:dimRestrictions} 
for details.
\end{remark}

\begin{remark}

Therefore we will only be considering pairs of reductive dual pairs of the form
$$
 \centerline{
\xymatrix{
& 
G_{ \fA } 
 \ar@{-}[dl]\ar@{-}[dr]
& 
\\
H_{ \fA } 
& 
& 
H\p_{ \fA } 
}
} 
$$ 
where
\begin{enumerate}
\item
\dcase nn
\item
\dcase { n + 2 }n
\item
\dcase { n + 2 }{ n + 2 }
\end{enumerate}

\end{remark}

\begin{defi}[Globally admissible quadruple]
Let $ \left ( q, V \right ) $ and $ \left ( q\p, V\p \right ) $ be quadratic spaces defining orthogonal groups $ H $ and $ H\p $, respectively, and let $ \xi$ and $\xip $ be automorphic quadratic characters of $ H\ad  $ and $ H\p\ad $, respectively.
\par
A quadruple $ \adata $ is called a \bbb{globally admissible quadruple} if and only if the quadruple
$ \adatanu  $
is locally admissible 
at every place $ \nu $
where $ \xileft= \otimes_{ \nu } ^{ \prime } \xi_{ \nu }  $ and 
$ \xiright= \otimes_{ \nu } ^{ \prime } \xi_{ \nu }\p  $.

\end{defi}

\newpage

\section{ From local morphisms to local pairings }\label{sec:morphismsToPairingsTop}
In this section we define what properties our local pairings should satisfy. We then show that the existence of local morphisms such as those appearing in Snitz's thesis can be equivalently recast as certain equivariant pairings. Finally we define a natural global pairing which will be the focus of this work. This global pairing is essentially an inner product of two global theta integrals.
A keep point to keep in mind is that the pairing we study is between Schwartz spaces instead of between co-invariants.
The advantage of working with Schwartz spaces over co-invariants is that the Schwartz spaces do not depend on the parameters defining the quadruple $ \adata $. 
\par
To sum up constructing morphisms between coinvarints is a very difficult problem. We replace this problem with finding pairings between coinvariant spaces. Finally we replace this problem with finding pairings between Schwartz spaces where the pairing satisfyies certain equivariance properties with respect to the actions of the orthogonal groups 
$ H $ and $ H\p $
and the symplectic (or metaplectic) group $ G $.

\subsection{ Equivariance properties of our pairings }
In this section we describe the equivariance properties of the pairings of interest.
\par
Let  $ \F $ be a local field.
The quadratic spaces, groups etc. will all be given locally therefore, for instance, $ V = V_{ \F }  = V_{ \nu }  $ denote one and the same thing in this section.
Given two morphisms between two local theta lifts we construct a good pairing between these two local theta lifts. That is we define
\begin{defi}[Good local pairing]
\label{sec:goodPairing}
A pairing
 $$
B \colon S  \left( \voxpr_{ \F }  \right) \times S  \left( \voxprp_{ \F }  \right) \rightarrow \fC
$$
is called \bbb{good} if it satisfies the following properties for every
 $ \phi\in  S \left ( \voxpr_{ \F }  \right ) ,\phip\in  S \left ( \voxprp_{ \F }  \right ) $
\begin{enumerate}
\item
\textit{(
 $  \left (  H_{ \F } ,\xi \right )  $-equivariance
)}
For every  $ h\in H_{ \F }  $
$$
B \left ( 
\omega_{ \psi, V }  \left (  h \right ) \phi
 ,\phip \right ) 
= 
 \xi \left ( h \right ) 
 B \left (\phi,\phip \right ) 
$$
\\
\item
\textit{(
 $  \left (  H_{ \F }\p,\xi\p \right )  $-equivariance
)}
For every  $ h\in H\p_{ \F }$
$$
B \left ( 
\phi ,
\omega_{ \psi^{ -1 }, V\p }  \left ( h\p \right )
\phip \right ) 
= 
 \xi\p \left ( h\p \right ) 
 B \left (\phi,\phip \right )  
$$\\
\item
\textit{(
 $ G $-invariance
)}
For every  $ g\in G $ 
$$
B \left ( 
\omega_{ \psi, V }  \left (  g \right ) \phi 
,
\omega_{ \psi^{ -1 }, V\p }  \left (  g \right ) \phip
\right )
= 
 B \left (\phi,\phip \right )
$$
\end{enumerate}
\end{defi}
\begin{remark}
Similarly we can define a good global pairing.
$$
B \colon S  \left( \voxpr \ad  \right) \times S  \left( \voxprp \ad  \right) \rightarrow \fC
$$
satisfying the above equivariance properties.
\end{remark}
\subsection{ Transition from morphisms to pairings }
\label{sec:morphismsToPairings}
In this section we describe how to construct good pairings from morphisms.
Let's suppose, as in Snitz's thesis, that we have constructed a local morphism between local theta lifts
$$
\Delta\,\colon \co{ \psi_{ \nu }  }{ 3 }  \to \cop{ \psi_{ \nu }  }{ 1 } 
$$
Note that there is a natural local pairing
$$
B \,\colon  \cop{ \psi_{ \nu }  }{ 1 } \times \widetilde{  \cop{ \psi_{ \nu }  }{ 1 } } \to \fC
$$
and from properties of the Weil representation, i.e.  $ \widetilde{ \omega_{ V, \psi_{ \nu }  } } \cong
\omega_{ V, \psi_{ \nu } ^{ -1 }  }
 $
(see Lemma \ref{sec:contraOfWeilRep} part (\ref{notation:contraWeilPsia}) page \pageref{sec:contraOfWeilRep})
  we have a canonical isomorphism
 $$
 \widetilde{  \cop{ \psi_{ \nu }  }{ 1}  }
\to
\cop{ \psi_{ \nu } ^{ -1 }  }{ 1 } 
$$
Therefore we actually have a pairing
$$
B_{ 0 }  \,\colon  \cop{ \psi_{ \nu }  }{ 1 } \times \cop{ \psi_{ \nu } ^{ -1 }  }{1 }  \to \fC
$$
Moreover, by composing the left coordinate of the pairing  $ B $ with  $ \Delta $ we obtain a pairing 
$$
B \,\colon  \co{ \psi_{ \nu }  }{ 3 } \times \cop{ \psi_{ \nu } ^{ -1 }  }{ 1}  \to \fC
$$
i.e. 
$$
B \left ( \phi, \phi\p \right ) = 
B_{ 0 }  \left (   \Delta \left ( \phi \right ) ,\phi\p \right ) 
$$
Now this is the same as constructing a pairing on the level of Schwartz spaces which we continue to denote by  $ B $
$$
B \colon S  \left( \voxpr_{ \F }  \right) \times S  \left( \voxprp_{ \F }  \right) \rightarrow \fC$$
which is equivariant in the left coordinate with the respect to the action of  $ H $ and on the right with the respect to the action of  $  H\p $ and invariant with respect to the diagonal action of  $ G $.

\par
Thus we went from a local morphism
$$
\Delta\,\colon \co{ \psi_{ \nu }  }{ 3 }  \to \cop{ \psi_{ \nu }  }{ 1 } 
$$
to a good local pairing
 $$
B \colon S  \left( \voxpr_{ \F }  \right) \times S  \left( \voxprp_{ \F }  \right) \rightarrow \fC
$$
\begin{remark}
Note that a key point in the transition from morphisms to pairings is the fact that the dual of the Weil representation is the dual of another Weil representation, namely $ \widetilde{ \omega_{ V_{ \nu } , \psi_{ \nu }  } } \cong
\omega_{ V_{ \nu } , \psi_{ \nu } ^{ -1 }  }
 $ (See Lemma \ref{sec:contraOfWeilRep} part (\ref{notation:contraWeilPsia}) page \pageref{sec:contraOfWeilRep}). This gives a canonical isomorphism
 $$
 \widetilde{  \cop{ \psi_{ \nu }  }{ 1}  }
\to
\cop{ \psi_{ \nu } ^{ -1 }  }{ 1 } 
$$
which is used in the transition from morphisms to pairings.
\end{remark}
To sum up, the focus of this work is to construct explicit non-trivial good local pairings and also show that certain natural global pairings are factorizable into a product of good local pairings.
\par
\subsection{ Natural global pairings }
In this section we define the actual global pairing of interest. If this pairing is non-zero then that will mean that the two representations $ \pi $ and $ \pi\p $ that we have constructed are isomorphic.
This follows from the fact that we will prove that the representations $ \pi $ and $ \pi\p $ are irreducible and then we applye Schur's lemma. 
We already know that a necessary condition for these two representations to be isomorphic is that the defining data of $ \pi $ and $ \pi\p $ is given by a globally admissible quadruple
 $  \adata  $.
\begin{defi}[Natural global pairing]
We define the \bbb{natural global pairing} 
$$
B \,\colon \schad \times \schadp  \to \fC
$$
given by
$$
B \left ( \phi, \phi\p \right ) = 
\int_{ \spquot } 
I_{ V } ^{ \phi }  \left (  \xi \right ) \left (  g \right ) 
 \, 
 \overline{ I_{ V\p } ^{ \phi }  \left (  \xi \right ) \left (  g \right )  }
\,dg
$$
\end{defi}

This pairing will be a central object in this work. One of our main goals will be to factor this into a product of explicit local pairings.
 \begin{lemma}[$ B $ is a good global pairing]
 \label{sec:BisAGoodGlobalPairing}
The global pairing
$$
B \,\colon \schad \times \schadp  \to \fC
$$
is a good global pairing, i.e. it satisfied the equivariance properties stated in Definition \ref{sec:goodPairing}.
\end{lemma}
 \proof
Trivial.
$ \square $

\subsection{ Definition of a local pairing }

In this section we define a family of local pairings $ B_{ \nu, s } $ which appear naturally from unfolding the auxiliary global pairing $ B_{ s } $.  We then present the formula for the local pairing of interest $ B_{ \nu }  $. Once again this pairing is naturally defined by the global results. Namely, globally, we take the residue of the auxiliary global pairings $ B_{ s }  $ and this naturally leads us to the definition of the good local pairing of interest.
\par

In this section  $ \F $ will denote a local field.
Recall that if $ g \in G $ where $ G $ is the symplectic or metaplectic group then by the Iwasawa decomposition 
$$
G = PK = MNK\cong GL_{ n } NK
$$
where $ P  $ is the Siegel parabolic of $ G $ and $ M $ is the Levi part of $ P $ and $ N $ is the unipotent radical of $ P $. The group $ M $ is of the form 
$$
M = \Set*{  m \left ( a \right )  = \begin{pmatrix} a & 0 \\ 0 & a^\star \end{pmatrix} } { a \in GL_{ n } }
$$
so for $ m =  m \left ( a \right )  \in M $ we define $  a \left ( g \right )  = a $, in other words $a \,\colon G \to GL_n$. Finally we denote
$$
\abs{  a \left ( g \right )  } 
 = 
 \abs{  \det \left ( a \right )  } 
$$
and 
$$
\begin{cases}
   \Psi_{ s }  \,\colon GL_n(\fA) \to \F \\
   \fPsi{ g }{ s }=\absA{  a \left ( g \right ) }^{ s + \rhon } 
   \\
\end{cases}
$$
where $ \rhon = \dfrac{ n+1 } { 2 } $.
\par
 In this section we define an explicit family of local pairings depending on a complex parameter $ s $.
\begin{defi}[The family of local pairings $ B_{ \nu, s }  $]
\label{sec:defLocalPairing}
\label{notation:localNaturalPairing}
We define a family of local pairings
$$
B_{ \nu, s } \,\colon \schloc  \times \schlocp   \to \fC 
$$
given by
\begin{align*}
B_{ \nu, s } \left ( \phi, \phip \right )  & =  
 \int_{ K_{ \nu }  } \int_{ M_{ \nu }  }  W_{ \phi }  \left (  mk \right )  W\p_{ \phip }  \left (  mk \right ) \delb m \fPsi{ mk }{ s } \,dm \,dk
\end{align*}
\end{defi}
\begin{remark}
This pairing naturally appears from unfolding the auxiliary pairing $ B_{ s } $.
\end{remark}

\begin{lemma}\label{meromorphic}
The local factors \( B_{ \nu, s }  \) is defined for \(  Re \left ( s \right ) >>0 \), i.e. the integral defining it converges in some right half plane. Moreover  \( B_{ \nu, s }  \) has a meromorphic continuation to \( \fC \). 
\end{lemma}
\proof
This follows from the asymptotics of the Whittaker-Fourier coefficient, i.e. the integrand of the local pairing is the Mellin transform of a finite linear combination of functions of finite type. 
This follows from \cite[Prop 2.2]{jacquet1979automorphic} or more generally from results of \cite{sakellaridis2017periods}.
From this it follows that the  \( B_{ \nu, s }  \) converges for \(  Re \left ( s \right ) >>0 \) and \( B_{ \nu, s }  \) meromorphically extends to all of \( \fC \).

$ \square $

\newpage
\newcommand{ \qv }[1]{ \left (  q_{ #1 } ,V_{ #1 } \right ) }
\newcommand{ \qvp }[1]{ \left (  q_{ #1 } \p ,V_{ #1 } \p \right ) }
\newcommand{ \oqv }[1]{ O \qv{ #1 } }
\newcommand{ \oqvp }[1]{ O \qvp{ #1 } }
\newcommand{ \slpair }[2]{ The pair $ \oqv{ #1 }\adpr \times \metasl $ and
the pair 
$ \oqvp{ #2 }\adpr \times \metasl $.
 }
\section{ Proof of the main statement }
In this section we prove the main statement for the low rank case of $ n = 1 $ although most of our proofs do generalize to the higher rank case. Namely we restrict ourselves to the following three low rank cases, i.e. pairs of dual pairs of the form
\begin{enumerate}
\item
The pair $ \oqv{ 1 } \times \metasl$ and
the pair 
$ \oqvp{ 1 }\times \metasl $.
\item
The pair $ \oqv{ 1 } \times \metasl$ and
the pair 
$ \oqvp{ 3 }\times \metasl $.
\item 
The pair $ \oqv{ 3 } \times \metasl$ and
the pair 
$ \oqvp{ 3 }\times \metasl $.
\end{enumerate}
where $ \qv n $ denotes an $ n $-dimensional quadratic space.
\par
We prove the main statement in the following steps.
\begin{enumerate}
\item
To begin with we unfold the auxiliary global pairing $ B_{ s } $. The key tool here is the Fourier expansion of the theta integrals.
It is important to note that the Fourier expansion of the theta integrals is relatively simple and sparse. Multiplying the theta integrals by the normalized auxiliary Eisenstein series $ \cE^{ \star }  $ greatly facilitates the unfolding, turning an impossible unfolding problem of the pairing of interest $ B $ which does not involve an Eisenstein series into a tractable problem of unfolding the auxiliary pairing which integrates over a product of theta integrals against a certain simple Eisenstein series. 
\par
Essentially we unfold a Rankin-Selberg integral. The unfolding is very similar to the derivation in Jacquet’s \cite{jacquet2006automorphic} application of the Rankin-Selberg method to the  group $  GL \left ( 2 \right ) \times  GL \left ( 2 \right )  $. 
\item
Next we take the residue of $ B_{ s } $ at $ s = \rhon = \frac{ n+1 } { 2 } $ to obtain a factorization of our original natural global pairing $ B $. The key points here are the fact that our theta lifts are cuspidal and that the residue of the normalized auxiliary Eisenstein series $ \cE^{ \star } $ is a constant function identically equal to one.
\item
Finally one can prove the main statement. Namely for an admissible quadruple
if a certain global partial $ L- $function
has a simple pole at $ \rhon $
then the representations $ \pi $ and $ \pi\p $ are isomorphic, i.e. $ \pi \cong \pi\p $.
Note that the global partial $ L- $function is obtained from the unramified computation.
\begin{remark}
Actually in the one dimensional case one can even prove $ \pi = \pi\p $ since multiplicity one is satisfied in this case. However in the higher rank case it is not clear that multiplicity one is satisfied for representations of the metaplectic group considered in this work.
\end{remark}
 
\end{enumerate}

\subsection{ Unfolding the auxiliary global pairing }

\label{sec:unfoldGlobal}
In this section we unfold the auxiliary global pairing $ B_{ s } $ which is an integral of the two theta integrals against a relatively simple normalized auxiliary Eisenstein series $ \cE^{ \star } $. 
The key property of $ \cE^{ \star } $ is that it is a normalization of a certain Eisenstein series $ \cE $, which we refer to as the auxiliary Eisenstein series, whose residue at $ \rhon $ is a constant function. Thus after normalizing $ \cE $ we obtain an Eisenstein series $ \cE^{ \star }  $ whose residue at $ \rhon $ is a constant function that is identically equal to one.
\par
The proof is an adaptation of Jacquet's \cite{jacquet2006automorphic} $ GL_{ 2 } \times GL_{ 2 }  $ Rankin-Selberg convolution unfolding of a similar pairing. The proofs are very similar although the context is different.
\par
Although slightly tedious most of the steps in the unfolding are straightforward (especially after reading Jacquet's work!) A crucial tool is the global Fourier expansion of each of the theta integrals in the pairing as described in 
Lemma \ref{sec:FourierExpLemma} page~\pageref{sec:FourierExpLemma}.

Our goal is to factor the global pairing
$$
B \,\colon \schad \times \schadp \to \fC 
$$

In order to do this we will consider a family of pairings.
\begin{remark}
We have restricted the result of the main statement Theorem \ref{sec:mainStatement} to the $ n = 1 $ case since we only did the unramified computation for this case. Besides the restriction of the unramified computation the proof is completely general. Also throughout this work we are considering globally anisotropic quadratic spaces and non-trivial automorphic characters of the orthogonal groups.
\par
Clearly the main results in this work can not be true for arbitrary globally anisotropic quadratic spaces and arbitrary non-trivial automorphic characters. This is precisely where the notion of admissible quadruples comes into play. In other words a central restriction in this work is to consider admissible quadruples $ \adata $. These restrictions come from relatively trivial local necessary conditions. Namely an admissible quadruple $ \adata $ satisfies two local necessary conditions which we have referred to as the central character condition and the Fourier coefficient condition.
See section \ref{sec:admissibleQuadrupleExamples} page~\pageref{sec:admissibleQuadrupleExamples} for more details on admissible quadruples.
\end{remark}

\begin{defi}[Auxiliary Eisenstein series]
\label{sec:auxEis}
We define the \bbb{auxiliary Eisenstein series} by\label{sec:defSimpleEisSeries}
$$
\aeis g s = 
\sum_{ \gamma \in \rquot P G } 
\fPsi{ \gamma g }{ s }
$$
where 
$
\fPsi{ g }{ s }
= 
\absA{  a \left ( g \right ) }^{ s + \rhon }
$ and $ \rho_{ n } = \frac{ n+1 } { 2 } $.
\end{defi}

\begin{defi}[Normalized Auxiliary Eisenstein series]

Define the \bbb{normalized auxiliary Eisenstein series} by\label{sec:defAuxEisSeries}
 $$
 \naeis
g s
= 
\glblZetaConst
\cdot
\aeis g s
$$
which is normalized so that $ \cE^{ \star }  $ has a residue of  $ 1 $ at the simple pole 
$ s  = \rhon $.
\end{defi}
\begin{remark}
The factor $ \glblZetaConst $ appearing in Definition \ref{sec:defAuxEisSeries} will not be given and is not important.
The only important fact is that such a normalization exists and with this normalization the residue of the normalized auxiliary Eisenstein series is a constant function identically equal to one.
For the curious reader $ \glblZetaConst $ is a certain product of L-functions and explicit formulae for $ \glblZetaConst $ can be found in \cite{kudla1990poles}.
\end{remark}
Now we can define a family of auxiliary pairings which will be relatively easy to unfold.
\begin{defi}[Auxiliary pairing]
\label{auxPair}
We define the auxiliary pairing to be the pairing
$$
\apair \,\colon \schad \times \schadp \to \fC 
$$
given by 
$$
B_{ s }  \left ( \phi, \phip \right ) 
=
\int_{ \aq G } 
\thLift \phi g
\overline{ \thLiftp \phip g }
\naeis g s
\,dg
$$
\end{defi}

Before we 
present the unfolding of the global pairing which is the content of Theorem 
\ref{sec:unfoldGlobal}
we recall some results from Fourier analysis.

\subsubsection{ Fourier expansion of the global theta integrals }
Throughout this section we assume our defining data is an admissible quadruple $ \adata $.
The following functions may be obtained from some standard Fourier analysis of the global theta integrals. See Section \ref{sec:GlobalFourierAnalysis} page~\pageref{sec:GlobalFourierAnalysis}
  for more details.
\par

Recall the Fourier functions from Definition \ref{sec:defFourierFunction}.
  $$
W_{ \phi }  \left ( g \right ) = 
\Hintegral
\left (
\wfunc  \left ( g, h \right ) 
\phi \right )
\left (  
\vv
\right ) 
\xi \left ( h \right ) 
\,dh
$$
where $ \phi \in
\voxpr_{ \fA } ,
g \in G\ad
 $ and we also defined
$$
W\p_{ \phip }  \left ( g \right ) = 
\Hintegralp
\left (
\wfuncpi \left ( g, h \p  \right ) 
\phip \right )
\left (  
\vv\p
\right ) 
\xi\p \left ( h \p \right ) 
\,dh\p
$$
where
$ \phip \in
\voxprp_{ \fA } ,
g \in G\ad
 $.

Recall the Fourier expansions proven in Lemma \ref{sec:globalFourierExpansion}. 
\begin{enumerate}
	\item We have the Fourier expansion
\begin{align*}
   \thLift \phi g     =  \sum_{ \alpha \in M_\F } 
	  W_{ \phi } \left ( \alpha g \right )
\end{align*}
\item Similarly, we have the Fourier expansion
\begin{align*}
   \thLiftp \phip g  & = \sum_{ \alpha \in M_\F }
	 W_{ \phip } \left ( \alpha g \right )
\end{align*}
\end{enumerate}
where $ M $ is the Levi part in the Siegel parabolic subgroup $ P $ of $ G $.

\begin{theorem}[Unfolding of the global auxiliary pairing]
\label{sec:unfoldGlobalThm}
For  $  Re \left ( s \right ) > \auxRange $ we have
$$
B_{ s }
\left ( \phi, \phip \right ) 
= 
\glblZetaConst
\cdot
\int_{ \ga K } \int_{ \ga M }  W_{ \phi }  \left (  mk \right )  W\p_{ \phip }  \left (  mk \right ) \delb m \fPsi{ mk }{ s } \,dm \,dk
$$
where our data $ \adata $ is given by a globally admissible quadruple.
\end{theorem}

\proof
See Appendix \ref{SomeProofs}.

$ \square $

\begin{remark}
The unfolding of the global pairing $ B_{ s } $ is almost identical to Jacquet’s application of the Rankin-Selberg method to the  group $  GL \left ( 2 \right ) \times  GL \left ( 2 \right )  $. A crucial tool in Jacquet's proof is the Fourier expansion of cusp forms \cite{jacquet2006automorphic}.
\par
Similarly, the Fourier expansion of the theta integrals $ \thlift {   } $ and $ \thliftp {   } $ with respect to the unipotent radical $ N $ of the Siegel parabolic subgroup $ P \subseteq G $  will be of central importance in the unfolding of a the auxiliary global pairing 
$ B_{ s }  
 \,\colon \schad  \times \schadp  \to \fC
$.
\end{remark}
\subsection{ Factorization of the global pairing }
In this section we obtain a factorization of our original natural global pairing $ B $ by taking a residue of $ B_{ s }  $ and using the cuspidality of the representations $ \pi $ and $ \pi\p $. 
\begin{remark}
For the argument in Theorem \ref{sec:pairing} to be valid it is enough that only one of the representations $ \pi $ or $ \pi\p $ are cuspidal.
\end{remark}


\begin{theorem}[Natural global pairing as a residue of the auxiliary pairing]
\label{sec:pairing}
Given a globally admissible quadruple $ \adata $.
For
large $  Re \left ( s \right ) $ the 
pairing
$ B_{ s }  $ is convergent and has meromorphic continuation. Moreover $ B_{ s }  $
has a pole of order no more than one.
Moreover
\begin{align*}
    B \left ( \phi, \phip \right )  & =  Res_{ s = \rhon }  B_{ s } \left ( \phi, \phip \right )   \\
\end{align*}
\end{theorem}
\proof
By Theorem \ref{sec:unfoldGlobal} we have 
$$
B_{ s }
\left ( \phi, \phip \right ) 
= 
\glblZetaConst
\cdot
\int_{ \ga K } \int_{ \ga M }  W_{ \phi }  \left (  mk \right )  W\p_{ \phip }  \left (  mk \right ) \delb m \fPsi{ mk }{ s } \,dm \,dk
$$
and using the fact that the global theta lifts are cuspidal for an admissible quadruple we have 
$$
B \left (  \phi,\phip \right ) 
 = 
 Res_{ s = \rhon }
 B_{ s }  \left (  \phi,\phip \right ) 
$$
together with the fact that we normalized $ \cE^{ \star } $ so that its residue at $ \rhon $ is a constant function identically equal to one.
\par
The convergence properties follow from properties of the Eisenstein series $ \cE^{ \star }  $.
\par
More precisely
\begin{align*}
    B \left ( \phi, \phip \right )  & =  Res_{ s = \rhon }  B_{ s } \left ( \phi, \phip \right )   \\
	 &   = 
Res_{ s = \rhon }
 \left (  
 \glblZetaConst
 \cdot	 
 \int_{ \ga K } \int_{ \ga M }  W_{ \phi }  \left (  mk \right )  W\p_{ \phip }  \left (  mk \right ) \delb m \fPsi{ mk }{ s } \,dm \,dk
 \right )
   \\
	 &   = 
Res_{ s = \rhon } \glblZetaConst
 \cdot	 
 \int_{ \ga K } \int_{ \ga M }  W_{ \phi }  \left (  mk \right )  W\p_{ \phip }  \left (  mk \right ) \delb m \fPsi{ mk }{ s = \rhon } 
 \\
 	 &   = 
\kappa
 \cdot	 
 \int_{ \ga K } \int_{ \ga M }  W_{ \phi }  \left (  mk \right )  W\p_{ \phip }  \left (  mk \right ) \delb m \fPsi{ mk }{ s = \rhon } 
\end{align*}
where
$ \kappa
 = 
 Res_{ s = \rhon } \glblZetaConst
$.
$ \square $

\begin{defi}[Normalized local pairing]
For $ n = 1 $ and for every unramified place $ \nu $ the normalized local pairing 
$$
B_{ \nu } ^{ \star }  \,\colon \schloc \times \schlocp \to \fC
$$
is defined to be
$$
 B_{ \nu } ^{ \star }
 \left ( \phi_{ \nu } , \phip_{ \nu } \right ) 
 = 
\dfrac{  B_{ \nu }
 \left ( \phi_{ \nu } , \phip_{ \nu } \right ) 
 } { L_\nu \left ( 1 \right )  }
$$
\end{defi}

\begin{theorem}[Factorization of pairing]
\label{sec:mainthm}
For an admissible quadruple $ \adata , \phi \in \schad,\phip \in \schadp  $ and for $ n = 1 $
 \begin{align*}
    B \left ( \phi, \phip \right )  & =
 \kappa
\cdot
\left \{ 
\prod_{ \nu \in S } 
 B_{ \nu, s = \rhon } 
\left ( \phi_{ \nu } , \phip _{ \nu } \right )  
\right \}\cdot
Res_{ s = \rhon } 
L_{ S }  \left ( s \right ) 
 \end{align*}

where 
$ \kappa
= 
Res_{ s = \rhon } \aeis g s
 $
is the residue of the unnormalized Eisenstein series used in the unfolding of the global auxiliary pairing.
\end{theorem}

\proof
This follows from 
Theorem \ref{sec:pairing} together with the fact that $ \pi $ is cuspidal and we have normalized the Eisenstein series so that 
$$
 Res_{ s = \rhon }\naeis g s = 1
$$
Namely Theorem \ref{sec:pairing} states that 
\begin{align*}
   Res_{ s=\rhon }  \left <  I_{ V,\xi }^{ \phi }   , I_{ V\p,\xi\p }^{ \phip } \cdot \cE^{ \star }   \right >_{ G\ad }    & =  \left <  I_{ V,\xi }^{ \phi }   , I_{ V\p,\xi\p }^{ \phip } \cdot Res_{ s=\rhon } \cE^{ \star }   \right >_{ G\ad }     \\
	 &   =  \left < I_{ V,\xi }^{ \phi }  , I_{ V\p,\xi\p }^{ \phip }    \right >_{ G\ad }   \\
	 &   =  B\left ( \phi, \phip \right )   
\end{align*}

Moreover by the unramified computation Lemma \ref{sec:finalUnramifiedComputation} at almost every place $ \nu $, in the unramified setting, outside of a finite set $ S $
we have
$$
B_{ s, \nu } \left ( \phi_{ \nu }, \phip_{ \nu } \right ) 
 = 
 L_\nu \left ( s \right ) 
 =  \left (  1-q_{ \nu } ^{ -s }  \right ) ^{ -1 } 
 ,
 \quad
 \forall v\not\in S
$$
Finally we have 
\begin{align*}
    B \left ( \phi , \phip \right )   & =    
\kappa
\cdot
\left \{ 
\prod_{ \nu \in S } 
 B_{ \nu, s = \rhon } 
\left ( \phi_{ \nu } , \phip _{ \nu } \right )  
\right \}\cdot
\prod_{ \nu \not \in S } 
 L_\nu \left ( s \right ) 
  \\
&   = 
\kappa
\cdot
\left \{ 
\prod_{ \nu \in S } 
 B_{ \nu, s = \rhon } 
\left ( \phi_{ \nu } , \phip _{ \nu } \right )  
\right \}\cdot
Res_{ s = \rhon } 
L_{ S }  \left ( s \right ) 
  \\
 \end{align*}
where in the last transition we completed the partial $ L- $function to the Dedekind zeta function and we used that fact that the Dedekind zeta function has a simple pole whose value is $ 1 $ at $ \rho_{ 1 }  = 1 $.

$ \square $

 \begin{remark}
Note that by Lemma \ref{meromorphic} the local pairings are defined in a right hand plane and have meromorphic continuation to all of \( \fC \).
\end{remark}

\subsection{ Proof of the main statement }
In this section we prove the main statement Theorem \ref{sec:mainStatement}. Namely for an admissible quadruple the representations $ \pi $ and $ \pi\p $ are isomorphic if and only if a certain global condition is satisfied where this global condition is given in terms of the existence of a pole at $ \rhon $ of a certain global partial $ L- $function.

\begin{customthm}{\ref{sec:mainStatement}}[Main statement]
Given a globally admissible quadruple $ \adata $
the global theta lifts
$$ \pi =  \gLift \tn{ and } \pi\p =  \gLiftp $$
are isomorphic, i.e. $ \pi\cong \pi\p $.
\end{customthm}
\proof
It is known that the Dedekind zeta function has a simple pole at $ s = 1 = \rho_{ 1 }  = \dfrac{ 1+1 } { 2 } $ whose value is equal to one.
This is a crucial fact and is implicit in the computations of Theorem \ref{sec:mainthm}.
 Let $ \phi_{ \nu }  \in \schloc,\phip_{ \nu } \in \schlocp $ be local data so that the integral defining $ B_{ \nu }
 \left ( \phi_{ \nu } , \phip_{ \nu } \right )  $ is nonzero at $ s = \rhon = 1 $.
 Now we have obtained a non-degenerate bilinear global pairing $ B $. 
For an admissible quadruple $ \adata $ we can choose certain Schwartz functions $ \phi \in \schad,\phip \in \schadp $ such that
the pairing $ B_{ \nu, s } $ is nonzero at $ s= \rhon $.
This essentially follows from the proof of first occurrence, namely from a similar argument to Theorem \ref{sec:nonzeroLift}.
 Let $ \phi, \phip $ be the choice of data such that $ B_{ \nu, s } $ is non-zero.
\par
Consequently, the fact that $ B_{ s }  $ has a pole implies that 
$$
B \left (  \phi,\phip \right ) \neq 0
$$
Now we have obtained a non-degenerate bilinear pairing $ B $. 
This proves that the representations $ \pi $ and $ \pi\p $  are isomorphic.
$ \square $
\begin{remark}
As stated in the introduction for an admissible quadruple $ \adata $ if a certain global $ L $-function has a pole at $ \rhon $ then $ \pi $ and $ \pi\p $ will be isomorphic. This global $ L $-function is precisely the Dedekind zeta function in our low rank setting. Note that this would have been more apparent if we would have taken a slightly more general Eisenstein series.
See for example \cite[Proposition $ 3.8.5 $]{bump1998automorphic}. Essentially the central character condition ensured that we would obtain the Dedekind zeta function.
\end{remark}

\newpage
\section{ The Unramified Computation }
\label{sec:unrSetting}
In this section we restrict the results to the case of $ n = 1 $. Namely from the dimension restrictions of Lemma \ref{sec:dimRestrictions}, i.e. that for an admissible quadruple $ \adata $ the dimensions of the quadratic spaces are must satisfy
 $ 
 dim \left ( V \right ) , dim \left ( V\p \right ) \in  \left \{ n, n + 2 \right \} $, i.e. $ V $ and $ V\p $ maybe either one or three dimensional where we are lifting to $ \widetilde{ Sp }_{ 2n }  = \widetilde{ SL_2 } $.

The goal of this section is to compute the local pairing $ B_{ \nu, s }  $ in the generic setting. 
This is referred to as the unramified computation.
In other words we would like to calculate the following integral in the unramified setting which is described in Section \ref{sec:unrs}.
$$
 B_{ \nu, s } 
\left ( \phio, \phiop \right ) 
= 
\int_{ \fa M }  W_{ \phio }   \left (  m \right )  W^{ \prime } _{ \phiop }  \left (  m \right ) \delb m \fPsi{ m }{ s } 
\,dm
$$
\begin{remark}
It is important to note that our local pairings depend on choices of vectors $ \vv $ and $ \vv\p $.
\end{remark}
\subsection{ Unramified Setting
\conf{ 95 }
 }
 \label{sec:unrs}
  In this section we describe the unramified setting.
 
 Let 
 $$
\left (  L \subseteq V, \left (  q,V \right ), \xi \,\colon H_{ \nu } \to \fC^\times  \right ) 
 \tn{ and } \left (  L\p \subseteq V \p, \left (  q\p,V\p \right ), \xi\p \,\colon V\p \to \fC^\times  \right )  
$$
 be the unramified data where $ L $ and $ L\p $ are lattices in $ V $ and $ V\p $, respectively.
 We assume that our characters, quadratic spaces and lattices in the corresponding vector spaces are unramified.
More precisely 
 \begin{enumerate}
 \item
Let $ \F $ be a non-Archimedean field. We assume 
that $ W_{ \phio }   , W_{ \phiop } $ are right $ K -$ invariant with respect to the maximal compact subgroup $ K $ if $ G $. In particular 
 $$
W_{ \phio }   \left (  e \right ) =  W_{ \phiop }  \left (  e \right ) = 1 
$$
\item
The vectors
 $ \vv\in  \left (  V \otimes X \right ) \left (  \cO_{ \nu } \right ) \cong V^{ n }  \left (  \cO_{ \nu } \right ),\vv\p\in  \left (  V\p\otimes X  \right )  \left (  \cO_{ \nu } \right ) \cong V^{\prime,  n }  \left (  \cO_{ \nu } \right )$ which are the choices of vectors such that 
 $$
\gram{ \vv } = \gramp{ \vv\p }
$$
which appear in the definitions of $ W_{ \phio }   , W_{ \phiop }  $. 
\begin{remark}
Note that $ X $ is a maximally isotropic subspace of $ W $ where $  \left (  W,\left < \cdot,\cdot  \right > \right )  $ is the symplectic space defining $ G $.
This is part of the typical construction of the Schrödinger model of the Weil representation (see Section \ref{sec:ScrodModel} for more details).
\end{remark}
 \item
 For every $ \ell \in L $ we have $  q \left ( \ell \right ) \in \cO_{ \nu }  $. Similarly, 
for every $ \ell\p \in L\p $ we have $  q\p \left ( \ell\p \right ) \in \cO_{ \nu }  $. Note that in our setting $ L =  V \left ( \cO_\nu \right )  $ and $ L\p =  V\p \left ( \cO_\nu\p \right ) $.
\item
We have  $ \phio\in  S \left ( V _{ \nu }  \right ) ,\phiop\in  S \left ( V\p _{ \nu }  \right )  $ are two Schwartz functions defined on the lattices
 $$
\phio
= 
1_{ V \left (  \cO_{ \nu } \right )  } 
\tn{ and }
\phiop
= 
1_{ V\p \left (  \cO_{ \nu } \right )  } 
$$
\item
 The discriminants of our quadratic spaces satisfy
 $ d,d\p\in \cO_{ \nu } ^{ \times }  $.
 \item
The characters $ \xi, \xip $ are unramified, that is, trivial on $  H \left ( \cO_\nu \right )  $ and $  H\p \left ( \cO_\nu \right )  $, respectively.
 \end{enumerate}

\subsection{ The computation }
In this section we actually compute the local pairing in the unramified setting for $ n = 1 $.
We compute each of the factors separately in the above pairing.
\begin{lemma}
\label{simpleUnr}
In the unramified setting consider the reductive dual pair $ H\p\times G $ where $ H\p $ is given by a one-dimensional quadratic space $ \left ( q\p, V\p \right ) $. Then 
$$
 W^{ \prime } _{ \psi^{ -1 } ,\phiop }  \left (  m \left (  a \right )  \right ) = 
\gamma
\cdot
\abs{ a } ^{ 1 / 2 } 
\cdot
1_{ \cO_\nu }  \left ( a \right ) 
$$
\end{lemma}
\proof
In the one dimensional case in the unramified setting $ \xi\p\equiv1 $.
Moreover recall the normalization $ W_{ \phiop }  \left (  e \right ) = 1  $.
Hence we have
\begin{align*}
1
& =
 W^{ \prime } _{ \psi^{ -1 } ,\phiop }  \left (  m \left (  a \right )  \right ) \\
& = 
\int_{ H\p_\nu } 
\left (
\wfuncpi \left (  m \left ( a \right )\right ) 
\phiop \right )
\left (
h^{ \prime, -1 } 
\vv\p
\right ) 
\xi\p \left ( h \p \right ) 
\,dh\p
 \\
& = 
\gamma
\cdot
\abs{ a } ^{ 1 / 2 } 
\underbrace{  \chi_{ V\p } \left ( a \right )  }_{ = 1 }
\cdot
\dfrac{ 1 } { 2 }\left [ 
\phiop \left ( a\vv\p \right )  +  \phiop \left ( a\vv\p \right )
\right ]
\\
&   = 
\phiop \left ( a\vv\p \right )
\\
& = 
1_{ \cO_\nu }  \left (  a \right ) 
\end{align*}
and this equals one if and only if $ a \in \cO_{ \nu }  $ and zero otherwise, since in the unramified setting at least one of the coordinates of $ \vv  $ is in $ \cO_{ \nu } $ and $ \vv \in V \left (  \cO_{ \nu } \right )  $.
\par
Therefore 
$$
W^{ \prime } _{ \psi^{ -1 } ,\phiop }  \left (  m \left (  a \right )  \right ) = 
\gamma
\cdot
\abs{ a } ^{ 1 / 2 } 
\cdot
1_{ \cO_\nu }  \left ( a \right ) 
$$
$ \square $
\begin{lemma}
In the unramified setting consider the reductive dual pair $ H\times G $ where $ H $ is given by a three-dimensional quadratic space $ \left ( q\p, V\p \right ) $. Then 
$$
W_{ \phio,\psi }   \left (  m \left (  a \right )  \right ) 
W_{ \phio\p,\psi^{ -1 }  }   \left (  m \left (  a \right )  \right )
  =
\abs{ a } \cdot
1_{ \cO_\nu }  \left ( a \right ) 
$$
\end{lemma}
\proof
Snitz proves that 
$$
W_{ \psi, \phio }  \left (   m \left ( a \right )  \right ) 
= 
W_{ \psi, \phiop }  \left (   m \left ( a \right )  \right ) 
$$
where the left hand side it the three dimensional lift and the right hand side is the one dimensional lift.
This is contained in the proof of the fundamental lemma in \cite{snitz}, Proposition $ 35 $.
\par
Together with Lemma \ref{simpleUnr} and the cancellation of $ \gamma\cdot\gamma^{ -1 } = 1 $ we obtain the result.
$ \square $

\begin{theorem}[Unramified computation $ n = 1 $]
\label{sec:finalUnramifiedComputation}
In the unramified setting 
$$
B_{ \nu, s }  \left ( \phio, \phiop \right ) 
= 
\left (  1-q_{ \nu }^{ -s }  \right ) ^{ -1 } 
= 
\zeta_{ \nu } \left (  s \right ) 
$$
We will refer to this as the unramified computation.
\end{theorem}
\proof
\begin{align*}
 B_{ \nu, s }  \left ( \phio, \phiop \right )  
& = 
\int_{ \fa M }  W_{ \phio }   \left (  m \right )  W^{ \prime } _{ \phiop }  \left (  m \right ) \delb m \fPsi{ m }{ s } 
\,dm
\\
& = 
\int_{ \F_{ \nu } } 
\abs{ a } \cdot
1_{ \cO_\nu }  \left ( a \right ) 
\delb \mm \fPsi{ \mm }{ s } 
\,da
\\
& = 
\int_{ \cO_\nu } 
\abs{ a } \cdot
\abs{ a } ^{ -2 } \abs{ a } ^{ s + 1 } 
\,da
\\
& = 
\int_{ \cO_\nu } 
\abs{ a } \cdot
\abs{ a } ^{ -2 } \abs{ a } ^{ s + 1 } 
\,da
\\
& = 
\int_{ \cO_\nu } 
\abs{ a }
\,da
\\
& = 
\left (  1-q_{ \nu }^{ -s }  \right ) ^{ -1 } 
\end{align*}

$ \square $
\begin{remark}
Note that the unramified computation does not depend on our characters \( \xi, \xi\p \). This is due to the central character condition. This would not be true for a quadruple of data \( \adata \) which is not admissible.
\end{remark}

\newpage
\section{ 
Fourier Analysis, cuspidality and irreducibility of the representations
 }

\subsection{ 
Results from Global Fourier Analysis
 }
 \label{sec:GlobalFourierAnalysis}
Fourier analysis plays two important roles in this work. It is used in the unfolding of our global auxiliary pairing $ B_{ s }  $ and it is also used in our application of Rallis' theory of towers (see Section \ref{sec:RallisTheoryOfTowers} for details), namely we use Fourier analysis over a locally compact abelian group $ N $ in order to determine whether or not a theta lift is non-zero in a tower of theta lifts. To ensure cuspidality of a global theta lift $ \sigma $, we will be concerned with the first time the global theta lift in a tower of lifts is non-zero. Hence we will need to understand the Fourier analysis of theta integrals in a tower of theta lifts and not just the two lifts that are the main concern of this work.
\par
It is known that the Fourier coefficients which we consider are factorizable. The goal of this section is to explicitly find such a factorization by unfolding the Fourier coefficient with respect to a character $ \psi_{ \beta } \,\colon N\ad \to \fC  $ of the unipotent radical $ N $ of the Siegel parabolic subgroup $ P $ in $ G $.
\par
Finally note that we can do Fourier analysis with respect to the unipotent radical of the Siegel parabolic $ N \subseteq P \subseteq G $ because $ N $ is a locally compact abelian group.
\par

In this section we consider the reductive dual pair $H  _\fA \times G_\fA $ where $  \left (  q, V \right )  $ is an  $ m $-dimensional globally anisotropic quadratic space and $ G $ is either the symplectic or metaplectic group given by a symplectic vector space $  \left (  W,\left < \cdot,\cdot  \right > \right )  $ of dimension $ n $ and $ X $ denotes a maximal isotropic subspace of $ W $.
\subsection{ Some definitions needed in the statement of this section }
\begin{remark}
Recall that 
$$
\thLift \phi g =
\int_{ \aq H }
\sum_{ \gamma \in  \left (  \vox \right ) _{ \F } } 
 \left ( 
\omega_{ \vpsi } 
 \left ( g, h \right ) 
\phi
  \right ) 
 \left (  \gamma \right ) 
\xi \left ( h \right )  
\,dh
$$
\end{remark}
\begin{remark}
Note that the Fourier coefficient is 
Eulerian
$$
 \fthLiftBeta \phi  \left (  g \right )
 = 
W_{ \phi, \beta }  \left (g \right ) = 
\prod_{ \nu } 
W_{ \phi, \beta, \nu }  \left (g_{ \nu }  \right )
$$
where  $ \beta\in  \SymF  $ and $ \phi \in S\voxpr_{ \fA }  ,
\phi
 = 
\otimes \p_{ \nu } \phi_{ \nu } $.
\end{remark}

\subsection{ 
Fourier functionals are Eulerian
}
The purpose of this section is to find an explicit factorization of the Fourier coefficients $ \thLiftBeta \phi { \cdot } $.

$$ 
\begin{dcases}
   S( (  \underbrace{ V \otimes X }_{ \cong V\ad  ^{ n }  } )\ad) \to \Fns( G_{ \fA }, \psi ) \\
\end{dcases}
 $$
 
which is given by applying the $\psi_{ \beta }$-Fourier coefficient to the theta lift of $\xi $ using the theta kernel.
\par
We know that this function should be factorizable.
In the following lemma we find an explicit factorization of the image of this function.
\begin{lemma}\label{sec:globalFourierExpansion}
Given a quadratic space $ \left ( q, V \right ) $. Let $ \beta\in  \SymF  $. Then for every $ \phi \in S \left (  \voxpr\ad \right )  $
\begin{enumerate}
\item
Suppose $ q $ does not represent $ \beta $ then $ I^{ \phi } _{ V,\xi,\beta } = 0 $.
\item
If $ \beta $ is degenerate then $ I^{ \phi } _{ V,\xi,\beta } = 0 $.
\item
If $ \beta $ is non-degenerate then
the  $ \psi_{ \beta }  $-th Fourier coefficient of the form  $ \fthLift \phi  $ is given by
$$
\boxed{ \thLiftBeta \phi { g }
= 
 \int_{ \nicequotm } \left (  \wni (  g, h ) \phi  \right ) (\vvbeta) \xi  (h)\,dh
 \cdot
 \int_{ \badquotm }  \xi  (h_{ 0 } )\,dh_{ 0 } }
$$
where $ \vvbeta \in  \voxpr_{ \F } \cong V_{ \F } ^{ n }  $ satisfies  $ \gram\vvbeta = \beta $.
\end{enumerate}

\end{lemma}
\seeapp
\subsection{ Fourier Expansion }
In this section we find the Fourier expansion of a theta lift corresponding to data $ \qcdata $
so that the theta lift satisfies first occurrence.
This Fourier expansion is crucial in the unfolding of the normalized auxiliary Eisenstein series.
\par
The proof in this section involves substituting the Fourier coefficients that we've previously calculated in 
Lemma \ref{sec:globalFourierExpansion}
together with choosing a $ \vv\in \voxpr_{ \F } $  which defines a single orbit for the vectors $ \vv^{ \beta } \in \voxpr_{ \F } $  satisfying $ \gram \vvbeta  = \beta $. 
\par
In this section we let  $ M \subseteq P $ denote the Levi part and $ N \subseteq P $ denote the unipotent radical of the Siegel parabolic subgroup $ P  = MN$ of $ G  $.
\begin{lemma}
\label{sec:FourierExpLemma}
For $ \qcdata $ such that
the corresponding theta lift $  \Theta \left ( \eta \right )  $ satisfies first occurrence
there exists  $ \vv \in \voxpr_{ \F } $ such that the Fourier expansion of the function 
$$
b\mapsto \thLift \phi {  n \left ( b \right ) g } 
$$
is given by
\begin{align*}
   \thLift \phi { \begin{pmatrix} 1 & b \\ 0 & 1 \end{pmatrix} g }  & =
 \sum_{ \alpha\in M_{ \F }  } W_{ \phi }  \left ( \alpha g \right )
\,\cdot
\psi \left ( b \gram { \alpha \vv } \right ) 
\end{align*}
for every  $ g \in G $ and $ b \in \SymF $.
\par
In particular we have 
$$
\boxed{
\thLift \phi { g } =
 \sum_{ \alpha\in M_{ \F }  } W_{ \phi }  \left ( \alpha g \right )
 }
$$
\end{lemma}
\seeapp

\subsection{ Cuspidality }
\label{section:globalFourier}
\label{sec:proofOfCuspidality}
In this section we prove cuspidality for theta lifts of the dual pairs and characters of interest. We use first occurrence and Fourier analysis. The main result in this section is Theorem \ref{sec:IrrCuspidalAutomorphic}.
\par
In this section we denote by $  \left (  q_{ k } ,V_{ k } \right )  $ a quadratic space of dimension $ k $.
\begin{theorem}
\label{sec:zeroLifts}
Given
$ \qcdata $
for
$  \left (  q, V \right )  $ and $ \xi \,\colon H_\fA \to \fC^\times $ as defined in this work, i.e. $ \xi $ is trivial on $ H^{ \vvbeta } $ for non-degenerate $ \beta \in \SymF $.
 Then 
$$
\pi_{ k } = 
\liftc{ k }
= 
0
$$
for every $ 1 \leq k < \N $.
\end{theorem}
\seeapp

\begin{theorem}
\label{sec:nonzeroLift}
Given
$ \qcdata $
such that $ \restr { \xi } { \vvbeta } \equiv 1$ for some non-degenerate symmetric matrix $ \beta $.
Then the corresponding global theta lift is non-zero.
\end{theorem}
\seeapp

\begin{remark}
Note that for the dual reductive pair $ H \left (  \fA \right )  \times G  \left (  \fA \right )  $, where $  \left (  q,V \right )  $ is $  \dim \left ( V \right ) = n $ dimensional and when $ W $ is $ 2n $-dimensional then Theorem \ref{sec:nonzeroLift} follows immediately from Theorem \ref{sec:zeroLifts} together with Rallis'
 tower property, i.e. Theorem \ref{sec:FO}.
\par
In any case we still need the proof of Theorem \ref{sec:nonzeroLift} for the 
case where
$  \dim \left ( V \right ) = n + 2 $ and $ W $ is $ 2n $-dimensional.
\end{remark}

\begin{theorem}
\label{sec:IrrCuspidalAutomorphic}
Given $  \left (  q, V \right )  $ and $ \xi \,\colon H_\fA \to \fC^\times $ as defined in this work the representation 
$ \pi_{ N } = 
\liftc{ N }
$
is an irreducible, cuspidal, automorphic representation of $ G_{ \fA }  $ where $ \N = n$ for the dual pair $  O \left ( q_n, V_n \right ) \times G $
and $ \N = n + 2 $ for the dual pair $  O \left ( q_{ n+2 } , V_{ n+2 }  \right ) \times G $.
\end{theorem}
\proof
The cuspidality of $ \pi $ follows from Theorem \ref{sec:zeroLifts}, Theorem \ref{sec:nonzeroLift} and Rallis tower property, i.e. Theorem \ref{sec:FO}. The irreducibility of $ \pi $ follows from a general theorem of Moeglin \cite{moeglin1997quelques} stating that if $ \pi $ is cuspidal then $ \pi $ is irreducible. The theorem is stated in terms of lifting a cusp form, but is still valid for a non-trivial automorphic character.

$ \square $

Next we state some dimension restrictions for first occurrence to take place for a nontrivial quadratic automorphic character.
\begin{lemma}
\label{sec:dimensionRestriction}
Given a globally anisotropic quadratic form $ ( \qmmm , \vmmm )  $ where $ \dim \vmmm = m $ and a $ k $-dimensional subspace $ U_{ k }  \subseteq \vmmm$. Let $ \ximmmnu =  \xi_{ \nu } = \xi_{ \lambda, \epsilon, \nu }  $ be a non-trivial automorphic quadratic character of $ \hmmm = O( \qmmm )( \Fnu )   $ such that

$$ \restr { \xi_{ \lambda, \epsilon, \nu }  } { O( U_{ k }  ) } \equiv 1 $$

If $ k \geq 3 $ then 
$$ \xi_{ \lambda, \epsilon, \nu } \equiv 1 $$

\end{lemma}

\proof

The main point is that the spinor norm is onto $ \F^{ \times } \backslash \left (   \F^{ \times }   \right )^{ 2 }  $ for spaces of 3 or more variables (see Proposition $ 3.1.iv, $ in \cite{kudla1996notes}). Together with the fact that the Hilbert symbol is non-degenerate one can easily show

$$ \restr { \xi_{ \lambda, \epsilon, \nu }  } { O( U_{ k }  ) } \equiv 1 \implies \xi_{ \lambda, \epsilon, \nu } \equiv 1 $$

$ \square $

\begin{cor}[Dimension restrictions]
\label{dimRestrictions}
The theta lift of a non-trivial quadratic automorphic character $ \xi_{ n }  $ of $ \oem  $ from $ \oem  $ to $ \spen $ for is zero if and only if $ m - n \geq 3 $.

\end{cor}

\proof
Follows from Lemma \ref{sec:dimensionRestriction}.
$ \square $

\begin{lemma}[Dimension restrictions]

\label{sec:dimRestrictions}
Suppose $ \adata $ is a globally admissible quadruple such that the corresponding theta lifts are non-zero. Then the dimensions of the quadratic spaces are
 $ 
 dim \left ( V \right ) , dim \left ( V\p \right ) \in  \left \{ n, n + 2 \right \} $.

\end{lemma}
\proof
Technically this is just a alternative formulation of Lemma \ref{sec:dimensionRestriction} and
Corollary \ref{dimRestrictions}.
A key tool in this proof is the global Fourier analysis of the theta integral with respect to the unipotent radical of the Siegel parabolic in $ G $ therefore we defer the proof to Section \ref{section:globalFourier}.

$ \square $




\appendix

\newpage
\section{ Some proofs }
\label{SomeProofs}
In this section we present proofs of results that we chose to defer until later in order not to hinder the general flow of this work.
\subsection{ Proofs from Fourier analysis }
We restate and prove Lemma \ref{sec:globalFourierExpansion}.
\begin{lemma}
Given a quadratic space $ \left ( q, V \right ) $. Let $ \beta\in  \SymF  $. Then for every $ \phi \in S \left (  \voxpr\ad \right )  $
\begin{enumerate}
\item
Suppose $ q $ does not represent $ \beta $ then $ I^{ \phi } _{ V,\xi,\beta } = 0 $.
\item
If $ \beta $ is degenerate then $ I^{ \phi } _{ V,\xi,\beta } = 0 $.
\item
If $ \beta $ is non-degenerate then
the  $ \psi_{ \beta }  $-th Fourier coefficient of the form  $ \fthLift \phi  $ is given by
$$
\boxed{ \thLiftBeta \phi { g }
= 
 \int_{ \nicequotm } \left (  \wni (  g, h ) \phi  \right ) (\vvbeta) \xi  (h)\,dh
 \cdot
 \int_{ \badquotm }  \xi  (h_{ 0 } )\,dh_{ 0 } }
$$
where $ \vvbeta \in  \voxpr_{ \F } \cong V_{ \F } ^{ n }  $ satisfies  $ \gram\vvbeta = \beta $.
\end{enumerate}

\end{lemma}
\proof

First not that all non-degenerate $ \beta $ lie on the same orbit. This follows from a suitable version of Witt's theorem.

Note that changing the order of summation and integration is justified by the fact that the quadratic space $ \left ( q, V \right ) $ is assumed to be globally anisotropic over $\F$, hence by reduction theory  (see \cite{borel2019introduction}) $H =  O \left ( q, V \right )  $ is compact. Moreover the unipotent radical $ N $ of the Siegel parabolic $ P $ of $ G $ forms an abelian group. Therefore by abelian Fourier analysis \cite{rudin2017Fourier} we may expand 
$
\thLiftBeta \phi { \nn g }
$
as a Fourier series along the unipotent radical $
\qFA NN
$
where 
$$
N= 
\Set*{ \nn } { b\in \Sym },
\quad
\nn= \begin{pmatrix} I_n & b \\ 0_n & I_n \end{pmatrix}
$$
For each $ g\in G_{ \fA }  $ we consider the continuous function on 
$$ 
\SymQuot
\to
\fC $$
given by 
$$
b\mapsto \thLiftBeta \phi { \nn g }
$$
\par
It is enough to calculate the Fourier coefficient for $ g = 1 $ since we are proving for a general Schwartz function $ \phi $ and we can always replace $ \phi $ with the Schwartz function
$ 
 \omega_{ \vpsi } \left (  g \right ) 
 \phi
 $.
\par
Also, without loss of generality we may assume that
$ \beta $ is a diagonal matrix since 
we know that there exists $ a \in GL_{ n }  \left ( \F \right )  $ such that 
$$
a\cdot \beta\cdot ^{ t } a = 
diag \left (  b_{ 1 } , \ldots, b_{ n }  \right ) 
$$
for some $ b_{ 1 } , \ldots, b_{ n }  \in \F $
together with the fact that 
$$
\thLiftBeta \phi {  m \left ( a \right )  }
 = 
I^{ \phi } _{ V, \xi, a\cdot \beta\cdot ^{ t } a } 
$$
where $  m \left ( a \right )  = \begin{pmatrix} a & 0 \\ 0 & a^{-1} \end{pmatrix} $.
In other words for the above chosen $ a \in GL_{ n }  \left ( \F \right ) $ we may assume
$$
\thLiftBeta \phi {  m \left ( a \right )  }
 = 
I^{ \phi } _{ V, \xi, diag \left (  b_{ 1 } , \ldots, b_{ n }  \right )  } 
$$
for some $ b_{ 1 } , \ldots, b_{ n }  \in \F  $.
Thus without loss of generality we may assume that $ \beta $ is a diagonal matrix.
Hence
\begin{align*}
\thLiftBeta \phi { 1 } & =
\int_{ \SymQuot } \thetaintni( n(\bp)) \psi_{ -\beta }( n(\bp) ) \,d\bp  \\
  & =  
  \int_{ \SymQuot } \int_{ \qqquotni } \nnsumiF{m} \underbrace{ \left (  \wni ( n( \bp ) , h) \phi    \right ) (\gamma) }_{ \psi( tr \left ( \bp \gram \gamma \right )  )\left (  \wni (  1, h  ) \phi   \right ) (\gamma) } \xi  ( h ) \,dh \cdot \underbrace{ \psi_{-\beta}( n(\bp) ) }_{ \psi(- \tr\left (  b\beta  \right )) } \,d\bp \\
	 &
	  =  \int_{ \qqquotni }
\nnsumiF{m}
	  \left (  \int_{ \SymQuot } \psi( tr\left ( b ( \gram \gamma  -\beta)  \right ))  \,d\bp  \right )\left (  \wni (  1, h ) \phi   \right ) (\gamma) \,dh \\
	 &
	  =  \int_{ \qqquotni } \nnsumiFbeta{m}
\left (  \wni (  1, h ) \phi   \right ) (\gamma) \,dh 
\qquad
 \left (  \star \right ) 
\\
	 & 
	 =  \int_{ \qqquotni } \int_{ \qqquotniii } \left (  \wni (  1, h ) \phi   \right ) ( h_{ 0 } ^{ -1 } \vvbeta) \xi  ( h ) \,dh_{ 0 }  \,dh 
	 \\
	 & =  \int_{ \qqquotni } \int_{ \qqquotniii } \left (  \wni (  1, h_{ 0 }h ) \phi   \right ) ( \vvbeta) \xi  ( h_{ 0 } h ) \,dh_{ 0 }  \,dh 
\com{Automorphy of $ \xi $}
	 \\
	 & =  \int_{ \qqquotniiii } \left (  \wni (  1, h ) \phi   \right ) ( \vvbeta) \xi  ( h ) \,dh \\
	 & = \int_{ \nicequotm } \left (  \wni (  1, h ) \phi  \right ) (\vvbeta) \xi  (h)\,dh \int_{ \badquotm }  \xi  (h_{ 0 } )\,dh_{ 0 } \\
\end{align*}
where $ H^{ \vvbeta }  $ is the stabilizer of  $
 \left (  v_{ 1 } , \ldots, v_{ n }  \right ) =
 \vv^{ \beta }  \in \voxpr_{ \F } \cong V^{ n }_{ \F }   $, that is
$$
H^{ \vvbeta } = 
\Set*{ h \in H } { hv_{ i } = v_{ i } , 1 \leq i \leq n }
$$
Therefore 
$$
\thLiftBeta \phi { 1 }
= 
 \int_{ \nicequotm } \left (  \wni (  1, h ) \phi  \right ) (\vvbeta) \xi  (h)\,dh
 \cdot
 \int_{ \badquotm }  \xi  (h_{ 0 } )\,dh_{ 0 }
$$
hence if we replace $ \phi $ with $ \wni \left (  g, 1 \right )  \left (  \phi \right )  $ we obtain the result
$$
\thLiftBeta \phi { g }
= 
 \int_{ \nicequotm } \left (  \wni (  g, h ) \phi  \right ) (\vvbeta) \xi  (h)\,dh
 \cdot
 \int_{ \badquotm }  \xi  (h_{ 0 } )\,dh_{ 0 }
$$
$ \square $

Next we restate and prove Lemma \ref{sec:FourierExpLemma}.

\begin{lemma}
\label{sec:FourierExpLemmaAppendix}
For $ \qcdata $ such that
the corresponding theta lift $  \Theta \left ( \eta \right )  $ satisfies first occurrence
there exists  $ \vv \in \voxpr_{ \F } $ such that the Fourier expansion of the function 
$$
b\mapsto \thLift \phi {  n \left ( b \right ) g } 
$$
is given by
\begin{align*}
   \thLift \phi { \begin{pmatrix} 1 & b \\ 0 & 1 \end{pmatrix} g }  & =
 \sum_{ \alpha\in M_{ \F }  } W_{ \phi }  \left ( \alpha g \right )
\,\cdot
\psi \left ( b \gram { \alpha \vv } \right ) 
\end{align*}
for every  $ g \in G $ and $ b \in \SymF $.
\par
In particular we have 
$$
\boxed{
\thLift \phi { g } =
 \sum_{ \alpha\in M_{ \F }  } W_{ \phi }  \left ( \alpha g \right )
 }
$$
\end{lemma}
\proof
By the above Fourier coefficient computation in Lemma \ref{sec:globalFourierExpansion} we have
\begin{align*}
\thLift \phi { \begin{pmatrix} 1 & b \\ 0 & 1 \end{pmatrix} g }  & =   \sum_{ \betarepn  } 
 \int_{ \nicequotm } \left (  \wni (  1, h ) \phi  \right ) (\vvbeta) \xi  (h)\,dh
  \cdot \psi_{\beta}( n(b) )  \\
	 & =  \sum_{ \betarepn  } 
 \int_{ \nicequotm } \left (  \wni (  1, h ) \phi  \right ) (\vvbeta) \xi  (h)\,dh
	  \cdot 
	  \psi(  tr \left ( b\beta \right )  )  \\
	 & =  \sum_{ \betarepn  }    \int_{ \nicequotm } \left (  \weilVn( g, h ) \phi    \right )(\vvbeta) \xi   ( h ) \,dh \cdot \psi(  tr \left ( b\beta \right )  )
\end{align*}
Where we are running over the image of the Gram matrix of $ \left ( q, V \right ) $ in the sum
$ \gram\cdot\,\colon \voxpr_{ \F } \cong V^{ n } _{ \F } \to \SymF $.
\par
If we choose  $ \vv \in \voxpr_{ \F }$ such that  $  \gram\vv  = \beta $ for some  $ \beta \in \SymF $ with $  det \left ( \beta \right ) \neq 0 $ that is represented by  $ q $ then 
\begin{align*}
   \thLift \phi { \begin{pmatrix} 1 & b \\ 0 & 1 \end{pmatrix} g }  & =    \sum_{ \alphat\in\glnf   }    \int_{ \nicequotmv } \left (  \weilVn( g, h ) \phi    \right )( \alphat\vv) \xi   ( h ) \,dh \cdot \psi(  tr \left ( bq  ( \alphat\vv ) \right )   )  \\
	 &   =  \sum_{ \alphat\in\glnf } W_{ \phi }  \left ( m \left ( \alphat \right ) g \right ) \cdot \psi(  tr \left ( bq  ( \alphat\vv ) \right )  )
\end{align*}
In particular if we take $ b = 0_{ n }  $  then this gives us
\begin{align*}
   \thLift \phi { g }  & =  \sum_{ \alphat\in\glnf } W_{ \phi }  \left ( m \left ( \alphat \right ) g \right ) \\
& =  \sum_{ \alpha\in M_{ \F }  } W_{ \phi }  \left ( \alpha g \right )
\end{align*}
$ \square $
We restate and prove Theorem \ref{sec:zeroLifts}.

\begin{theorem}
Given
$ \qcdata $
for
$  \left (  q, V \right )  $ and $ \xi \,\colon H_\fA \to \fC^\times $ as defined in this work, i.e. $ \xi $ is trivial on $ H^{ \vvbeta } $ for non-degenerate $ \beta \in \SymF $.
 Then 
$$
\pi_{ k } = 
\liftc{ k }
= 
0
$$
for every $ 1 \leq k < \N $.
\end{theorem}
\proof
We will prove that for every $ \beta \in \SymF $ the Fourier coefficient of the $ k $-th dual pair in the tower is zero. Namely we will show that 
$$
\fthLiftBeta \phi \equiv 0
$$
By Lemma \ref{sec:globalFourierExpansion} we have
$$
\thLiftBeta \phi { g }
= 
 \int_{ \nicequotm } \left (  \wni (  g, h ) \phi  \right ) (\vvbeta) \xi  (h)\,dh
 \cdot
 \int_{ \badquotm }  \xi  (h_{ 0 } )\,dh_{ 0 }
$$
and one can easily show that the non-trivial character $ \xi $ is non-trivial on $ H^{ \vvbeta } $, therefore
$$
\int_{ \badquotm }  \xi  (h_{ 0 } )\,dh_{ 0 } = 0
$$
hence we have shown that for every $ \beta\in \SymF $ we have 
$$
\fthLiftBeta \phi \equiv 0
$$
hence 
$$
\fthLift \phi \equiv 0
$$
hence 
$$
\pi_{ k } = 
\liftc{ k }
= 
0
$$
for every $ 1 \leq k < \N $ where $ \N = n$ for the dual pair $  O \left ( q_n, V_n \right ) \times G $
and $ \N = n + 2 $ for the dual pair $  O \left ( q_{ n+2 } , V_{ n+2 }  \right ) \times G $.
(By Corollary \ref{dimRestrictions} these are the only possible dual pairs of interest).
\par
Just to clarify the remark that "one can easily show that" notice that if $  rank \left ( \beta \right )  = r $ then $ \beta $ is of the form 
$$
\beta = diag \left (  b_{ 1 } , \ldots, b_{ r },\underbrace{ 0,0 , \ldots, 0 }_{ n-r\tn{ times} } \right ) 
$$
therefore in step $  \left (  \star \right )  $ of the above computation the condition that the Gram matrix of $ \gamma $  is equal to $ \beta $, i.e. the condition that 
$$
\gram \gamma = \beta
$$
means that there exists 
$$
\gamma =  \left (  \gamma_{ 1 } , \ldots, \gamma_{ n }  \right )  \in \voxpr_{ \F } \cong
V^{ n } _{ \F } 
$$
such that
$$
\begin{cases} b_{ q }  \left (  \gamma_{ i } ,\gamma_{ j }  \right )  = 0,&\quad\tn{ for all }1 \leq  i \neq  j \leq  n \\ \\  
q \left ( \gamma_i \right )  = b_{ q }  \left (  \gamma_{ i } ,\gamma_{ i }  \right )  = 0,&\quad\tn{ for all }r + 1 \leq  i \leq n \\ \\ 
 q \left ( \gamma_i \right )  =  b_{ q }  \left (  \gamma_{ i } ,\gamma_{ i }  \right )  = b_{ i } ,&\quad\tn{ for all }1 \leq  i \leq  r 
 \end{cases}
$$
Notice that if 
$$
\vvbeta = 
 \left (  
 w_{ 1 } , \ldots, w_{ r } 
 ,
\underbrace{  0 , \ldots, 0 }_{ n - r \tn{ times} }
  \right ) 
   \in  \left (  \vox \right ) _{ \F } \cong V_{ \F } ^{ n } 
$$
and if we let 
$$
V^{ 0 }  = \Span_{ \F } 
\left \{ 
 w_{ 1 } , \ldots, w_{ r } 
\right \}
$$
and complete $  \left \{  w_{ 1 } , \ldots, w_{ r } \right \} $ to a an orthogonal basis of $ V $ given by 
$$
\left \{  w_{ 1 } , \ldots, w_{ r }
,
w_{ r+1 } , \ldots, w_{ n } 
\right \} 
$$
and let 
$$
V_{ 0 }  = \Span_{ \F } 
\left \{ 
w_{ r+1 } , \ldots, w_{ n } 
\right \}
$$
then we have $ V = V^{ 0 } \oplus V_{ 0 }  $ and we have an isomorphism 
$ H^{ \vvbeta }\cong  O \left ( q, V_0 \right )  $
given by the restriction map 
$$
g \in H^{ \vvbeta }\mapsto
\restr { g } { V_0 }
$$
Since $ \xi $ is a non-trivial character when restricted to $  O \left ( q, V_0 \right ) \ad $ we get that 
$$
\int_{ \badquotm }  \xi  (h_{ 0 } )\,dh_{ 0 } = 0
$$
as claimed earlier.
$ \square $

We restate and prove Theorem \ref{sec:nonzeroLift}.

\begin{theorem}
\label{nonvanish}
Given
$ \qcdata $
such that $ \restr { \xi } { \vvbeta } \equiv 1$ for some non-degenerate symmetric matrix $ \beta $.
Then the corresponding global theta lift is non-zero.
\end{theorem}

\proof
Without loss of generality we may take $ \beta = \betao $ to be the matrix representing the quadratic space $ \left ( q, V \right ) $.
The matrix $ \betao $ is represented by the form $ q $
We will show that
$$
\fthLiftBetao \phi\neq 0
$$
By Lemma \ref{sec:globalFourierExpansion} we have 
$$
\thLiftBetao \phi { g }
= 
 \int_{ \nicequotOm } \left (  \wni (  g, h ) \phi  \right ) (\vvbetao) \xi  (h)\,dh
 \cdot
 \int_{ \badquotOm }  \xi  (h_{ 0 } )\,dh_{ 0 } 
$$
and since $ \xi $ is trivial on the stabilizer 
$ H^{ \vvbetao } $
we have 
$$
\thLiftBetao \phi { g }
= 
 \int_{ \nicequotOm } \left (  \wni (  g, h ) \phi  \right ) (\vvbetao) \xi  (h)\,dh
$$
Assume $ \phi = \otimes \phi_{ \nu } $ is decomposable. Then 
$$
\int_{ \nicequotOm } \left (  \wni (  g, h ) \phi  \right ) (\vvbetao) \xi  (h)\,dh
= 
\prodnu
\int_{ \nicequotOnnu } \left (  \wni (  g, h ) \phi_{ \nu } \right ) (\vvbetao) \xi  (h)\,dh
$$
Set $ g = 1 $.
Choose a small neighborhood $ U_{ \nu }  $ of $ I $ such that $  \xi \left ( h \right ) = 1 $ for every $ h\in U_{ \nu } $. Furthermore, we will restrict the $ \phi_{ \nu }  $ such that 
$$
 supp \left ( \phi_\nu \right )  \subseteq H_{ \nu } 
$$
for every $ v\in S $.
Moreover denote 
\begin{enumerate}
\item
$ \cpfInf $ the set of infinite places of $ \F $.
\item
$ \cpffin $ the set of finite places of $ \F $.
\item
$ \cpf = \cpfInf\cup \cpffin $ the set of place of $ \F $.
\end{enumerate}
and set $ \phi_{ \nu } = 1_{ \V_\nu } $ for every $ \nu\not\in S $. In addition 
\begin{enumerate}
\item
For $ \nu\in \cpfInf\cap S $, let $ \phi_{ \nu } $ be any positive Scwartz function supported in $ U_{ \nu } $. 
\item
For $ \nu\in  \cppf \cap S $, let $ \phi_{ \nu } = 1_{ U_\nu\cdot \vvbeta }  $.
\end{enumerate}
Hence 
\begin{align*}
 \prodnu
\int_{ \nicequotOnnu } \phi_{ \nu } ( h^{ -1 } \vvbetao) \xi  (h)\,dh
& = 
\prodnuIn S
\int_{ \nicequotOnnu } \phi_{ \nu } ( h^{ -1 } \vvbetao) \xi  (h)\,dh
\\
&
\quad
\cdot 
\underbrace{ \prodnuNotIn S
\int_{ \nicequotOnnu } \phi_{ \nu } ( h^{ -1 } \vvbetao) \xi  (h)\,dh }_{ = 1 }
 \\
&   = 
\prodnuIn S
\int_{ \nicequotOnnu } \phi_{ \nu } ( h^{ -1 } \vvbetao) \xi  (h)\,dh
 \\
&   = 
\prodnuIn \si
\int_{ \nicequotOnnu } \phi_{ \nu } ( h^{ -1 } \vvbetao) \xi  (h)\,dh
 \\
& 
\quad
\underbrace{ \cdot
\prodnuIn \sii
\int_{ \nicequotOnnu } \phi_{ \nu } ( h^{ -1 } \vvbetao) \xi  (h)\,dh }_{ = 1 }
 \\
&   = 
\prodnuIn \si
\int_{ \nicequotOnnu } \phi_{ \nu } ( h^{ -1 } \vvbetao) \xi  (h)\,dh
 \\
& >
0
\end{align*}
for a suitable choice of normalized measures.
$ \square $

\subsection{ Proof of Theorem \ref{sec:unfoldGlobalThm} }

For  $  Re \left ( s \right ) > \auxRange $ we have
$$
B_{ s }
\left ( \phi, \phip \right ) 
= 
\glblZetaConst
\cdot
\int_{ \ga K } \int_{ \ga M }  W_{ \phi }  \left (  mk \right )  W\p_{ \phip }  \left (  mk \right ) \delb m \fPsi{ mk }{ s } \,dm \,dk
$$
where our data $ \adata $ is given by a globally admissible quadruple.

\proof
Let $ P  \subseteq G$ be the Siegel parabolic subgroup of $ G $.
Let $ P = MN $ be the Levi decomposition of $ P $ where $ M $ is the Levi subgroup of $ P $ and $ N $ is the unipotent radical in $ P $.
\par

Let us unfold  $ B_{ s } $ where in the following we will assume that $ Re \left ( s \right ) > \auxRange $ so that $ \cE^{ \star } $ is absolutely convergent.
Hence 
\begin{align*}
\glblZetaConst^{ -1 } 
 \cdot B_{ s } \left ( \phi, \phip \right )  & = 
\int_{ \aq G } 
\thLift \phi g
\overline{ \thLiftp \phip g }
\naeis g s
\,dg 
\\
&  =
\glblZetaConst^{ -1 } 
\int_{ \aq G }  \thLift \phi g \overline{ \thLiftp \phip g } \cdot 
\glblZetaConst
\aeis g s \,dg  \\
\\
& = 
\int_{ \aq G } 
\thLift \phi g \overline{ \thLiftp \phip g }
\sum_{ \gamma \in \rquot P G } 
\fPsi{ \gamma g }{ s }
\,dg
\\
&  =
\int_{ \aq G }
\sum_{ \gamma \in \rquot P G } 
\thLift \phi g \overline{ \thLiftp \phip g }
\fPsi{ \gamma g }{ s }
\,dg
\\
&  =
\int_{ \aq G }
\sum_{ \gamma \in \rquot P G } 
\thLift \phi { \gamma g } \overline{ \thLiftp \phip { \gamma g } }
\fPsi{ \gamma g }{ s }
\,dg
\com{ Automorphy }
\\
&  =
\int_{ \qFA P G }
\thLift \phi { g } \overline{ \thLiftp \phip { g } }
\fPsi{ g }{ s }
\,dg
\\
&  =
\int_{ \qAA P G }
\int_{ \qFA P P }
\thLift \phi { pg } \overline{ \thLiftp \phip { pg } }
\fPsi{ pg }{ s }
\,dp
\,dg
\\
&  =
\int_{ \qAA P G }
\int_{ \qFA MM}
\int_{ \qFA N N }
\thLift \phi { mng } \overline{ \thLiftp \phip { mng } }
\fPsi{ mng }{ s }
\,dn
\,dm
\,dg
\\
\end{align*}
Note however that
\begin{align*}
   \fPsi{ mng }{ s }  & =  \fPsi{ \underbrace{ mnm^{ -1 } }_{ \in N } \cdot mg }{ s }  \\
	 &   =  \fPsi{ mg }{ s }    
\end{align*}
Hence continuing our computation thus far we have obtained
\begin{align*}
\glblZetaConst
^{ -1 } 
 \cdot B_{ s } \left ( \phi, \phip \right )  & = 
\int_{ \qAA P G }
\int_{ \qFA M M }
\int_{ \qFA N N }
\thLift \phi { mng } \overline{ \thLiftp \phip { mng } }
\,dn
\fPsi{ mg }{ s }
\,dm
\,dg
\end{align*}
Next we will substitute the Fourier expansion of 
\begin{align*}
   \thLift \phi g  & = \sum_{ \alpha \in M_\F }  \Hintegral \left ( \wfunc  \left ( g, h \right )  \phi \right ) \left (   \alpha \vv \right )  \xi \left ( h \right )  \,dh  \\
	 &   =  \sum_{ \alpha \in M_\F } 
	  W_{ \phi } \left ( \alpha g \right )
\end{align*}

to obtain

\begin{align*}
\glblZetaConst
^{ -1 } 
\cdot B_{ s } \left ( \phi, \phip \right )  & = 
\int_{ \qAA P G }
\int_{ \qFA M M }
\int_{ \qFA N N }
  \sum_{ \alpha \in M_\F }  \constiINV
   W_{ \phi } \left ( \alpha mng \right )
\overline{ \thLiftp \phip { mng } }
\,dn
\fPsi{ mg }{ s }
\,dm
\,dg
\\
&  
\nquad
=
\int_{ \qAA P G }
\int_{ \qFA M M }
\int_{ \qFA N N }
  \sum_{ \alpha \in M_\F } 
  W_{ \phi } \left ( \alpha mng \right )
\overline{ \thLiftp \phip { \alpha mng } }
\,dn
\fPsi{ \alpha mg }{ s }
\,dm
\,dg
{ \tiny 
\,\,\,\left [ \tn{Automorphy} \right ]
  }
\\
& 
\nquad
 =  
\int_{ \qAA P G }
\int_{ \ga M }
\int_{ \qFA N N }
 W_{ \phi }  \left ( 
mng
\right )
\overline{ \thLiftp \phip { mng } }
\,dn
\fPsi{ mg }{ s }
\,dm
\,dg
\end{align*}
Now make a change of variables 
$$
n\mapsto m^{ -1 } nm
$$
which multiplies the meaure  $ dn $ by 
 $ \delb m $ to obtain
\begin{align*}
\glblZetaConst
^{ -1 } 
 \cdot B_{ s } \left ( \phi, \phip \right )  & = 
\int_{ \qAA P G }
\int_{ \ga M }
\int_{ \qFA N N }
 W_{ \phi }  \left ( 
nmg
\right )
\delb m
\cdot
\overline{ \thLiftp \phip { nmg } }
 \,dn
\fPsi{ mg }{ s }
\,dm
\,dg
\\
&  =
\int_{ \qAA P G }
\int_{ \ga M }
 W_{ \phi }  \left ( 
mg
\right )
\int_{ \qFA N N }
\overline{ \thLiftp \phip { nmg } }
\psi \left ( n \right ) 
 \,dn
 \cdot
\delb m
\fPsi{ mg }{ s }
\,dm
\,dg
\\
&  =
\int_{ \qAA P G }
\int_{ \ga M }
 W_{ \phi }  \left ( 
mg
\right )
 W\p_{ \phip }  \left ( 
mg
\right )
\delb m
\fPsi{ mg }{ s }
\,dm
\,dg
\end{align*}
Note that conjugation does not appear on $  W\p_{ \phip } $ because it was defined in terms of the character $ \psi^{ -1 } $.
We note that by the Iwasawa decomposition
$$
G_{ \fA } = P_{ A } \cdot K
$$
therefore
\begin{align*}
   \glblZetaConst ^{ -1 }  \cdot B_{ s } \left ( \phi, \phip \right )  & =  \int_{ \ga K } \int_{ \ga M }  W_{ \phi }  \left (  mk \right )  W\p_{ \phip }  \left (  mk \right ) \delb m \fPsi{ mk }{ s } \,dm \,dk   
\end{align*}
hence
\begin{align*}
B_{ s } \left ( \phi, \phip \right )  & = 
\glblZetaConst
\cdot
 \int_{ \ga K } \int_{ \ga M }  W_{ \phi }  \left (  mk \right )  W\p_{ \phip }  \left (  mk \right ) \delb m \fPsi{ mk }{ s } \,dm \,dk   
\end{align*}

\par

By the unramified computation for \( n = 1 \) the unramified local pairing is given by an Euler factor. For \( n>1 \) it is not clear if this is true since we have not done the unramified computation in this case. 
From general results the bad local factoRs \( B_{ s,\nu }  \left (  \phi_{ \nu } ,\phi_{ \nu } \p \right )  \).
Moreover for certain local data
\( \phi_{ \nu }  \)  and \( \phi_{ \nu } \p \) we can choose these with small enough support such that the pairing is non-zero.
\par
Note that the product 
\[
L^{ s }  \left (  s \right )  = 
\prod_{ v\not\in S } 
L_{ \nu }  \left ( s \right )
\]
has an analytic continuation with respect to \( s \) in a neighborhood of \(  s = \rho_{ 1 }  = 1 \)  with a simple pole at \( s = \rho_{ 1 }  = 1 \).
This follows from the fact that we can write
\[
L_{ \nu }  \left ( s \right )
 = 
b_{ n } ^{ s }  \left (  
\rho_{ 1 } 
 \right ) 
^{ -1 } 
 \left (  
\prod_{ v\not\in S } 
B_{ s=\rho_1, \nu } 
\left (  \phi_{ \nu } ,\phi_{ \nu } \p \right ) 
  \right ) 
  ^{ -1 } 
\cdot
Res_{ s=\rho_1 } 
L^{ S }  \left ( s \right ) 
\]

$ \square $

\newpage

\section{ The Weil representation and the Theta correspondence }\label{sec:weil}
\label{sec:weilRepThetaCorre}
In this section we define the Heisenberg group and then we state the Stone-von Neumann theorem. We then construct a family of representations of the Heisenberg group such that every two representations are isomorphic. 
From Schur's lemma we obtain intertwining operators between every two isomorphic representations. This will give us a projective representation of the symplectic group which will lead us to the definition of the Weil representation of the metaplectic group. We also state basic properties of the Weil representation that will be used later in this work.
\par
In this section we closely follow the description in \cite{kudla1996notes} and \cite{moeglin1987correspondances}.
We wish to define the local and global theta correspondence. To do this we must describe the Heisenberg group and the Stone von Neumann theorem in order to define the Weil representation. Furthermore we describe reductive dual pairs and the Schr\"{o}dinger model of the Weil representation.
\subsection{ Heisenberg group, Stone von Neumann theorem, and the Weil representation }\label{sec:HeisStoneWeilRep}
In this section $ \F $ will denote a non-Archimedean local field  such that the residue characteristic of $ \F $ is not equal to two and $  \left ( \fW, \left < , \right > \right )  $ will denote a non-degenerate symplectic $ \F $-vector space. Fix a nontrivial additive character $ \psi \,\colon \F \to \fC $.
\par
 In order to define the Weil  representation we will first briefly recall the definition of the Heisenberg group and the Stone von Neumann theorem.  This will give rise to a projective representation of the symplectic group  and this in turn will allow us to define the Weil  representation of the metaplectic group determined by a fixed nontrivial additive character $ \psi \,\colon \F \to \fC $.



\begin{defi}[Heisenberg group]
Define the \bbb{Heisenberg group} as the set  $  \heisww = \fW\oplus \F $
\label{notation:heisGroup}
 together with the group operation given by 
$$
 \left (  w,t \right ) 
 \circ
 \left (  w\p ,t\p \right ) 
 = 
  \left (  
  w+ w\p,
  t+ t\p+\frac{ 1 } { 2 }\left < w, w\p  \right >
   \right ) 
$$
where
$  \left (  w,t \right ) 
,
 \left (  w\p ,t\p \right )  \in \heisww $.
\end{defi}
 \par
We have a natural action of  $  \spww  $ on  $ \heisww $ given by 
 $  \left (  w,t \right )^{ g } =  \left (  wg, t \right )   $ where  $ g \in \spww, \left (  w,t \right )\in \heisww $.
 \begin{remark}
The center of  $ \heisww $ is given by
\label{notation:ZheisGroup}
$$
 Z \left ( \heisww \right ) = 
\Set*{ \left ( 0, t \right )  } { t\in \F }\cong \F
$$
Moreover,  $ \spww $ acts trivially on  $  Z \left ( \heisww \right ) $.
\end{remark}

\begin{theorem}[Stone-von Neumann]
Up to isomorphism, there is a unique smooth irreducible representation
\label{notation:RepHeis}
  $ \rep \heisww \rhop S  $  with central character  $ \psi $. Namely 
  $$
\rhop \left ( 0, t \right ) =  \psi \left ( t \right ) \cdot id_{ S } , \forall t \in \F
$$
\end{theorem}
\proof
See \cite[pp. 28-31]{moeglin1987correspondances}.
$ \square $

Note that for every  $ g\in \spww $  we can define a representation
 $ \rep \heisww \rhopg S $
\label{notation:RepHeisg}
  given by the action
 $$
 \rhopg \left ( h \right ) =  \rhop \left ( h^g \right ) 
$$
 and it is important to note that this new irreducible representation $ \rep \heisww \rhopg S $ has the same central character
  $ \psi \,\colon \F \to \fC $ as the representation $ \rep \heisww \rhop S $. Therefore by the Stone-von Neumann theorem the representations $ \rep \heisww \rhopg S $ are all isomorphic to  $ \rep \heisww \rhop S $ therefore  for every $ g\in \spww $ there exists an automorphism 
  \label{notation:HeisIntertwine}
   $ 
 A \left ( g \right )  \,\colon S \to S
    $
 satisfying
 $$
\weilcond
$$
Note that these automorphisms are not unique however
 by Schur's Lemma they are uniquely determined up to a scalar in  $ \fC^{ \times }  $  therefore the automorphisms  $  A \left ( g_1 \right )  A \left ( g_2 \right )  $ and  $  A \left ( g_1 g_2 \right ) $  are  equal up to a scalar.  Therefore we obtain a homomorphism
 $$
\begin{cases}
   \spww\to GL \left ( S \right ) / \fC^{ \times } \\
   g \mapsto  A \left ( g \right ) \\
\end{cases}
$$
\label{notation:centralExtSpW}
 One can then define a central extension of $ \tspwwp $, that is 
 $$
1
\to
\fC^{ \times } 
\to
\tspwwp
\to
\spww
\to
1
$$
 where we define
 $$
\tspwwp = 
\Set*{  \left ( g, A \left ( g \right )  \right ) 
 \in \spww\times  GL \left ( S \right )  } { 
 \weilcond
 ,
 h \in \heisww
  }
$$
 this gives a representation  $ \rep \tspwwp { \omega_{ \psi }  } S $. We will state the following results without proof.  For proofs one can consult \cite{kudla1996notes}.
 \label{sec:defWeili}
\begin{theorem}[Definition and properties of metaplectic group]
\hspace{2em}
\label{sec:defWeil}
\begin{enumerate}
   \item The central extension $ \tspwwp $  does not depend on the character  $ \psi $. 
   \item Therefore we can define the metaplectic group \label{notation:metGroup}$  Mp \left ( \fW \right ) = \tspwwp $. 
   \item The metaplectic group is isomorphic to the extension which is obtained from a nontrivial twofold topological central extension $$1\to\mu_{ 2 } \to Mp \left ( \fW \right ) ^{ (2) } \to\spww\to1$$ 
   \item Therefore we can define a representation $ \rep {  Mp \left ( \fW \right ) } { \omega_{ \psi } } S $.
\end{enumerate}
\end{theorem}
\begin{remark}
For details see \cite[p. 1-6]{kudla1996notes} 
\end{remark}
%
Recall from Section \ref{addChar} that for $ a \in \F^{ \times } $ we denote 
$
\psi_{ a }  \left (  x \right )  =  \psi \left ( ax \right ) 
$.
\begin{lemma}[Properties of the Weil representation]
\hspace{2em}
\label{sec:contraOfWeilRep}
\begin{enumerate}
   \item  We have an isomorphism \label{notation:weilPsia}$$ \repwpsi{ a } \cong \repwpsi{ ab^{ 2 } } $$ for every $ a,b\in \F^{ \times }  $. 
   \item More generally $$ \repwpsi{ a } \not\cong \repwpsi{ ab^{ 2 }  } $$  for every  $ a\F^{ \times, 2 } \neq b\F^{ \times, 2 } $ in  $  \F^{ \times } /  \F^{ \times, 2 } $. 
   \item 
\label{notation:contraWeilPsia}
   $ \wwpsi a^{ \vee } \cong \wwpsi{ -a }  $ for every  $ a\in \F^{ \times } $
where $ \wwpsi a^{ \vee } $ denotes the contragradient of the Weil Representation with respect to $ \psi_a $.
\end{enumerate}
\end{lemma}

\begin{remark}
Part $ 3 $ of the previous lemma  will be particularly important later when we take the contragradient  of a certain representation of coinvariants. More precisely, it will be used in the transition from local morphisms between theta lifts to local pairings of theta lifts.
\end{remark}

\begin{remark}
 In this work we will primarily use a very specific model of the Weil representation which is referred to as the Schr\"{o}dinger model. We will present precise formulae when restricting to a certain reductive dual pair in the following section.
For a general exposition of models of the Weil representation see \cite[Chapter $ 1 $]{moeglin1987correspondances}.
\end{remark}


\subsection{ Reductive Dual Pairs }
\label{sec:LTCandRedDualPairs}

In this section we will define reductive dual pairs $ H\times G $.
Such pairs will appear in the local and global theta correspondence. Namely, we will "lift" representations from $ H $ to $ G $.
\label{sec:redDualPair}
If  $ H < K $ is a subgroup of  $ K $ then define 
$$
 Cent_K \left ( H \right ) = 
\Set*{  k\in K } { kh = hk, \forall h \in H }
$$
\label{notation:centralizer}
\begin{defi}[Mutual commutants]
A pair of subgroups  $ H $ and  $ G $ of  $ K $ are said to
be mutual \bbb{commutants} if 
 $  Cent_K \left ( H \right ) = G $  and
 $  Cent_K \left ( G \right ) = H $.

\end{defi}
\begin{defi}[Reductive dual pair]
\label{sec:ReductiveDualPairs}
A \bbb{reductive dual pair} $  \left (  H  , G \right )  $ in  $ Sp \left (  \fW \right )  $ is a pair of reductive subgroups  $ H  , G \subseteq Sp \left (  \fW \right )  $ such that  $ H   $ and  $ G $  are mutual commutants of  $ Sp \left (  \fW \right )  $,  i.e.
$$
\begin{cases}
   Cent_{ Sp(\fW) }  \left (  H   \right )  = G\\
   Cent_{ Sp(\fW) }  \left (  G \right )  = H  \\
\end{cases}
$$\label{fm:cent}
\end{defi}

%

\begin{remark}
We will sometimes refer to a reductive dual pair $ 
H\times G
  $
simply as a dual pair.
\end{remark}

\subsubsection{ Example }
Let  
$ \left ( q, V \right )  $
be a finite dimensional quadratic space over $ \F $ with 
an associated nondegenerate symmetric bilinear form 
$$
 b_{ q } \left ( \cdot , \cdot \right )  \,\colon V \times V \to \F
$$
Let 
$
H = O \left ( q, V \right ) 
= 
\Set*{ h \in  GL \left ( V \right )  } { 
b_{ q }  \left ( hx , hy \right ) = b_{ q }  \left ( x , y \right ) 
 \tn{ for every }
 x,y \in V
 }
$
be the isometry group of the quadratic space $ \left ( q, V \right )  $. Likewise, let $  \left ( W, \left <  , \right > \right )  $ be a finite dimensional  $ \F $-vector space with 
a nondegenerate skew-symmetric bilinear form 
$$
 \left < \cdot , \cdot \right >  \,\colon W \times W \to \F
$$
Let 
$
G = Sp \left ( W \right ) 
= 
\Set*{ g \in  GL \left ( W \right )  } { 
 \left < xg , yg \right > =  \left < x , y \right > 
 \tn{ for every }
 x,y \in W
 }
$
be the isometry group of the symplectic space $  \left ( W, \left <  , \right > \right )  $. Then we can define a nondegenerate $ \F $-skew-symmetric bilinear form on the space  $ \fW = V \otimes _{ \F } W $ given by 
\label{notation:TensPairing}
$$
\begin{cases}
   \left < \left <  \cdot , \cdot \right >  \right > \,\colon \fW \times \fW \to \F \\
   \left < \left < \ts 1,\ts 2  \right >  \right >  =  b_{ q } \left (  x_{ 1 } ,x_{ 2 }  \right ) \cdot \left < x_{ 2 } ,y_{ 2 } \right > \\
\end{cases}
$$
Then there is a natural map 
$$
\begin{cases}
   H\times G\to Sp \left ( \fW \right ) \\
  H\times G\mapsto h \otimes g \\
\end{cases}
$$


\subsection{ The local theta correspondence}
\label{sec:localThetaAndHoweDuality}
In this section we describe the local theta correspondence.
\par
Let  $ \F $ be a local field. Let  $ H $ be an orthogonal group and let $ G $ be the symplectic or metaplectic group.

\begin{defi}[The local theta correspondence]

\label{notation:LocalThetaLift}
The local theta correspondence is a function 
$$
\begin{cases}
\Theta = \Theta_{ \psi }  = \Theta_{ V,W,\psi }  \,\colon Rep ( \tH ) \to Rep ( \tG ) \\
   \Theta_{ V,W,\psi }  \left ( \pi \right ) = \left (  \omega_{ \psi } \otimes \widetilde{ \pi }  \right )_{ \tH } \\
\end{cases}
$$

\end{defi}\label{notation:bigTheta}
This is referred to as the "big theta" lift in the literature.
Note that  $  \Theta_{ V,W,\psi }  \left ( \pi \right ) $ is not necessarily irreducible.
One can define a map \label{notation:littleTheta}$ \theta = \theta_{ \psi } = \theta_{ V,W,\psi } \,\colon Irr ( \tH ) \to Irr( \tG ) $
referred to as the "little theta" lift. This is the content of part two of the following theorem.

\begin{theorem}[Properties of the local theta correspondence]
If  $ \F $ has residue characteristic not equal to $ 2 $ then for any irreducible admissible representation  $ \pi $ of  $ \tH $ we have the following
\begin{enumerate}
   \item Either $  \Theta_{ \psi }  \left ( \pi \right ) = 0 $  or $  \Theta_{ \psi }  \left ( \pi \right )  $  is an admissible representation of $ \tG $ of finite length. 
   \item If  $  \Theta_{ \psi }  \left ( \pi \right ) \neq 0 $, then there exists a unique $ \tG $  invariant submodule $  \Theta_{ \psi } ^0 \left ( \pi \right )  $\label{notation:invarSubTheta} of  $  \Theta_{ \psi } \left ( \pi \right )  $  such that$$ \theta_{ \psi } \left ( \pi \right ) =  \Theta_{ \psi } \left ( \pi \right ) /  \Theta_{ \psi }^0 \left ( \pi \right ) $$ is irreducible.  The irreducible admissible representation $  \theta_{ \psi } \left ( \pi \right )  $ of  $ \tG $  is uniquely determined by $ \pi $. If  $  \Theta_{ \psi } \left ( \pi \right ) = 0 $  then we define $  \theta_{ \psi } \left ( \pi \right ) = 0 $. 
   \item If $  \theta_{ \psi } \left ( \pi_1 \right ) \cong \theta_{ \psi } \left ( \pi_2 \right ) \neq 0 $ then  $ \pi_{ 1 } \cong \pi_{ 2 }  $. 
\end{enumerate}
\end{theorem}


\subsection{ The Schr\"{o}dinger model of the Weil representation }
\label{sec:ScrodModel}
Consider the reductive dual pair $ H\times G $  embedded inside of  $ \spww $ where $ \fW = V \otimes _{ \F } W $.
\label{notation:nontriviqlladdchar}
For a fixed nontrivial additive character  $ \psi \,\colon \F \to \fC $  we now give explicit formulas for the Weil  representation
\label{notation:Schrod}
$$
\omega_{ \vpsi } \defeq
\restr { \omega_{ \psi }  } { \widetilde{ Sp } \left (  W \right )  \left (  \fA \right ) \otimes O \left (  V \right ) \left (  \fA \right )   }
$$

\label{section:weil} 
Let  $ W = X\oplus Y $ be a complete polarization of  $ W $ over  $ \F $. Then  $ \omega_{ \vpsi } $ acts on 
$
S \left (  \vox  \right )
$
because  $ W\otimes V \cong  \left (  \vox  \right ) \oplus  \left (  Y\otimes V \right )  $ is a complete polarization of the symplectic space $  \left (  W\otimes V, \doubleTri \right )  $.
Note that $ \vox \cong V^{ n } $.

\begin{lemma}[Weil representation restricted to a reductive dual pair]
\label{lemma:weilFormulas} 
Assume  $ \dim_{ \F } W = 2n $ and  $ \dim_{ \F } V = m $. Then for  $ \phi  \in S \left (  \vox  \right ) , v 
 \in  \left (  \vox  \right )  \ad  $
we have the following formulas
\begin{enumerate}
	\item 
The action of the Levi part is given by
	\begin{align*}
   \left (  \omega_{ \vpsi }  \left (  \ma, 1 \right ) \phi \right )  \left (  v \right ) & = 
 \abs{ \det \left (  a \right ) } ^{ \frac{ m } { 2 } }
\widetilde{ \chi }_{ V }  \left ( a \right )
     \phi \left (  va \right )
     \label{eq1}
\end{align*}
where  $ a \in GL_{ n }  \left (  \fA \right )  $ and  $ d = d( V )  $ denotes the discriminant of the quadratic space  $  \left (  q, V \right )  $. Also recall that 
$$
\widetilde{ \chi }_{ V }  \left ( a \right )
 = 
   \gammafactor
    \left (  \det \left (  a \right ),  \left (  -1 \right )^{ \frac{ m \left (  m-1 \right )   } { 2 } }d( V ) \right ) 
$$
for $ m $ odd and 
$$
\widetilde{ \chi }_{ V }  \left ( a \right )
 = 
    \left (  \det \left (  a \right ),  \left (  -1 \right )^{ \frac{ m \left (  m-1 \right )   } { 2 } }d( V ) \right ) 
$$
for $ m $ even.
\item 
The action of the unipotent radical is given by
\begin{align*}
   \left (  \omega_{ \vpsi }  \left (  \nb, 1 \right ) \phi \right )  \left (  v \right ) & = \psi \left ( 
   \tn{tr} \left ( b \gram v \right )  \right ) \cdot \phi( v ) \\
\end{align*}
where  $ b  \in Sym_{ n }  \left (  \fA \right )  $,  $ b_{ q } $ is the symmetric bilinear associated to the quadratic form $ q $ and  $ \gram v =  \left (  b_{ q }  \left (  v_{ i } , v_{ j }  \right )  \right )_{ i, j }  $ is the Gram matrix of $ q $.
\item
The action of the long Weyl element  $ \omega_{ 0 } 
 $ is given, up to a constant, by the Fourier transform, more precisely
\begin{align*}
   \left (  \omega_{ \vpsi }  \left ( \omega_{ 0 } , \eps \right ) \phi \right )  \left (  v \right ) & =  \eps  \gamma \left (  \psi\circ V \right ) ^{ -1 }     \int_{  \voxpr\ad  }  \phi( y ) \psi \left (  -b_{ q }  \left (  v, y \right )  \right ) dy \\
\end{align*}
where $ \eps \in  \left \{ 1, -1 \right \} $.
\end{enumerate}

\end{lemma}

\subsection{ Global theta correspondence }
In this section we describe an important invariant distribution called the \bbb{theta distribution}. Using this we define the theta kernel which will be essential in transferring automorphic forms from $ H\ad $ to $ G\ad $.
\par
\label{notation:AutomorphicForms}
Let $  \cA \left ( G_{ \fA }  \right )  $ be the space of automorphic forms on the group $ G $.
Let 
 $ \fA $
be the \adele ring of a number field  $ \F $ and fix a nontrivial additive character 
 $ \psi = \otimes _{ \nu } \psi_{ \nu }  \,\colon \fA / \F \to \fC $. 
 Let  $ V $ and  $ W $ be nondegenerate symmetric and symplectic quadratic spaces, respectively.
Let  $ W = X\oplus Y $ be a complete polarization of the symplectic vector space $ W $ over  $ \F $ namely  $ X, Y $ are maximal isotropic subspaces and  $ X, Y $ are dual to each other via the symplectic form  $ \left < \cdot, \cdot \right > $.
Consider the Schwartz space
 $ S( X\ad  ) = \otimes_{ \nu } S( X_{ \nu }  ) $.
We can construct a projective representation $ \rep{ Sp( W\ad  ) }{ \omega_{ \psi }  }{ S( X\ad  )  }  $
 which becomes an ordinary representation of the metaplectic group  $ Mp( W\ad  ) $. 
\begin{defi}[Theta distribution]
The linear functional 
$$
\Theta \,\colon S( X\ad  )  \to \fC
$$
\label{fm:theta functional}given by 
$$
 \Theta( \phi ) = 
 \sum_{ \gamma \in X_{ \F }   } 
  \phi( \gamma )  \, \, \, \, \, \, \, , \phi\in S( X\ad  ) 
$$
is referred to as the \bbb{theta distribution}.
\end{defi}
We have the following fundamental result of A. Weil.
\begin{theorem}[Invariance of the theta distribution]
The theta distribution is  $ Sp( X_{ \F }  )  $-invariant,  i.e.
$$
\Theta \left (  \omega_{ \psi }  \left (  g \right )  \right ) \left (  \phi \right ) =  \Theta( \phi )
$$
  for every $ g\in
Sp \left (  V\otimes W \right )  _{ \F }  
 $
 and every
  $ \phi\in S( X\ad  )  $.

\end{theorem}

\begin{remark}
The above series converges absolutely,  since  $ \phi $ is rapidly decreasing.
\end{remark}
\begin{defi}[Theta kernel]
\begin{enumerate}
\item
We can form a natural equivariant function
$$
\theta_{ \vpsi } \colon \omega_{ \vpsi } \to
\cA \left (  G\ad \times H\ad \right ) 
$$
which we will refer to as the theta kernel. 
\item
Moreover this function can be given explicitly by
$$
\theta_{ \vpsi } ^{ \phi }  \left (  g,h \right ) = 
\sum_{ \gamma \in  \left (  V \otimes X \right )_{ \F }  } 
 \left (  
\omega_{ \vpsi } \left (  g,h \right )
\phi
  \right ) 
 \left (  \gamma \right ) 
$$

\end{enumerate} 
\end{defi}
\subsection{ Transferring automorphic forms 
from $ \cA( H\ad  )  $ to $ \cA( G\ad  ) $
and global theta lifts
}
\label{sec:TransferringAutomorphicForms}
One of the main goals of this work is to lift a very simple automorphic representation of the orthogonal group to a representation of the symplectic or metaplectic group. This new representation will be called the \bbb{global theta lift} and it is the image of an equivariant map.
\par

We may use the theta kernel  
to transfer automorphic forms from $ \cA( H\ad  )  $ to automorphic forms in  $ \cA( G\ad  ) $. 

\begin{defi}[Theta integral]
\label{sec:ThetaMap}
The \bbb{theta integral} is a function
defined as follows. For every  $ f^{ H } \in \pi \subseteq  \cA \left ( H\ad  \right )  $ let 
$$
I \left ( \phi, f^{ H }  \right ) \left (  g \right ) = 
\int_{ H_{ \F } \backslash H_{ \fA }  } 
\theta_{ \psi } ^{ \phi }  \left (  g,h \right )\cdot  f^{ H }  \left ( h \right ) \,dh
$$
where  $ dh $ denotes the Tamagawa measure.
We call the function $ I \left ( \phi, f^{ H }  \right ) $ the theta integral lifting the cusp form $ f^{ H }  $.\label{notation:genThetaIntegral}
\end{defi}
 If the integral in Definition \ref{sec:ThetaMap} is well-defined (for example if  $ \pi \subseteq  \cA \left ( H\ad  \right )  $ is $  \left (  q,V \right )  $ is anisotropic) then we may define a global theta lift.
\begin{defi}[Global Theta Lift]
The image of the map
$$
\Theta \,\colon 
\omega_{ \vpsi } \otimes \pi
 \to  \cA \left ( G\ad \right ) 
$$
denoted by  $  \Theta_{ V, W, \psi } \left ( \pi \right )  $ is called the \bbb{global theta lift of  $ \pi $}. Namely this is the space generated by 
$$
\Set*{ I \left ( \phi, f \right ) } { \phi\in \omega_{ \vpsi } , f \in \pi }
$$
\end{defi}

\begin{remark}
 In this work we will consider reductive dual pairs  $  H\times G  $ and $ H\times G\p $ where  $ H, H\p $ are orthogonal groups and
  $ G $ is a symplectic group (or metaplectic group).  We will be working with characters  $ \xi$ and $\xip $  therefore we will denote the theta integrals in our context by  $ \thlift {   } $ and $ \thliftp {   } $.
So in the global setting we will mainly consider integrals of the form
$$
\thlift {} \left (  g \right ) 
= 
\int_{   H _{ \F }  \backslash  H \ad  }
\theta_{ \vpsi } ^{ \phi }  \left (  g,h \right ) 
 \xi( h ) \,dh
$$
and
$$
 \,
\thliftp {} \left (  g \right ) 
= 
\int_{   H\p _{ \F } \backslash  H\p \ad }
\theta_{ \vpsip } ^{ \phip }  \left (  g,h \right ) 
 \xi\p( h ) \,dh
$$

\end{remark}

\label{fm:thetakernel}

\subsection{ Rallis' theory of towers and first occurrence }
\label{sec:RallisTheoryOfTowers}
This section gives us conditions to understand when a global theta lift is cuspidal which is closely tied with when a theta lift is nonzero in a "tower" of lifts.
This question was first addressed by Rallis in \cite{rallis1984howe}. In order to state this result we must consider a sequence of global theta lifts instead of just one global theta lift.
Usually this is referred to as a tower of theta lifts.
We will consider the tower of theta lifts of the following form

\vsp
\boo{ 2n  }
\vsp
\boo{ 4  }
\boo{ 2  }

where $ \xi \,\colon H\ad \to \fC $ is a nontrivial, quadratic, automorphic character of the orthogonal group $  H\ad  $ and
$  \left (  W_{ k } ,\left < \cdot, \cdot \right > \right )  $ is a symplectic space of dimension $ k $ and $ G_{ k }  $ is the corresponding symplectic or metaplectic group and $ H=  O \left ( q, V \right ) $ for a quadratic space $  \left (  q, V \right )  $.
\par
 We recall the following result of Rallis in our context, namely in the case where the quadratic spaces are globally anisotropic, i.e. the Witt index of the quadratic space is zero.
\begin{theorem}[First occurrence and cuspidality]
\label{sec:FO}
Let  $  \sigma  $ be an irreducible automorphic cuspidal representation of the orthogonal group $  H\ad  $. The smallest $  k \in \fN $ such that
$$
\lift k
\neq 0
$$
satisfies that the representation $ \lift k $ is cuspidal. Moreover
$$
\lift{ 2n } \neq 0
$$
\end{theorem}
\proof
See Rallis \cite{rallis1984howe}.

$ \square $
\begin{remark}
\begin{enumerate}
\item
The smallest integer $  k \in \fN $ such that
$$
\lift k
\neq 0
$$
is called the \bbb{first occurrence index} of $ \sigma $ in the tower of theta lifts.
\item
Although Theorem \ref{sec:FO} is stated for cusp forms, it remains valid for liftings of non-trivial automorphic characters of the orthogonal group.
\item
We will use this result in Section \ref{sec:proofOfCuspidality} in order to prove that the automorphic representations considered in this work are cuspidal.
\end{enumerate}
\end{remark}

\newpage
\newcommand{ \pin }{  \tn{Pin}_V \left ( k \right ) } 
\newcommand{ \sqrK }{ \K^{ \times } 
/
\K^{ \times, 2 } 
 } 
 \newcommand{ \hb }[1]{ Hom \left ( #1, \fC^{ \times }  \right ) }
 \newcommand{ \ovk }{ O \left ( V_\K \right ) } 

\section{ Quadratic characters and quadratic forms } 

In the following we describe notations for quadratic characters of the orthogonal group both globally and locally and we set notations for the quadratic spaces and characters appearing in this work. We recall the definitions and notations of the Hilbert symbol, Hasse invariant and Weil index and we recall basic facts about orthogonal groups and their structure in terms of reflections.
We also recall definitions and set notations of invariants of quadratic spaces over a local field.
\par
\label{section:characters} 
In this section $ \F $ denotes a local field. We essentially state that local quadratic characters  $ \xi $ of  $ H_{ \nu }    $ are given by two parameters  $ \lam \in \sqrclass $ and  $ \eps \in  \left \{ \pm1  \right \} $. Automorphic quadratic characters of  $ H_{ \fA }  $ are given by a product of the local characters where  $  \left (  q  , V  \right )  $ denotes an odd $ m $-dimensional quadratic space.
\begin{remark}
We will denote the local characters depending on the parameters
$ \lam \in \sqrclass $ and  $ \eps \in  \left \{ \pm1  \right \} $
by 
$$
\xi_{ \nu } \defeq \xi_{ \lam, \eps, \nu } 
$$
However, in order to make the notation less tedious at times we will drop the  $ \lam $ and  $ \eps $.

\end{remark}
We will describe automorphic quadratic characters
$$ \xi \,\colon H_{ \fA }  \to \fC^{ \times }   $$
of the odd orthogonal group.
\subsubsection{ Local characters of  $ O( q )( \Fnu ) $ }
We define local characters of the special  orthogonal group $ O( q )( \Fnu ) $ as the composition of the Hilbert symbol $  \alpha _{ \lam, \nu } =  \left (  \cdot, \lam \right )_{ \nu }   $ with the spinor norm $ SN $.
\[
\xymatrix{
O( q )( \Fnu )    \ar[r]^-{SN} \ar@/_2pc/[rr]_{ \chil\circ SN }  \ar[r] &   \sqrclass \ar[r]^-{ \chil }  & \mu_{ 2 } 
} 
\]

Namely the character is given by

$$ \begin{cases}
   \xil \,\colon O( q )( \Fnu )  \to \fC^{ \times } \\
   \xil( h ) = ( SN(h) , \lambda )_{ \F_\nu } \\
\end{cases} $$

i.e. it is the composition of the character 
\label{fm:quad char}
$$ \begin{cases}
   \chil \,\colon \sqrclass \to \mu_{ 2 } \\
   \chil( x ) = ( x , \lambda )_{ \F_\nu } \\
\end{cases} $$
\label{fm:quadchar}
given by the Hilbert symbol
where $ \lambda\in\sqrclass $
and the spinor norm
$$ \begin{cases}
   SN \,\colon O( q )( \Fnu )  \to \sqrclass\\
   SN( h ) = q  ( u_1 ) \cdots q  ( u_n ) \\	
\end{cases} $$

where $$ h = \refl{ u_{ 1 } }\cdots \refl{u_{ n } }  $$
where $ \refl { u_{ k }   } $ are reflections for all  $ 1 \leq k \leq m $.
\begin{remark}
By the Cartan-Dieudonne theorem every element of the orthogonal group can be written as a product of at most $m=\dim V$ reflections. See theorem \ref{theorem:Cartan}.
\end{remark}
 Let us denote $  H = O \left ( q \right )  $.
Since we are assuming $ m $ is odd we have $$ H \cong O( q  ) \times \mu_{ 2 }  $$
therefore we can extend the character $ \xil $ to $ H  $ by defining it on $ \mu_{ 2 }  $ so we denote this extended character by
$$  
\xileps \,\colon H ( \Fnu )  \to \fC^{ \times }  $$
\label{fm:fullchar}where
$$ \restr { \xileps } { SO(q )( \Fnu )  } = \xil $$
and
$$ \eps = \xileps( -I_m )  $$
 
\subsubsection{ Automorphic characters of the orthogonal group $  H_{ \fA } $ }

We will describe nontrivial automorphic quadratic characters
$$ \xi \,\colon H_{ \fA }  \to \fC^{ \times }   $$
of the orthogonal group.
 Such a character is given by the product of local characters as given in the previous section, namely
 $$
\xi = \otimes_{ \nu }  \xi_{ \nu } 
$$
 That this character is automorphic according to the product formula
 $$
\prod_{ \nu } \kh a b = 1
$$
See Lemma \ref{sec:Hilbert} for general properties of the Hilbert symbol.


\good

 \label{section:hilberthasse} 

\subsection{ The Hilbert Symbol, Hasse Invariant and Weil Index }
\good
Let  $ \F $ be a local field. For proofs and generalizations see \cite{fesenko1993local}, \cite{neukirch1986class}, and \cite{serre2012course}.
\label{fm:hilbsymb}
\begin{defi}[Hilbert symbol]\label{def:hilbertsymbol} 
Let  $ a, b\in \F^{ \times }  $. Then we define  $  \left (  a, b \right )_{ \F } = 1  $ 
 if and only if the equation
 $$
z^{ 2 } -ax^{ 2 } -by^{ 2 } = 0
$$
has a non-zero solution  $  \left ( x, y, z \right ) \in\F^{ 3 }  $, otherwise we define  $  \left (  a, b \right )_{ \F } = -1 $. The number  $  \left (  a, b \right )_{ \F } \in \mu_{ 2 }  $ is called the Hilbert symbol of  $ a $ and  $ b $ relative to  $ \F $ where 
 $ \mu_{ 2 } =  \left \{ \pm1 \right \} $.
\end{defi}

We summarize some of the basic properties of the Hilbert symbol in the following
lemma.
\label{sec:Hilbert}
\begin{lemma}
Let  $ \F $ be a local field with  $ a, b, c\in\F^{ \times }  $. Then the Hilbert symbol satisfies the following:
\begin{enumerate}
   \item  $ \kh{ ab }{ c }= \kh{ a }{ c }\kh{ b }{ c } $ \\
   \item  $ \kh a b = \kh b a $ \\
   \item  $ \kh a a =  \kh{ a }{ -1 } $ \\
   \item The Hilbert symbol is non-degenerate, i.e if for a given $ a\in\F^{ \times } , \kh a b = 1 $ for all  $ b\in\F^{ \times }  $,  then  $ b\in \left (  \F^{ \times }  \right )^{ 2 }   $ \\
   \item  $  \left (  a, b \right ) _{ \fC } = 1 $ for all  $ a, b\in\fC^{ \times }  $ \\
   \item  $ \kh{ a }{ -a } = \kh{ a }{ 1 - a } $ where  $ a \in \F^{ \times } , a\neq1 $ \\
   \item Let  $ \F $ be a nonarchimedian local field. Then  $ \kh a b = 1 $,
for all  $ a, b\in\cO_{ \F } ^{ \times }  $ where the characteristic of the residue field  $ \cK_{ \F }  $ is not equal to  $ 2 $. \\
   \item 
If  $ \nu $ ranges over all places,  $ \kh a b = 1 $, for almost all places.
\\
\item
 The Hilbert symbol satisfies the product formula
 $$
\prod_{ \nu } \kh a b = 1
$$
\end{enumerate}

\end{lemma}
\begin{defi}[Hasse invariant]\label{def:hasseinvariant}\label{fm:hasseinv}
Let  $ q \left (  x_{ 1 } , \ldots, x_{ n }  \right ) = a_{ 1 } x_{ 1 } ^{ 2 } + \ldots a_{ n } x_{ n } ^{ 2 }  $ be a quadratic form. The Hasse invariant of  $ q $ is given by
$$
 h_{ \F }  \left (  q \right ) = \Pi_{ i<j } \kh{ a_{ i }  }{ a_{ j }  } \in \mu_{ 2 } 
$$
\end{defi}

\begin{defi}[Weil Index]\label{def:weilindex}\label{fm:weilenden}\label{fm:weilindex}
The Weil index (or Weil factor)  $ \gamma_{ \psi }  \left (  \cdot  \right )  $ is given by
 $$ 
 \gamma_{ \psi }  \left (  a \right )  
 \int_{ \F } 
  f( t ) 
  \dd t
\abs{ a } ^{   \frac{ 1 } { 2 } } 
= 
 \int_{ \F^{ 2 }  } 
  f( t - z ) 
  \psi \left (  \frac{ 1 } { 2 }az^{ 2 }  \right ) 
  \dd z \dd t
 $$
Let  $ \eta $ be a nontrivial continuous character of the additive group of  $ \F $. If  $ a\in\F^{ \times } $,
 the character  $ \eta_{ a }  $ of  $ \F $ is given by
$$
\eta_{ a }  \left (  x \right ) = \eta \left (  ax \right ) , \forall x \in \F
$$

Define:
\label{fm:qweilindex}
$$
\begin{cases}
   \gamma_{ \F }  \left (  \eta \right ) \defeq \tn { Weil index of:  }  x\mapsto\eta \left (  x^{ 2 }  \right )\\
   \gamma_{ \F }  \left (  a, \eta \right ) \defeq\dfrac{ \gamma_{ \F }  \left (  \eta_a \right ) } { \gamma_{ \F }  \left (  \eta \right ) }\\
\end{cases}
$$

\end{defi}
For properties of the Weil index see \cite[p. 176]{weil1964certains} or \cite{ranga1993some}.

\subsection{ Basic facts about quadratic spaces and their invariants }\label{section:quadratic} 
\good
Let  $  \left (  q, V \right )  $ be a quadratic space over a field  $ \F $.


\begin{defi}[Gram matrix]
The Gram matrix associated to a quadratic space $  \left (  q, V \right )  $ of a vector  $ v = \left (  v_{ 1 } , \ldots, v_{ n }  \right )  $ is defined by 
$$
\gram\vv  = \left (  b_{ q }  \left (  v_{ i } , v_{ j }  \right )  \right )_{ i, j } 
$$
where  $ b_{ q }  $ is the bilinear form associated to a quadratic space $  \left (  q, V \right )  $.
\end{defi}

\begin{defi}[Discriminant]
The discriminant of a nondegenerate quadratic space  $  \left (  q, V \right )  $ represented by the matrix $ A $ is the image of  $ \det \left (  A \right )  $ in  $ \F^{ \times } / \F^{ \times, 2 }  $.

\end{defi}


\begin{theorem}
Let  $ \F $ be a local field. Two quadratic spaces  $  \left ( q, V \right ) ,  \left ( q\p , V\p \right ) `$ over  $ \F $ are isometric if and only if 
$$
\begin{cases}
   \dim\,q = \dim\,\qp \\
   d( q ) =  d( \qp ) \\
   h \left ( q \right ) = h \left (  \qp \right )  \\
\end{cases}
$$
where  $ d \left (  \cdot  \right )  $ and  $ h \left (  \cdot  \right )  $ denote the discriminant and Hasse invariant of the form,  respectively.

\end{theorem}
\proof
See \cite[p. 158]{lam2005introduction}.

$ \square $

\subsubsection{ Reflections }
\begin{defi}
Let  $  \left (  q, V \right )  $ be a quadratic space. For an anisotropic vector $ v \in V $ define the hyperplane reflection with respect to  $ v $ by
$$
\label{fm:refl}
\begin{cases}
   \tau_{ v }  \,\colon V \to V \\
   \tau_{ v }  \left (  w \right ) = w - 2 \cdot \dfrac{ b_{ q }  \left (  v, w \right )  } { q \left (  v \right )  }v \\
\end{cases}
$$
\end{defi}
 
\begin{lemma}
The hyperplane reflection satisfies the following
\begin{enumerate}
\item
 $ \tau_{ v } \in  O \left ( q, V \right )  $.
\item 
 $ det \left ( \tau_{ v } \right ) = -1 $ .
   \item  $  \tau_{ v }  \left (  v \right ) = -v $.
   \item  $ \tau_{ v }  \left (  w \right ) = w  $ 
such that   
 $ b_{ q }  \left ( v, w \right ) = 0 $.
\end{enumerate}
\proof
Trivial. $ \square $
\end{lemma}

\begin{theorem}[Cartan-Dieudonne] \label{theorem:Cartan} Let  $  \left (  V,q \right )  $ be a regular quadratic space of dimension  $ n $. Then every  $ h \in O \left (  q \right )  $ is a product of at most  $ n $ hyperplane reflections.

\end{theorem}
\proof
See \cite[p. 18]{lam2005introduction} $ \square $

\subsection{ Characters }
\label{sec:charClassificationAdmQuad}

\subsubsection{ Some character computations. }
In this section we compute present conditions for our non-trivial characters to be trivial on stabilizers of orthogonal groups.
In this section for $ \vv =  \left (  v_{ 1 } , \ldots, v_{ n }  \right )  \in V^{ n }  $ we denote the stabilizer of $ \vv $ by 
$ H^{ \vv }  $. Moreover if $ H =  O \left ( q, V \right )  $ where $ \dim V = m $ then at times we will denote $ \left ( q_{ m } , V_{ m }  \right ) ,H_{ m } ,H_{ m } ^{ \vv }  $ to stress the dimensions we are considering.

\begin{theorem}
\label{theDataDescription}
Suppose $ m = n + 2 $ and $ \vvbeta = \left (  v_{ 1 } ,\ldots , v_{ n }   \right )\in V_{ n+2 }^{ n }   $ is an orthogonal sequence of (linearly independent) vectors such that 
$$ \gramn{ n + 2 }{ \vvbeta } = \beta = \diag \left (  \beta_{ 1 } ,\ldots, \beta_{ n }  \right ) \in \Sym( \Fnu )  $$
and suppose

$$ \lambda \in d( V_{n+2} ) d( \Vp ) \cdot \squarex $$
and
$$ \eps = ( q_{ n+2 } ( u ) d( V_{n+2} )  , \lambda )  , \tn{ for some  }u \perp \vvbeta $$
then
$$ \restr { \xileps } { H_{ n+2 } ^{ \vvbeta }  } \equiv 1$$

\end{theorem}
\proof

Let $ u_{ 1 }, u_{ 2 } \in \left (  \vvbeta  \right )^{ \perp }  $ be orthogonal linearly independent, i.e.
$ \left (  v_{ 1 }, \dots, v_{ n } ,u_{ 1 } ,u_{ 2 }   \right ) $ is an orthogonal basis of $ V_{ n+2 } $. We denote the sub quadratic space spanned by $ \vvbeta $ by $ ( \Vp , q\prime_{ n }  ) $ where $ \Vp = \Span  \left \{ v_{ 1 } , \ldots, v_{ n }   \right \} $ and $ q^{ \prime } _{ n } = \restr { q_{ n+2 }  } { \Vp }$

Let $$ h_{ 0 } = \refl{u_{ 1 } } \refl{u_{ 2 } }\in H_{ m } ^{ \vvbeta }  \cong H( \Span \left \{  u_{ 1 } , u_{ 2 }   \right \}  )  $$ then since $ \refui\refuii \in SO( q_{ n+2 }  )  $ by definition of the character we have

\begin{align*}
   \xileps( \refui\refuii )  & = \xileps( \refvv\refvv\refui\refuii )\\
    & =  \chil( \underbrace{ q_{ m } ( \refl \vi ), \ldots, q_{ m } ( \refl \vn ) }_{ d(\Vp) } \cdot \underbrace{ q_{ m } ( \refl \vi ), \ldots, q_{ m } ( \refl \vn )  q_{ m } ( u_{ 1 }  ) q_{ m } ( u_{ 2 }  ) }_{ d(V_{n+2}) }  )\\
    & = \chil( d(\Vp)\cdot d(V_{ n+2 }) ) \\
    & = (  d(\Vp)\cdot d(V_{ n+2 }) , \lambda ) _{ \nu } \\
    & = (  d(\Vp)\cdot d(V_{ n+2 }) , \lambda ) _{ \nu } \\
    & = (  -\lambda  , \lambda ) _{ \nu } \\
    & = 1\\
\end{align*}
Now consider $ h_{ 0 } = \refu \in H_{ m } ^{ \vvbeta }  $ for $ u \perp v_{ i } , 1 \leq i \leq n $.

Then we have
$$ 
\xil( \refu )  $$
$$ 
\xil( -\refu\underbrace{ \refu\refup\refvv }_{ -id } )  $$
$$ 
\xil( -1 ) ( q_{ n+2 } ( u ) d( V_{n+2} )  , \lambda )  
 $$
 Note for linearly independent $ u_{ 1 } , u_{ 2 }  $ not necessarily orthogonal.
 
 $$
 \xileps( \refui ) 
 \xileps( \refuii ) =  \xileps( \refui \refuii ) = 1 \implies
  \xileps( \refui ) = 
 \xileps( \refuii )
  $$

\begin{align*}
   \xil( -\refvv )  & =  \xil( \refvv\underbrace{ \refvv\refui\refuii }_{ -id } ) \\
    & = \xil( \refui\refuii )  \\
    & = 1\\
\end{align*}

$ \square $

\newpage
\section{ Examples of Global Theta Lifts }
\label{sec:ExampleLiftingCharacters}
\subsection{ Example of lifting a trivial automorphic character }
For example one could lift a trivial automorphic character $ \xi \,\colon H_\fA \to \fC^\times $ (that is $ \xi\equiv 1 $) of the orthogonal group to an automorphic representation of the symplectic (or metaplectic) group. This lift is characterized by what is known as the Siegel-Weil formula and has deep arithmetic significance, generalizes work of Siegel on representation numbers of quadratic forms and has many important applications in the theory of automorphic forms.
\par
In the Siegel-Weil formula one equates a global theta integral with a certain Eisenstein series, so roughly speaking one obtains a formula of functions on $ G_{ \fA } $ of the form 
$$
I^{ \phi }_{ V, \xi \equiv 1 }  = Eis^{ \phi }
\qquad
 \left (  ! \right ) 
$$
where $ \phi $ is a test function, $ E= Eis $ is an Eisenstein series such that the section defining it depends on the test function $ \phi $ and $ I^{ \phi }_{ V, \xi \equiv 1 } $ is a theta integral.
\subsection{ Example of lifting a non-trivial automorphic character }
Another natural example would be to lift a non-trivial automorphic character $ \xi \,\colon H_\fA \to \fC^\times $ of an orthogonal group $ H_{ \fA }  $. For instance $ \xi $ could be taken to be trivial at almost every place outside of a finite set $ S $ and equal to the determinant character for every place $ \nu\in S $, i.e. 
$
\xi =
 \left (  \otimes_{ \nu\in S } det_{ \nu }  \right ) 
\otimes 
 \left (  \otimes_{ \nu \not \in S } 1_{ \nu }  \right ) 
$. Lifts of non-trivial automorphic quadratic characters are the primary concern in this work and this is the question Snitz addressed for certain low rank orthogonal groups of dimensions $ 3 $ and $ 1 $.
\par
It is important to note that a representation of the group $ G_{ \fA }  $ obtained in this manner is a cuspidal representation and it is known that such a representation can not be given by an Eisenstein series, hence a priori one could not expect to obtain a formula of the form $  \left (  ! \right )  $. Instead Snitz considers \underline{another} orthogonal group $ H\p_{ \fA }  $ together with a very specific non-trivial automorphic character $  \xi\p \,\colon H_\fA\p \to \fC^\times $ and lifts it via the global theta correspondence to $ G_{ \fA }  $. Essentially Snitz constructs a morphism between these two global theta lifts given by a morphism 
$ \Delta \,\colon  \schad   \to  S \left ( V_\fA\p \right )  $
between the test function spaces.
\par
In this manner, Snitz obtained an analog of the Siegel-Weil formula, so roughly speaking Snitz obtained a formula of functions on $ G_{ \fA } $ of the form 
$$
I^{ \phi }_{ V, \xi }  = I^{  \Delta \left ( \phi \right )  }_{ V\p, \xi\p }
\qquad
 \left (  !! \right ) 
$$
where $ \phi $ is a test function and $  \Delta \left ( \phi \right )  $ is a test function matching the test function $ \phi $. Henceforth the space of test functions on $ V_{ \fA } $ will be called the space of {\bfseries{Schwartz functions}} $ \schad  $.
\begin{remark}
On the surface no Eisenstein series appears in formula $  \left (  !! \right )  $, however it turns out that the automorphic representations obtained via the global theta lifts in both Snitz's and this work are what are referred to as CAP (\bbb{C}uspidal \bbb{A}ssociated to \bbb{P}arabolic) representations, colloquially referred to as "shadows of Eisenstein series", hence formula $  \left (  !! \right )  $ is similar to $  \left (  ! \right )  $, hence an honest analog of the Siegel-Weil formula.
\par
It is interesting to note that CAP representations provide counterexamples to the Generalized Ramanujan conjecture.
By work of H. Jacquet and J. Shalika \cite{jacquet1981euler}, it can be shown that the CAP phenomenon never occurs for $ GL(n) $. At the same time the Ramanujan conjecture for the general linear group $ GL(n) $ is believed to be true. CAP representations are not really the focus of this work, however the fact that CAP representations naturally appear in our study is a further indication of the arithmetic significance of these results.
For more on CAP representations and the Ramanujan conjecture see appendix \ref{sec:CAPRamanujan}.
\end{remark}

\newpage
\section{  CAP representations and the Generalized Ramanujan conjecture }
\label{sec:CAPRamanujan}

This exposition closely follows Eitan Sayag’s M.A. thesis \cite{sayag} and also Wee Teck Gan's paper on Saito–Kurokawa representations of $ PGSp_{ 4 }  $\cite{gan2008saito}.
The notion of a CAP representation was first coined by \piat \cite{piatetski1983saito}. In this paper he gave a characterization of CAP representations for the group $ PGSp_{ 4 }  $  in terms of global theta lifts for the Metaplectic group $ \widetilde{ SL_2 } $.
\par
What exactly is a CAP representation? We can think of such a representation as a cuspidal representation which is pretending to be an Eisenstein series. Hence it is also call a shadow of an Eisenstein series. We would like to present a more formal definition. First note that for an algebraic group $ G $ over a global field $ \F $,  two automorphic representations of the group $ G $  are called equivalence if at almost every place their local components are isomorphic .
\begin{defi}[CAP representation]
Given an automorphic, cuspidal representation  of the group $ G\ad $  we call such a representation CAP is there exist a parabolic subgroup $ P $  with Levi decomposition $ P = MU $ of $ G $  and an automorphic, cuspidal representation representation $ \sigma $  of the Levi part $ M\ad $ such that the representation $ \pi $  is  equivalent to a constituent of $ Ind_{ P\ad } ^{ G\ad }
\sigma
  $.
\end{defi}
Do CAP for presentations exist? It's impossible to embed a cuspidal representation in an induced representation. Therefore CAP representations are considered to be an anomaly in representation theory. Indeed, CAP representations are quite rare since according to strong multiplicity one we can conclude that the general linear group $ GL_{ n }  $ does not have any CAP representations.
\par
One of the most interesting achievements appearing in \cite{piatetski1983saito} was to show the connection between counter-examples of the generalized Ramanujan conjecture (GRC) and between CAP representations. Piatetski noticed that all of the counter-examples to GRC were CAP representations and this motivated research on CAP representations for every algebra group $ G $. This gave a “ classification” of the counter-examples of the generalized Ramanujan conjecture.
\par
\subsection{ Generalized Ramanujan-Peterson conjecture }
We would like to recall the Ramanujan-Peterson conjecture.
Let us quickly recall the setup for the classical theory of modular forms.
 The modular group $ \Gamma= SL_{ 2 }  \left (  \fZ \right )  $  acts holomorphically on the upper half plane
$ \cH= 
\Set*{ z\in\fC } {  Im \left ( z \right ) >0 }
 $
 with respect to the action
$$
\gamma
\cdot z
= 
\dfrac{ az+b } { cz+d }
$$
where $ z\in \cH, \gamma= \abcd\in \Gamma $.

\begin{defi}
A modular form of weight $ k $  with respect to $ \Gamma $ is a holomorphic function $ f \,\colon \cH \to \fC $ which satisfies
\begin{enumerate}
\item
$ f \left (  \gamma \cdot z \right ) = 
 \left (  cz+ d \right ) ^{ k }  f \left ( z \right ) 
 $
for every $ \gamma \in \Gamma ,z\in \cH $.
\item 
$ f $ Is holomorphic at the cusp of $ \Gamma  $, i.e.  at infinity, in other words the Fourier expansion of $ f $  does not have negative powers, i.e. 
$$
f \left ( z \right ) =
\sum_{ n\geq0 } a_{ n } e^{ 2\pi i n z } 
$$ 

\end{enumerate}

\end{defi}
\begin{defi} 
\begin{enumerate}
\item
Let $ k\geq0 $. We denote the space of modular forms of weight $ k $ by 
$
M_{ k }  \left (  \Gamma  \right ) 
$
\item
Next we define the space of cusp forms of weight $ k $ by 
$$
S_{ k }  \left (  \Gamma  \right ) = 
\Set*{ f\in M_{ k }  \left (  \Gamma  \right )  } { a_0 = 0 }
$$
\end{enumerate}
\end{defi}
\begin{theorem}
\begin{enumerate}
\item
The vector space $ M_{ k }  \left (  \Gamma  \right )  $  is finite dimensional and has the structure of a Hilbert space.
\item
There exists a family Hermitian operators $ T_{ n }  $ acting on $ M_{ k }  \left (  \Gamma  \right )  $. These operators are called Hecke operators.
\end{enumerate}

\end{theorem}

\newpage
\label{sec:SnitzsWork}
In this section we briefly recall Snitz's results.
 \par
 Let  $ \F $ be a number field with \adele ring  $ \fA = \fA_{ \F } $ and  $ B $ be a quaternion division algebra over  $ \F $. Let  $ S =  S_{ B }  $ be the set of places in which $ B $ is ramified. Let  $  \left (  V_{ S } , q_{ S }  \right )  $ be the  $ 3 $-dimensional quadratic space where $ V $ is the space of trace zero elements in  $ B $ and $ q_{ V_{ S }  }  $ is the reduced norm of  $ B $. Snitz considers the reductive dual pair  $  \left (  SO \left (  V_{ S }  \right ), \widetilde{ SL }_{ 2 }  \right )   $, with its associated Weil representation realized on the space  $ S( V_{ S } ( \fA ) )  $ of Schwartz functions on  $  V_{ S } ( \fA )  $. Hence one may consider the theta lift of  $ \xi_{ \kappa }  $ where  $ \xi_{ \kappa }  $ is a quadratic character of  $ \fA^{ \times } / \F^{ \times }  $ which is given locally by the product of local characters given by the composition with the Hilbert symbol with the reduced norm (more precisely the spinor norm of  $  SO( V_{ S }  )  $). Namely, if  $ \phi  \in SO( V_{ S } ( \fA ) )  $, set 
 
 $$
I_{ V_{ S }  } ^{ \phi }  \left (  \xi_{ \kappa } 
 \right )   \in Fn \left ( \slsl \left (  \F \right ) \backslash \slsl \left (  \fA \right )   \right ) 
$$

given by 
$$
I_{ V_{ S }  } ^{ \phi }  \left (  \xi_{ \kappa }  \right ) \left (  g \right ) = 
\int_{   SO( V_{ S }  )( \F )  \backslash  SO( V_{ S }  )( \fA )  } 
\theta_{ V_{ S } , \psi } ^{ \phi }  \left (  h, g \right ) 
\cdot
 \xi_{ \kappa } ( h ) 
 \,dh
$$
and the space that these lifts generate.
\par
For trivial  $ \xi_{ \kappa }  = 1 $, the Siegel-Weil formula associates to each  $  \phi  \in SO( V( \fA ) ) $ a certain Eisenstein series  $  E_{ V } ( f )  $ such that 
$$
I_{ V } ^{ \phi }  \left (  \xi \right ) = E_{ V } ( f )
$$
This gives an alternative construction of  $ I_{ V } ^{ \phi }  \left (  \xi \right )  $. In \cite{snitz}, Snitz considered the analogous problem when  $ \xi $ is a non-trivial character. In this setting,
the function  $ I_{ V } ^{ \phi }  \left (  \xi \right )  $ is contained in the space of cusp forms and one would like an alternative construction of this cusp form. Moreover, we would like to generalize this procedure for higher rank orthogonal groups. We focus on lifting of nontrivial automorphic characters of the odd orthogonal group. The even case should be similar.
\par
In order to do this Snitz constructs a certain one-dimensional quadratic space  $  \left ( q_{ \kappa } , V_{ \kappa }   \right )  $ and a certain nontrivial automorphic quadratic character 
$$ 
\xi_{ S }  \,\colon O \left (  q_{ \kappa }  \right ) \left (  \fA \right )   \to \fC^\times
 $$
Snitz then considers the reductive dual pair  $  \left (  O \left (  q_{ \kappa }  \right ) , \slsl \right )  $
 with its associated Weil representation realized on  $ S( V_\kappa( \fA ) )  $. Thus for each Schwartz function  $ \phi^{ \prime }  \in S( V_\kappa( \fA ) )  $, one has the analogous theta lift
$$
I_{ V_{ \kappa }  } ^{ \phi }  \left (  \xi_{ S }  \right )\left (  g \right ) = 
\int_{   O( V_{ \kappa }  )( \F )  \backslash  O( V_{ \kappa }  )( \fA )  } 
\theta_{ V_{ \kappa } , \psi } ^{ \phi }  \left (  h, g \right ) 
\cdot
 \xi_{ S } ( h ) 
 \,dh
$$
and the space that these lifts generate.
Now the main result of \cite{snitz} is:

\begin{theorem}
There exists a factorizable morphism  $ \Delta \,\colon S( V_S( \fA ) ) \to S( V_\kappa( \fA ) )  $ such that 
the following diagram commutes 
\[
\xymatrix{
\sva \ar[rr]^{\Delta} \ar[dr]^{\theta}&      & \svoa \ar[dl]_{\theta\p}\\
	     & \arithFunc   &		\\
} 
\]
that is
$$
I_{ V_{ S }  } ^{ \phi }  \left (  \xi_{ \kappa }  \right )
= 
I_{ V_{ \kappa }  } ^{  \Delta( \phi )  } \left (  \xi_{ S }  \right )
$$
Moreover the morphism  $ \Delta = \otimes_{ \nu }  \Delta_\nu
$ is factorizable and for each place  $ \nu $ we have local morphisms  $ \Delta_{ \nu }  \,\colon S( V_S( \F_{ \nu }  ) ) \to S( V_\kappa( \F_{ \nu }  ) ) $ given by the following orbital integrals 
$$
\Delta_{ \nu }  \left (  \phi \right )  \left (  x \right ) = 
\abs{ x } _{ \nu }  \left (  x, -\kappa \right ) _{ \nu } 
\int_{  SO( V ) ^{ \vv }  \left (  \F_{ \nu }  \right ) \backslash  SO( V ) \left (  \F_{ \nu }  \right ) }
 \phi_{ \nu }  \left (  xh^{ -1 } \cdot \vv \right )\xi_{ S }  \left (  h \right ) 
\,dh
$$
where  $ \vv \in V_{ \kappa }  \left (  \F \right )  $ is a fixed choice of a vector such that
 $ q_{ V_S }  \left (  \vv \right ) = d \left (  q_{ \kappa }  \right )  $, that is, we have an embedding of quadratic spaces  $ V_{ \kappa } \hookrightarrow V_{ S }  $ where 
  $ V_{ \kappa } = \Span \left 
 \{ \vv  \right \}
 $.

\end{theorem}s

Let  
\begin{enumerate}
   \item $ \sigma \cong \otimes \sigma_{ \nu } = \Theta_{ V_S }  \left (  \xi_{ \kappa }  \right )  $  \\
   \item $ \sigma\p \cong \otimes \sigma\p_{  \nu }  = \Theta_{ V_\kappa }  \left (  \xi_{ S }  \right )  $  \\
\end{enumerate}
be the global theta lifts  $ \theta $ and  $ \theta\p  $. Then we have the following corollary
\begin{cor}
 $ \sigma \cong \sigmap $
\end{cor}

\newpage
\section{ Local pairing }

In this section we present the actual explicit formulas for the local pairing of interest.
One must keep in mind that the actual formula for this local pairing is of marginal significance. The important thing to understand is that this pairing is equivariant with respect to the diagonal action of $ G $ and the actions of $ H $ and $ H\p $ and it is a pairing between two functions which are elements of $ \pi_{ \nu } $ and $ \pi\p_{ \nu } $.
Moreover these local pairings will depend on a choice of certain vectors $ \vv $ and $ \vv\p $. This is an important point to keep in mind.
\par
Before we define the local pairing we must set some notations.
Let
$ P  $ be the Siegel parabolic of $ G $ and $ M $ be the Levi subgroup in the Siegel parabolic $ P $. The Levi subgroup $ M \subseteq P \subseteq G $ of the Siegel parabolic $ P $ is given by
$$
M = \Set*{  m \left ( a \right )  = \begin{pmatrix} a & 0 \\ 0 & a^\star \end{pmatrix} } { a \in GL_{ n } }
$$

Fix a maximal compact subgroup $ K \subseteq G\ad $  so that we have an Iwasawa decomposition 
$$
G\ad = P\ad K
 = M\ad N\ad K
$$
 and for 
$$
h = n\cdot m \left ( a \right ) k
,\qquad
n \in N\ad,\,\,a \in GL_{ n ,\fA} ,\,\,k \in K
$$
Let 
$$
\abs{  a \left ( g \right )  }  = 
\abs{  \det \left ( a \right )  } 
$$ 
and for
$$
\rhon = \dfrac{ n+1 } { 2 }
$$
let the section
$$
\Psi_{ s } \,\colon G\ad \to \fC
$$
be given by
$$
\Psi_{ s } \left (  g \right )  = 
\abs{  a \left ( g \right )  } ^{ s + \rhon } 
$$
Likewise, locally we can define the character
$$
\Psi_{ s,\nu } \left (  g \right )  = 
\abs{  a \left ( g \right )  }_{ \nu }  ^{ s + \rhon } 
$$
\begin{remark}
Usually we will not write the place $ \nu $ since it should be clear from the context whether we are working locally or globally.
Thus when working locally we will continue to denote $ \Psi_{ s } $ instead of $ \Psi_{ s,\nu } $.
\end{remark}

\par
Before we present the definition of our local pairing we must present a few more definitions.
We make certain choices of vectors $ \vv $ and $ \vv\p $ which our local pairing will depend upon. Namely let $ \vv \in  \left (  \vox \right )   _{ \F }  \cong V_{ \F } ^{ n }
\cong
\homj
,
\,
\vv\p \in  \left (  \voxp \right )   _{ \F }  \cong V_{ \F } ^{\prime, n }
\cong
\homjp
 $
 and
 $$
\gram\vv  =  \gramp{ \vv\p  }
$$
i.e. the Gram matrices of $ \vv $ and $ \vv\p $ are equal. 
\par
 We will consider vectors $ W_{ \phi } $ and $ W\p_{ \phip } $ coming from the local theta correspondence.
These vectors are given by
$$
W_{ \phi }  \left ( g \right ) = 
\HintegralF
\left (
\wfunc  \left ( g, h \right ) 
\phi \right )
\left (  
\vv
\right ) 
\xi \left ( h \right ) 
\,dh
$$
and
$$
W\p_{ \phip }  \left ( g \right ) = 
\HintegralFp
\left (
\wfuncpi \left ( g, h \p  \right ) 
\phip \right )
\left (  
\vv\p
\right ) 
\xi\p \left ( h \p \right ) 
\,dh\p$$

\begin{defi}[The local pairing]
\label{sec:defLocalPairingTopii}
\label{notation:localNaturalPairingTopii}
We define a local pairing
$$
B_{ \nu,s } \,\colon \schloc  \times \schlocp   \to \fC 
$$
given by
\begin{align*}
B_{ \nu,s } \left ( \phi, \phip \right )  & =  
 \int_{ K_{ \nu }  } \int_{ M_{ \nu }  }  W_{ \phi }  \left (  mk \right )  W\p_{ \phip }  \left (  mk \right ) \delb m \fPsi{ mk }{ s } \,dm \,dk
\end{align*}

$$
\localPairingSigs
$$

\end{defi}

\bibliographystyle{plain} 

\vspace{0.5cm}
\noindent
{\scshape{School of Mathematics, Tel Aviv university, Tel Aviv 69978, Israel}}
\\
\textit{E-mail address:} erezron@tauex.tau.ac.il

\end{document}